\documentclass[11pt]{article}
\usepackage[margin=1in]{geometry}
\usepackage{amsmath,amssymb,amsthm,mathtools,mathrsfs}
\usepackage{booktabs,array,multirow,graphicx,enumitem}
\usepackage[authoryear,round]{natbib}
\usepackage{hyperref}
\hypersetup{colorlinks=true,linkcolor=blue,citecolor=blue,urlcolor=blue}
\newtheorem{theorem}{Theorem}[section]
\newtheorem{proposition}[theorem]{Proposition}
\newtheorem{lemma}[theorem]{Lemma}
\newtheorem{corollary}[theorem]{Corollary}
\theoremstyle{definition}
\newtheorem{definition}{Definition}[section]

\newtheorem{remark}[definition]{Remark}
\newcommand{\Pp}{\mathbb P}\newcommand{\E}{\mathbb E}
\newcommand{\1}{\mathbf 1}\newcommand{\op}{\operatorname}
\newcommand{\logit}{\operatorname{logit}}
\newcommand{\cI}{\mathcal I}\newcommand{\cJ}{\mathcal J}
\newcommand{\btheta}{\boldsymbol\theta}
\newcommand{\cQ}{\mathcal Q}
\newcommand{\cA}{\mathcal A}\newcommand{\suppsec}[2]{Appendix~\ref{#2}}

\title{Multiscale Localized Inference for Networks\\with a Measured Vertex Coordinate}
\author{Marios Papamichalis \and Regina Ruane}
\date{Preprint}
\begin{document}
\maketitle
\begin{abstract}
In many networks each vertex has a position measured from outside the network:
a neuron's location along the body axis, a residue's index along a protein
sequence, a genomic bin's position in base pairs.  The chance that two
vertices connect is then a surface over pairs of positions, and questions
about the network become questions about small regions of that surface.  A
region on the diagonal covers pairs inside one stretch of the axis, and a
departure there is a community with a boundary.  A region off the diagonal
covers pairs spanning two separated stretches, and a departure there is a
bridge.  Existing methods address one part of this at a time: community
detection returns groups without placing them on the axis, block models fix
the width in advance, and scan statistics test one window at one scale.  We
expand the surface in a wavelet dictionary whose elements are exactly such
regions, at every position and width.  We derive in closed form the
modularity, edge length, transitivity and degree spread each element
produces, and the generating element is recovered from those summaries.  A
scan estimates every coefficient against a background fitted on separate
edges and controls the error rate across all positions and widths at once.
Applied to a connectome, a protein contact map and a chromatin contact map, it
recovers known anatomy in the first and is calibrated in all three.  When positions are estimated with error near
the width sought, the location of a departure is not identified at any signal
strength, and a check decides this before any analysis.
\end{abstract}

\section{Introduction}                                          \label{sec:intro}

Many networks come with a measured coordinate for each vertex.  Neurons in a
connectome have a position along the body axis, recorded in an anatomical
atlas.  Residues in a protein have an index along the sequence.  Bins of a
chromosome have a position in base pairs.  In each case the coordinate is
external to the network and exact, and the questions asked of the network are
questions about places on it.  Does this stretch of the sequence fold into a
unit?  Do these two distant stretches contact each other?  Does the connection
pattern change at this point and not at that one?

Each such question has three parts: what kind of structure is present, where
along the coordinate it sits, and how wide it is.  Existing tools answer at
most one part at a time.  Community detection finds groups but does not place
them on an axis, so its answer is a partition of vertices and not a location.
A stochastic block model fixes one resolution in advance, so the width is
assumed and not estimated.  A scan statistic tests one window shape at one
scale.  None of them reports which kind of structure is present, and none
controls the error across the many locations and widths that are examined.

The chance that two vertices connect depends on where each sits on the axis,
through a surface $p(u,v)$ over pairs of positions.  Most of that surface is
background: vertices near each other connect more, and some regions are
busier than others.  A localized departure from the background is one
confined to a single stretch of the axis, or to a single pair of stretches, at
a single width, together with a rule for which pairs inside it connect more
and which fewer.  We call such a departure an atom, and the collection of
atoms at every place and width the dictionary.  Each atom is a product of two
wavelets on the coordinate, and each generates a distinct kind of network.  A diagonal atom raises the connection probability inside a
stretch and lowers it across the middle, which produces two adjacent
communities with a boundary between them.  An off-diagonal atom links two
separated stretches, which produces a bridge.  Atoms at several scales
together produce a hierarchy.  Terms that are not wavelets generate the things
that are not localized structure: additive terms generate hubs and a function
of coordinate separation generates a smooth gradient, and both belong to the
background.  The inferential question is which atoms a given network contains,
at which positions and widths, and with what confidence.

The paper makes four contributions.  First, we prove which network summary
each atom generates.  Modularity, edge length, transitivity and degree
heterogeneity are each an explicit functional of the atom, diagonal and
off-diagonal atoms have opposite signs in the first two, and the shape of the
wavelet enters only through its fourth moment (Theorem~\ref{thm:map}).  The
generating atom is recovered from these summaries in every one of $100$
simulated networks, and on five real networks the fitted model reproduces
the observed modularity and edge length within $0.05$
(Section~\ref{sec:gencheck}).  Second, we give a scan that estimates every coefficient
of the dictionary against a background fitted on an independent set of edges
and controls the familywise error over all atoms at once, and we prove, for a
regular design with an exact coordinate, that it attains the boundary at which
a localized departure becomes detectable, with every rate derived from
primitive conditions (Theorem~\ref{thm:uncond}).  Third, we show that the
canonical generative mechanisms carry no localized information: constant,
additive, separation and coarse-block surfaces all have dictionary
coefficients equal to zero up to a term of order $h^{d+1}$
(Proposition~\ref{prop:dict}), so the question the scan answers cannot be
posed within those models.  Fourth, we delimit the method.  Resolving a
feature of relative width $h$ in a coordinate of dimension $d$ requires
$n\bar d\ge30\cdot4^{d}h^{-2d}$, with $n$ the number of vertices and $\bar d$
the mean degree (Theorem~\ref{thm:info}), and when the
coordinate is estimated with error comparable to $h$ the location of a feature
is not identified at any amplitude although its presence is
(Proposition~\ref{prop:obstruct}).

Three applications use a physically measured coordinate.  On seven
\emph{C.~elegans} connectomes that meet the information requirement, out of
eight reconstructed, the boundary statistic tracks the ganglion partition of
the body axis in every one, with a combined $p$-value below $10^{-5}$.  On thirty protein chains the bridge atoms contain annotated beta-sheet
pairings at $1.79$ times the base rate, with exact $p=7.5\times10^{-7}$.  On a chromatin contact map
the scan is calibrated to nominal level and returns nested boundaries with no
bridges.  Two networks with estimated coordinates, a friendship network and a
legislative cosponsorship network, are refused by the feasibility check and
confirm the refusal when run.  The restriction to measured coordinates is the
content of the fourth contribution, and the applications are chosen to
respect it.

The paper is organized as follows.  Section~\ref{sec:model} defines the model,
the measured coordinate and the dictionary of atoms.  Section~\ref{sec:method}
gives the estimator and its calibration.  Section~\ref{sec:results} contains
the theory: the properties of the dictionary as a basis, the network summaries
each atom generates, the familywise error control and detection boundary at an
exactly measured coordinate, and the two limits on coordinate dimension and
coordinate accuracy.  Section~\ref{sec:sim} reports the simulation study, with
one subsection for each of these results.  Section~\ref{sec:apps} applies the
procedure to a connectome, thirty protein contact maps and a chromatin contact
map, and compares it against the domain callers of those fields.
Section~\ref{sec:related} states the relation to multiscale scanning, to graph
wavelets and to the detection and localization literature.
Section~\ref{sec:discussion} states what is established, what is assumed and
what remains open.  Proofs and protocols are in the Appendix.

\section{Model and Estimand}                              \label{sec:model}

\subsection{Sparse logistic graphon with a canonical coordinate}

Let $U_i\in[0,1]$ denote the target coordinate of vertex $i$.  The coordinate is external and is not a seriation of $A$.  It may be directly measured, obtained as the rank of an external scalar trait, or estimated from data independent of the inference dyads.  We orient the coordinate by a prespecified anchor.  If the external measurement identifies only an unoriented axis, the target is the two-element equivalence class $[U]=\{U,1-U\}$ until that anchor is imposed.

Conditional on $U=(U_1,\ldots,U_n)$, observed dyads satisfy
\begin{equation}
 A_{ij}\mid U,M_{ij}=1\sim\op{Bernoulli}(p_{ij}),\qquad
 p_{ij}=\sigma\{\alpha_n+f(U_i,U_j)\},                       \label{eq:model}
\end{equation}
independently for $i<j$.  The observation indicator $M_{ij}$ is conditionally
independent of $A_{ij}$ given design variables and has propensity bounded below on
the target region.  The logit surface $f$ is symmetric, $\|f\|_\infty\le B$, and
satisfies $\int_{[0,1]^2}f(u,v)\,du\,dv=0$, so its Lebesgue-constant component is
absorbed by $\alpha_n$.  Writing $\rho_n=\sigma(\alpha_n)$, bounded $f$ implies
$p_{ij}\asymp\rho_n$ in the sparse regime.

The theory is conditional on the realized coordinate design and then averaged for unconditional statements.  Random coordinates may have a nonuniform density bounded above and away from zero.  Design correction targets Lebesgue location and scale instead of the frequency with which the observed vertices populate a region.

\subsection{Background projection and localized alternatives}

Let $V_{J_0}$ be a fixed symmetric coarse tensor-wavelet space.  Define
\begin{equation}
  f_0=\Pi_{V_{J_0}}f,\qquad \delta=(I-\Pi_{V_{J_0}})f,          \label{eq:projection}
\end{equation}
where $\Pi_{V_{J_0}}$ is the Lebesgue $L^2([0,1]^2)$ projection.  This convention
defines the background relative to the identified coordinate: without it, a smooth
term can be transferred arbitrarily between $f_0$ and $\delta$.  The surface $f$
is unknown and may vary over a H\"older or Besov ball.  Its coarse part $f_0$ is by
this definition piecewise constant on the blocks of $V_{J_0}$ and contains no
component in the tested detail space, and everything smooth that $f$ carries
beyond $V_{J_0}$ is by the same definition part of the departure $\delta$.
Section~\ref{sec:nuisance} takes up what this means for a background that
varies smoothly along the coordinate.

Let $\psi_{j\ell}$ be the boundary-aligned Haar wavelet with support width
$h_j=2^{-j}$.  The theory and code use only isotropic atoms.  For
$\lambda=(j,\ell_1,\ell_2)$, with $j\ge J_0$, define
\begin{equation}
 \Psi_\lambda^S(u,v)=
 \frac{\psi_{j\ell_1}(u)\psi_{j\ell_2}(v)
      +\psi_{j\ell_1}(v)\psi_{j\ell_2}(u)}
 {\|\psi_{j\ell_1}\otimes\psi_{j\ell_2}
  +\psi_{j\ell_2}\otimes\psi_{j\ell_1}\|_2} .             \label{eq:symatom}
\end{equation}
Thus $\Psi_\lambda^S\in V_{J_0}^{\perp}$, its support has area of order
$h_j^2$, and $\|\Psi_\lambda^S\|_\infty\asymp h_j^{-1}$.
The primary coefficient is
\begin{equation}
  \theta_\lambda=\langle\delta,\Psi_\lambda^S\rangle_2.       \label{eq:theta}
\end{equation}
Its support gives two coordinate intervals and its resolution gives scale.  After
fixing the basis orientation, its sign is the direction of a contrast between signed
lobes, not a statement that the entire support has higher or lower connectivity.  Its
pointwise log-odds amplitude is of order $\theta_\lambda/h_j$ under the single-atom
alternative below.  The inferential target is
$\btheta_n=(\theta_\lambda:\lambda\in\cI_n)$ over a prespecified isotropic
multiscale dictionary.

For the headline detection theory we use single-atom alternatives
\begin{equation}
 f=f_0+\theta\Psi_\lambda^S.                                  \label{eq:alternative}
\end{equation}
In the lower-bound and simulation parameterizations the same alternative is written
$f=f_0+\varsigma ah\Psi_\lambda^S$ with sign $\varsigma\in\{-1,1\}$ and amplitude
$a>0$.  The two are the same object under $\theta=\varsigma ah$, and every rate below
may be read in either.
Multiple or misshapen departures are stress tests, but they are outside the proposed
upper/lower comparison unless a localization margin and nuisance class are proved.

\subsection{Relative null}                                 \label{sec:scinull}

Equation~\eqref{eq:theta} fixes a formal target: $\theta_\lambda$ is nonzero
whenever $f$ has any component outside $V_{J_0}$ in the direction of
$\Psi^S_\lambda$, whether that component is a localized departure or an ordinary
smooth variation of connection probability along the coordinate.  The two are
not the same object, and the paper's hypotheses must say which one
is being tested.

\begin{definition}[Nuisance class and relative null]\label{def:scinull}
Let $\mathcal N$ be a closed linear subspace of $L^2([0,1]^2)$ of symmetric
surfaces.  The \emph{relative null}
relative to $\mathcal N$ is $H_0(\mathcal N)$: $f\in\mathcal N$.  A
\emph{departure} relative to $\mathcal N$ is $\delta_{\mathcal N}=f-\Pi_{\mathcal N}f$,
with $\Pi_{\mathcal N}$ the $L^2$ projection onto $\mathcal N$, and the relative
target of atom $\lambda$ is $\theta_\lambda(\mathcal N)=\langle\delta_{\mathcal N},\Psi^S_\lambda\rangle$.
\end{definition}

Two classes are used.  The \emph{coarse class} $\mathcal N_0=V_{J_0}$ is the one
for which every theorem of Section~\ref{sec:results} is proved.  Under it
smooth distance decay is a departure, and for a diagonal atom of width $h$ on a
surface whose logit falls with slope $\beta$ in $|u-v|$ the target is of order
$\beta h^2$, not zero.  The \emph{combined class} $\mathcal N_1$ of
\eqref{eq:union} contains $\mathcal N_0$ together with smooth distance decay and
additive main effects.  Under it that same decay is background.  A rejection of
$H_0(\mathcal N_0)$ on a smooth surface is correct and empty.  A
rejection of $H_0(\mathcal N_1)$ is the claim a practitioner wants.  Every
table in this paper names the class its null refers to, and the comparisons of
Section~\ref{sec:sim} are comparisons of estimands before they are
comparisons of estimators.  Enlarging a nested linear class $\mathcal N$ reduces the norm of the
residual $\delta_{\mathcal N}$ but does not act coefficient by coefficient: a
single residual coefficient can be created or reversed by the enlargement
(take $\phi_1,\phi_2$ orthonormal and orthogonal to $\mathcal N_0$,
$\mathcal N_1=\mathcal N_0+\operatorname{span}(\phi_1+\phi_2)$ and
$f=\phi_2$.  The $\phi_1$-coefficient is $0$ under $\mathcal N_0$ and $-1/2$
under $\mathcal N_1$).  $H_0(\mathcal N_1)$ is the weaker hypothesis, and the
power against it is the relevant power.  The individual atoms that carry it can
change with the class.

\subsection{Nuisance classes}
\label{sec:nuisance}

The construction above takes $\mathcal N=\mathcal N_0=V_{J_0}$, where the
projection identity $\Pi_{J_0}^\top c_\lambda=0$ removes the background from
every scanned coefficient exactly, so the background cannot absorb a
departure.  The class assumes that connection probability is constant within a
coarse cell.  Ordered regimes and categorical coordinates produce block
structure, and for them $\mathcal N_0$ is the right relative null.  Where the
coordinate is continuous and connection probability varies smoothly along it,
as in both applications here, a block surface cannot represent the variation,
the residual decay is indistinguishable from a localized departure on the
diagonal, and the scan reports it as one.  Refining $J_0$ does not repair
this, because the approximation error concentrates on the diagonal at every
resolution.  A smooth class fails in the mirror direction, since a spline
expansion cannot represent block jumps.  Section~\ref{sec:sim}
measures both failures.

We therefore fit the background in the linear span of the coarse block
indicators and smooth terms,
\begin{equation}
\operatorname{logit}p_0(u,v)=\alpha+s(|u-v|)+m(u)+m(v)+\textstyle\sum_{a\le b}\beta_{ab}\1\{(u,v)\in B^S_{ab}\},
\label{eq:union}
\end{equation}
with $B^S_{ab}$ the coarse symmetric blocks of $V_{J_0}$ (not the fine cells
$Q^S_{ab}$ of Section~\ref{sec:method}), $J_0=2$ throughout the executed
work, and $s$ and $m$ piecewise-linear spline expansions whose knots are
spaced more coarsely than the finest scanned atom.  We call \eqref{eq:union}
the combined class.  It nests both alternatives, so it is correctly specified
whenever either is.  The additive terms $m(u)+m(v)$ are orthogonal to every
same-scale interaction atom, so main effects cannot absorb a departure, and
$s$ is orthogonal to second order in the atom width against the $O(1)$ error a
block class leaves.  The background is fitted on its own dyad fold, so it is
independent of the inference fold, and the atom coefficient is
\begin{equation}
\widehat\theta_\lambda=\sum_k c_{\lambda k}\{\operatorname{logit}\widetilde p_k-\operatorname{logit}\bar p_{0k}\},
\label{eq:unioncoef}
\end{equation}
with $\bar p_{0k}$ the cell-average fitted probability.  The delta method gives
the first-order estimation contribution $g_\lambda^\top\Sigma_\beta g_\lambda$,
added to the Bernoulli term, with a second-order remainder that is not
controlled here.  Comparing the cell-average fitted logit instead of the logit
of the cell-average probability leaves a Jensen gap that is largest on the
diagonal and is by itself enough to fire the diagonal atoms.

Two limits of the construction are stated once.  Definition~\ref{def:scinull}
uses the Lebesgue $L^2$ projection, while a likelihood-fitted background
targets the population logistic projection under the design measure, which
coincides with the $L^2$ projection when the class contains the truth and
differs from it under alternatives, where the fitted background can absorb
part of a departure.  And the span is finite and contains neither every smooth
background nor every block one.  Every theorem in Section~\ref{sec:results},
the level study, the regime map, the panels and both applications use the
coarse class, whose simulated backgrounds are block-constant by construction.
The combined class is reported in Section~\ref{sec:sim}, in the
regime map and on both applications, measured and not proved, and for a real
network it is the class to use.

\subsection{The dictionary}                                  \label{sec:dict}

\begin{figure}[t]\centering
\includegraphics[width=\textwidth]{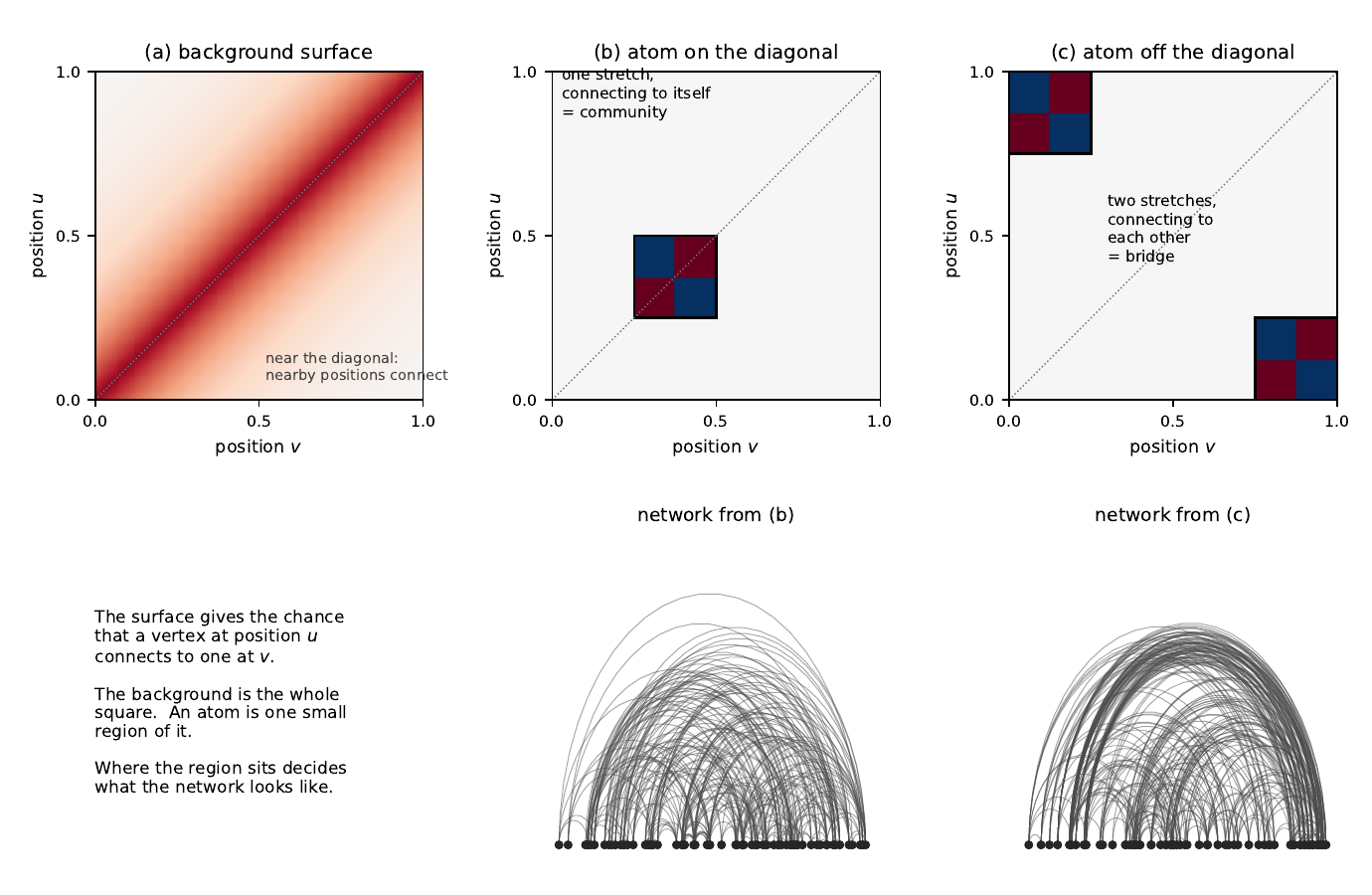}
\caption{The connection surface over pairs of positions, with red for a raised
and blue for a lowered probability.  (a) A background in which nearby
positions connect more.  (b) An atom on the diagonal, covering pairs within
one stretch, which produces a community split at its middle.  (c) An atom off
the diagonal, covering pairs that span two separated stretches, which produces
a bridge.  The lower row shows one network drawn from each, with vertices
placed on the axis and edges as arcs.}
\label{fig:plane}
\end{figure}

The connection probability is a surface over pairs of positions.  With a
one-dimensional coordinate that surface lives on the square $[0,1]^2$, and its
height at $(u,v)$ is the chance that a vertex at $u$ connects to one at $v$.
A localized departure from the background is a departure confined to a small
region of that square, and where the region sits decides what it produces
(Figure~\ref{fig:plane}).  A region on the diagonal covers pairs whose two
vertices lie in the same stretch of the axis, so raising the surface there
makes that stretch connect to itself, which is a community with a boundary at
the middle of the region.  A region off the diagonal covers pairs whose
vertices lie in two separated stretches, so raising the surface there makes
those two stretches connect to each other, which is a bridge.  A community and
a bridge are therefore the same kind of object, a departure confined to a
region of the square, and they differ only in where that region lies.

Write $\psi_{j,l}$ for the Haar function on the dyadic interval of scale $j$
and index $l$, and for a dyadic box $a$ at scale $j$ in dimension $d$ put
$\phi_a=\bigotimes_{k\le d}\psi_{j,l_k}$, normalized so that $\int\phi_a^2=1$.
The dictionary is the symmetrized family
\begin{equation}
\Psi_{a,b}(u,v)=\frac{\phi_a(u)\phi_b(v)+\phi_b(u)\phi_a(v)}{\sqrt{2(1+\1\{a=b\})}},
\qquad a,b\ \text{boxes at a common scale},
\label{eq:dict}
\end{equation}
together with the coarse space $V_{J_0}$ carrying the background.  We call
$\Psi_{a,b}$ an \emph{atom}.  An atom with $a=b$ is \emph{diagonal} and covers
the region $a\times a$, an atom with $a\ne b$ is \emph{off-diagonal} and covers
$(a\times b)\cup(b\times a)$, and $h=2^{-j}$ is the width of the atom.  The
generative model used throughout Section~\ref{sec:results} is
\begin{equation}
\operatorname{logit}p(u,v)=f_0(u,v)+\theta\,\Psi(u,v),
\label{eq:gen}
\end{equation}
with $f_0$ the background and $\theta$ the amplitude.

\subsection{Choice of basis}                                 \label{sec:whybasis}

Four properties are required of the functions that carry the departures, and
the Haar system has all four.  Each is used by name in
Section~\ref{sec:results}.

\emph{Localized support.}  A coefficient is to be a statement about one place,
so the function carrying it must vanish outside that place.  A basis whose
functions are spread over the whole coordinate, such as the Fourier system,
polynomials or global splines, produces coefficients that describe the whole
coordinate, and no coefficient of such a basis can be read as a claim about a
stretch.

\emph{Indexed by width.}  The dictionary contains every dyadic width at once,
so the width of a departure is estimated and is not fixed in advance.  A
single-bandwidth kernel or a block model at one resolution requires the width
to be chosen before the data are seen.

\emph{Orthogonality.}  Distinct atoms are orthogonal, so their coefficients do
not contaminate one another and a departure at one place and width does not
induce a spurious coefficient elsewhere.  This is what makes testing the whole
dictionary at once meaningful, and it costs only the factor
$\sqrt{2\log 2m_n}$ in the critical value.

\emph{Zero mean.}  Every wavelet integrates to zero, and this makes every atom
exactly orthogonal to the additive surfaces $m(u)+m(v)$ that carry degree
heterogeneity, by Proposition~\ref{prop:dict}(ii).  A basis without this
property, for instance the scaling function of the same wavelet family, gives
atoms whose coefficients respond to degree at first order, so a network in
which some regions are simply better connected would produce discoveries.

A fifth consideration decides which wavelet to use among those that satisfy
the four.  The zero-mean property annihilates a constant background within an
atom's support, and a wavelet with $p$ vanishing moments,
$\int u^{k}\psi=0$ for $k<p$, annihilates a polynomial of degree $p-1$.  The
Haar wavelet has $p=1$: it kills constants and leaves a linear trend, and
$\int u\psi=-0.031$ at width $h=1/4$.  The Daubechies wavelet with two
vanishing moments has $p=2$, with $\int u\psi$ zero to machine precision and
only the quadratic term surviving at $-0.001$.  The consequence is the size of
the bound in Proposition~\ref{prop:dict}(iv): the leakage of a smooth
separation kernel into a diagonal atom is $O(h^{d+1}\|s'\|_\infty)$ for $p=1$
and $O(h^{d+2}\|s''\|_\infty)$ for $p=2$.  Measured on an exponential kernel
$s(t)=e^{-\lambda t}$ at $h=1/4$, the inner product $|\langle s,\Psi\rangle|$
is $8.68$, $14.57$ and $21.02$ for the Haar wavelet at $\lambda=2,4,8$ against
$1.59$, $3.09$ and $5.68$ for the Daubechies wavelet, a factor between $3.7$
and $5.5$.

We use the Haar wavelet throughout for three reasons.  Its support is exactly
the dyadic box, so a coefficient refers to a stretch with no spill into
neighbouring boxes, which matters for a procedure whose output is a location.
It is constant on cells, so the coefficient is an exact linear combination of
cell counts and the studentizer is exact.  And its integrals are available in
closed form, which is what makes Theorem~\ref{thm:map} explicit.  A wavelet
with two vanishing moments is the better choice when the background carries a
trend within an atom's support that the fitted background class does not
capture, and the theory covers that case unchanged: the functionals in
Theorem~\ref{thm:map}(ii) and (iv) do not depend on the shape of the mother
wavelet, and only $\kappa_4$ in part (iii) changes.

\section{Estimator and Calibration}                       \label{sec:method}

\subsection{Dyad folds and cells}

Before observing either missingness or outcomes, hash every eligible unordered dyad
directly into three mutually exclusive masks $R^{B^T},R^{B^E},R^I$ for baseline
tuning, baseline estimation, and inference, with proportions $0.05$, $0.15$ and
$0.80$.  The hash, proportions, and seed are fixed across the theory,
simulations, and applications.  External coordinate data are processed before
these masks are opened, and the coordinate is never replaced by a degree or
spectral order of $A$.  An earlier version of the procedure reserved a fourth
fold for an ordering correction.  No step consumed it, and it is gone.  This
section describes the procedure that was run.  The pilot configuration it
replaced, and the evidence for the replacement, are in Section~\ref{sec:sim}
and the supplement.

Let $I_{Ja}=[a2^{-J},(a+1)2^{-J})$ and let $Q_{ab}^S$ be the symmetric dyadic cell generated by $I_{Ja}\times I_{Jb}$.  Its width $2^{-J}$ is smaller than the narrowest scanned feature and tends to zero.  For $F\in\{B^T,B^E,I\}$, let $\pi_{ij}$ be the observation propensity, $W_{ij}^F=R_{ij}^FM_{ij}/\pi_{ij}$, and define
\begin{align}
 D_{ab}^{F}&=\sum_{i<j}W_{ij}^{F}\1\{(\widehat U_i,\widehat U_j)\in Q_{ab}^S\},\nonumber\\
 N_{ab}^{F}&=\frac{(D_{ab}^{F})^2}
 {\sum_{i<j}(W_{ij}^{F})^2\1\{(\widehat U_i,\widehat U_j)\in Q_{ab}^S\}},\nonumber\\
 \bar A_{ab}^{F}&=(D_{ab}^{F})^{-1}
 \sum_{i<j}W_{ij}^{F}A_{ij}\1\{(\widehat U_i,\widehat U_j)\in Q_{ab}^S\}.
                                                                    \label{eq:cellcounts}
\end{align}
Here $N_{ab}^F$ is the inverse-weight effective count.  Set $N_{ab}^F=0$ when
its denominator is zero.  The proportion
$\bar A_{ab}^F$ is undefined when $D_{ab}^F=0$.  Any atom requiring such an
unsupported inference cell is inadmissible.  Cell ratios estimate observed-design
means.  They do not by themselves correct to Lebesgue measure.  The target requires a
uniform shrinking-cell occupancy and quadrature argument.  Lebesgue cell areas enter
only at projection.  Unknown propensities must be estimated without inference-fold
outcomes, and their error enters the remainder.

Select the bandwidth from a finite prespecified grid using $B^T$ deviance, then hold it fixed and fit the symmetric local-linear smoother to $\bar A_{ab}^{B^E}$.  Smoothing occurs before any logit transform.  We then stabilize, take logits, project onto exactly $V_{J_0}$, and transform back, giving $\widehat p_{0,ab}^{B^E}$.  This final projection makes the implemented baseline target agree with \eqref{eq:projection}.  Without it, nonlinear smoothing could retain tested details.  Empty estimation cells borrow information through the frozen smoother.  A cell unsupported by its neighborhood is inadmissible instead of assigned zero.

On the inference fold use the stabilized proportion
\begin{equation}
 \widetilde p_{ab}^I=
 \frac{N_{ab}^I\bar A_{ab}^I+\kappa_0\widehat p_{0,ab}^{B^E}}
 {N_{ab}^I+\kappa_0},                                                       \label{eq:stabilize}
\end{equation}
where $\kappa_0$ is fixed before analysis.  Clip both $\widetilde p_{ab}^I$ and $\widehat p_{0,ab}^{B^E}$ to $[\epsilon_n,1-\epsilon_n]$, with $\epsilon_n$ tied to the sparse envelope and reported.  A completed proof must control the clipping-hit probability and the bias ratios $\epsilon_n/\rho_n$ and $\kappa_0/(N_{ab}^I\rho_n)$.  The cellwise residual logit is
\begin{equation}
 \widehat d_{ab}=\ell_\epsilon(\widetilde p_{ab}^I)-\logit(\widehat p_{0,ab}^{B^E}),\label{eq:reslogit}
\end{equation}
with $\ell_\epsilon$ the clipped logit of Section~\ref{sec:method}.  The clip
is inactive on the good event (Supplement Lemma~8.2(d) for the instance), where
$\widehat d_{ab}$ agrees with the unclipped residual, and the supplement's
estimator contract is stated in terms of these residuals, so \eqref{eq:thetahat}
is their cellwise aggregation and no second object is introduced.

\subsection{Coefficients and maximum statistic}

Index symmetric cells by $q\in\cQ_J$, write
$c_{\lambda q}=|Q_q^S|\bar\Psi_{\lambda q}$, and define the clipped logit
$\ell_\epsilon(x)=\logit\{\min(1-\epsilon_n,\max(\epsilon_n,x))\}$.
The coefficient estimator is
\begin{equation}
 \widehat\theta_\lambda
 =\sum_{q\in\cQ_J}c_{\lambda q}
  \{\ell_\epsilon(\widetilde p_q^I)
    -\logit(\widehat p_{0q}^{B^E})\}.                           \label{eq:thetahat}
\end{equation}

The standard error $\widehat s_\lambda$ is the model-based two-fold sandwich
built from the implemented Jacobian of the frozen baseline map along both
stochastic paths, with fitted Bernoulli cell variances $\widehat V_r^F$.  The
exact gradient, loading, and meat formulas, verified against finite differences,
constitute the estimator-to-code contract in Appendix~\ref{sec:supp-finite}.
Theorem~\ref{thm:linear} isolates the curvature and studentizer
errors this construction must control, and
Supplement Propositions~S2 and~S3 bound them by checkable
envelopes.

Let $Z_\lambda=\widehat\theta_\lambda/\widehat s_\lambda$ for testing $H_{0\lambda}:\theta_\lambda=0$, and let
\begin{equation}
 T=\max_{\lambda\in\cI_n}|Z_\lambda|.                              \label{eq:T}
\end{equation}
For bootstrap draw $b$, generate independent Gaussian cell scores
$\xi_r^{F,(b)}\sim N(0,\widehat V_r^F)$ for $F\in\{I,B^E\}$ and reuse each score
for every atom.  Set
\begin{equation}
 Z_\lambda^{*(b)}=\widehat s_\lambda^{-1}
 \left\{\sum_rG_{\lambda r}^{I}\xi_r^{I,(b)}
 +\sum_rG_{\lambda r}^{B}\xi_r^{B,(b)}\right\},
 \quad T^{*(b)}=\max_\lambda|Z_\lambda^{*(b)}|.                 \label{eq:multiplier}
\end{equation}
Because every derivative is constant within a cell, this cell-aggregated Gaussian
draw has the same conditional distribution as independent dyad-level Gaussian
multipliers applied to the model-based influence terms.  It is also the exact
implementation used in the simulations.
This bootstrap is intended to approximate the centered maximum
$\max_\lambda|\widehat\theta_\lambda-\theta_\lambda|/s_\lambda$.
the observed $T$ in \eqref{eq:T} equals that maximum only under the global null.
With $B_*\to\infty$, let $c_{1-\alpha}^*$ be the empirical
$(1-\alpha)$ quantile of $T^{*(1)},\ldots,T^{*(B_*)}$.  It yields the candidate
simultaneous intervals
\begin{equation}
 \widehat\theta_\lambda\ \pm\ c_{1-\alpha}^*\widehat s_\lambda.    \label{eq:simci}
\end{equation}

\subsection{Admission rule and computation}

Let $S_{\lambda,+}=\{\Psi_\lambda^S>0\}$ and
$S_{\lambda,-}=\{\Psi_\lambda^S<0\}$.  Using only $B^T$ outcomes and
non-outcome denominators, compute the same plug-in influence map as above separately
on each signed lobe.  Let $I_{\lambda,+}$ and $I_{\lambda,-}$ be the inverses of
the resulting pilot variances and let $L_\lambda$ be the largest absolute
standardized pilot dyad contribution.  An atom is scanned only if \emph{every fine cell} in its support carries at
least five pilot expected successes and five expected failures, both lobes
carry enough information, no single dyad dominates, and its width exceeds the
reported coordinate uncertainty.  In compact notation, with $E_{k}^{\pm}$ the
pilot expected successes and failures of cell $k$,
\begin{equation}
 \cI_n=\left\{\lambda:
 \min_{k\in\operatorname{supp}\lambda}E^{+}_{k}\wedge E^{-}_{k}\ge5,\quad
 I_{\lambda,+}\wedge I_{\lambda,-}\ge I_{\min,n},\quad
 L_\lambda\le L_n,\quad
 r_n/h_\lambda\le c_r
 \right\}.                                                         \label{eq:gates}
\end{equation}
Every quantity in \eqref{eq:gates} is computed from pilot probabilities and
denominators, so the family is fixed before inference outcomes are read.  The
executed runs added one further exclusion that is \emph{not} outcome-free: an
atom touching a cell whose realized inference probability was clipped is
dropped.  It fires on about $0.01$ atoms per fit and the calibration theory
does not cover it.  Table~\ref{tab:calib} reports size with it \emph{off},
at eight configurations including every one at which power is reported.  The
power studies were run with it on, a difference of about one excluded atom
per hundred fits that we treat as immaterial and say so.
The order gate is used at the frozen value $c_r=0.25$.  That constant is an
empirical rule fixed before the executed runs, not a consequence of the theory:
Appendix Theorem~\ref{thm:order} requires $r_n/h_{\min}\to0$, and sign and magnitude
transfer require in addition an information-dependent bias condition
$Br_nn\sqrt{\rho_n\ell_n}\to0$, so a fixed ratio is insufficient as $a_n\to0$
(Appendix~\ref{sec:supp-detect}, condition (t3-13)).  Deriving a vanishing, information-aware gate is an open
problem stated in Section~\ref{sec:thy-limits}.  All gates use non-outcome
denominators and pilot baseline predictions from $B^T$,
never $B^E$ or inference outcomes.  The largest admissible $J$ is therefore data
adaptive while the index set remains fixed for the two stochastic folds.  A proof
must show that pilot admission implies the corresponding population information and
curvature conditions with probability tending to one.

\begin{center}
\fbox{\begin{minipage}{0.93\linewidth}
\textbf{Algorithm 1: Symmetric multiscale localization.}
Given $A,M$, the external $\widehat U$, a fixed fold hash, candidate grids,
and atoms:
\begin{enumerate}[leftmargin=1.6em,itemsep=1pt,topsep=3pt]
\item orient $\widehat U$ by the external anchor and obtain $r_n$ from the order model.
\item construct the symmetric fold-specific counts in \eqref{eq:cellcounts}.
\item select tuning on $B^T$ and fit the symmetric design-corrected baseline on $B^E$.
\item stabilize and clip inference probabilities and compute \eqref{eq:reslogit}.
\item form $(\widehat\theta_\lambda,\widehat s_\lambda)$ for atoms passing \eqref{eq:gates}.
\item generate shared $B^E$- and $I$-dyad multiplier maxima and obtain $c_{1-\alpha}^*$.  And
\item return simultaneous intervals, signed supports, scales, amplitudes, and diagnostics.
\end{enumerate}
\end{minipage}}
\end{center}

The dense reference implementation runs in $O\{m+LQ^2+MQ^2+B_*MQ\}$ time and
$O(Q^2+MQ)$ memory for $m$ observed dyads, $Q$ cells, $M$ atoms, $L$ bandwidth candidates, and $B_*$
multiplier draws (Appendix~\ref{sec:supp-finite}).  No fast wavelet or separable shortcut is
claimed.

\section{Theory}                                          \label{sec:results}

Five results are stated here and proved in the Appendix, which also contains
the supporting conditional reduction and multiplier calibration
(Appendix~\ref{sec:supp-finite} and~\ref{sec:supp-mult}), the assumed-rate
form of the detection boundary, the ordering-transfer theorem, the all-radius
localization floor and the rank-identifiability bound.  Proposition~\ref{prop:dict}
concerns the dictionary as a basis, Theorem~\ref{thm:map} the networks its
atoms generate, Theorem~\ref{thm:uncond} the calibration and the detection
boundary of the scan at an exact coordinate, and Theorem~\ref{thm:info} and
Proposition~\ref{prop:obstruct} the two limits, on coordinate dimension and on
coordinate accuracy.

\subsection{Assumptions}
This section states four theorems and one proposition, with proofs in the
Supplement: Theorem~\ref{thm:uncond} in Appendix~\ref{sec:supp-detect},
Appendix Theorem~\ref{thm:order} and Proposition~\ref{prop:obstruct} in
Appendix~\ref{sec:supp-order} together with Appendix Theorem~\ref{thm:floor}, and
Theorem~\ref{thm:uncond} in Appendix~\ref{sec:supp-uncond}.  The calibration
engine on which Theorem~\ref{thm:uncond} rests, a conditional
Bernstein--Taylor reduction and a shared-dyad multiplier calibration, is
stated and proved in Appendix~\ref{sec:supp-finite}, Appendix~\ref{sec:supp-mult}
and Appendix~\ref{sec:supp-mult}.  Theorem~\ref{thm:uncond} assume displayed linearization, studentizer, covariance,
packing and score-transfer rates.  Those rates are derived from primitive
conditions only in the instance of Theorem~\ref{thm:uncond}, and population
inference requires an explicit finite-design-to-Lebesgue bias condition,
without which the inferential centre is the frozen finite-design target.

Write $\ell_n=\log(2m_n\vee3)$.  The supplement's per-section variants
$\log(2m_n)$, $\log(e\bar m_n)$, and $\log\bar d_n$ agree with $\ell_n$ up to
bounded ratios on the admission event, every rate below is insensitive to the
choice, and each supplement section states its local variant at first use.
Similarly $m_n$ is the deterministic packing cardinality here and the realized
admitted cardinality where a supplement section says so.  Several symbols are
overloaded across sections and are always disambiguated by context or subscript:
$\delta$ is the departure $(I-\Pi_{V_{J_0}})f$, $\delta_n$ the fine-cell width,
$\delta_{V,n},\delta_{G,n}$ studentizer errors, and grid ``$\delta$'' in
Section~\ref{sec:sim} the simulation cell width.  $\Gamma_{\lambda,n}$ is the
curvature envelope, $\Gamma_n$ the transfer margin of
Appendix Theorem~\ref{thm:order}, and $\Gamma_{ab,n}$ the multiplier covariance.  $H_0$
is the baseline Jacobian and $H=nh_*$ the block length of
Proposition~\ref{prop:obstruct}.  $R^F$ are fold masks, $R$ in Section~\ref{sec:sim}
a replication count and in Section~\ref{sec:thy-limits} a coordinate range.  $k$ indexes
dyads in Theorem~\ref{thm:linear} and cells elsewhere.  $\Delta_n$ is the realized
coordinate error in Section~\ref{sec:model} and the realized rank error in
Appendix Theorem~\ref{thm:order}, and $\kappa$ is the standardized amplitude of the
simulations, distinct from $\kappa_0$ and $\kappa_n^{\mathrm{str}}$.  The finite-design and calibration results permit an
$\mathcal F_{0n}$-measurable admitted family $\cI_n$.  The global minimax result instead
uses a deterministic master packing and an all-pass admission event, as stated below.
All atoms are boundary-aligned, isotropic, symmetric Haar atoms.  Their discontinuity
is handled by the boundary-layer inequality
\begin{equation}
 \|\Psi_\lambda^S(\,\cdot-t)-\Psi_\lambda^S\|_1
 \le C\|t\|_\infty,
 \qquad \|t\|_\infty\le c h_\lambda,                            \label{eq:bvtranslate}
\end{equation}
not by differentiating the atom.  Off-diagonal atoms use the normalization in
\eqref{eq:symatom}, with transposed duplicates removed.

\subsection{Requirements on the basis}                    \label{sec:thy-basis}

Let $\mathcal D$ denote the family \eqref{eq:dict} over all scales and let
$\mathcal G_{J_0}$ be the span of constant surfaces, additive surfaces
$m(u)+m(v)$, separation surfaces $s(\|u-v\|)$ with $s$ continuously
differentiable, and block surfaces whose blocks are unions of dyadic boxes at
resolution $J_0$.

\begin{proposition}[Dictionary properties]                  \label{prop:dict}
Let $\mathcal D$ denote the family \eqref{eq:dict} over all scales.
\begin{enumerate}[label=\textup{(\roman*)},leftmargin=*]
\item \textup{Completeness.}  $\mathcal D\cup V_{J_0}$ is an orthonormal basis
of the symmetric subspace of $L^2([0,1]^d\times[0,1]^d)$.  Every symmetric
square-integrable departure has an expansion in the dictionary, and a departure
orthogonal to every atom lies in $V_{J_0}$, that is, it is background.
\item \textup{Exact orthogonality to the nuisance.}  For every
$m\in L^2([0,1]^d)$ and every atom, $\langle m(u)+m(v),\Psi_{a,b}\rangle=0$.
Degree heterogeneity therefore has zero coefficient on every atom of the
dictionary and cannot be mistaken for localized structure.
\item \textup{Sparsity of localized alternatives.}  If the departure $\delta$
is bounded by $B$ and supported on a set $S$ of measure $A$, then at scale
$h=2^{-j}$ at most $2A/h^{2d}+O(h^{-d}|\partial S|)$ atoms have a nonzero
coefficient, and $\sum_{a,b}\langle\delta,\Psi_{a,b}\rangle^2\le B^2A$.  The
energy of a localized departure is carried by a number of coefficients
proportional to the measure of its support, not to the size of the dictionary.
\item \textup{Invisibility of global surfaces.}  Let $\Psi_{a,b}$ be a atom at scale $h=2^{-j}$ with $j>J_0$, and
let $\delta\in\mathcal G_{J_0}$ with $s$ continuously differentiable.  Then
\[
\langle\delta,\Psi_{a,b}\rangle=0\quad\text{for }a\ne b,
\qquad
\bigl|\langle\delta,\Psi_{a,a}\rangle\bigr|\le C_d\,h^{\,d+1}\|s'\|_\infty .
\]
Consequently every member of $\mathcal G_{J_0}$ has the same dictionary
coefficients, namely zero up to $O(h^{d+1})$, and no statement of the form
``a departure of width $h$ is present at position $a$'' is identifiable within
that class.  By contrast $\theta_\lambda=\langle\delta,\Psi_\lambda\rangle$ is
exactly the coefficient the scan estimates, so within the dictionary the
position and the width are identified.
\end{enumerate}
\end{proposition}

Part (i) answers the objection that the structure of interest might lie
outside the dictionary: the atoms and the coarse space together span the
symmetric $L^2$, so a departure invisible to every atom is background by
definition.  Part (ii) is exact and is what separates signal from nuisance,
since degree heterogeneity has zero coefficient on every atom.  Part (iii) is
the formal content of the word localized.  Estimating a surface at resolution
$h$ costs $h^{-2d}$ parameters whatever is present, whereas a departure of
support measure $A$ has $O(A/h^{2d})$ nonzero coefficients, so a scan over the
dictionary pays for the unused remainder only through the multiplicity
factor $\sqrt{2\log2m_n}$.  Part (iv) states the consequence for the
canonical generative mechanisms: an Erd\H os--R\'enyi surface, a
degree-corrected surface, a latent-distance or small-world kernel, and a
coarse-block model all have the same dictionary coefficients, namely zero up
to $O(h^{d+1})$, so a statement of the form ``a departure of width $h$ sits at
position $a$'' is not identifiable within those classes.  Section~\ref{sec:nuisancestudy}
measures parts (ii) and (iv) on simulated networks and
Section~\ref{sec:apps} relies on them in every application.

\subsection{Network summaries generated by an atom}                        \label{sec:thy-gen}

Throughout this subsection the coordinates are uniform on $[0,1]^d$, the
background $f_0$ in \eqref{eq:gen} is constant, and we write $\rho=\sigma(f_0)$,
$w=\rho(1-\rho)$ and $w'=\rho(1-\rho)(1-2\rho)$, so that
\begin{equation}
p_\theta=\rho+\theta w\Psi+\tfrac12\theta^{2}w'\Psi^{2}+O(\theta^{3}).
\label{eq:expand}
\end{equation}
In the sparse regime $\rho<\tfrac12$ and $w'>0$.  We write
\begin{equation}
\kappa_4=\int_0^1\psi^{4},\qquad\text{so that}\qquad\int\phi_a^{4}=\kappa_4^{\,d}
\label{eq:kappa4}
\end{equation}
by Fubini.  Fix a partition of $[0,1]^d$ into congruent dyadic blocks
$B_1,\dots,B_K$ at a resolution at least as fine as $h$.  Let $Q$ be the
modularity of that partition, $\mathrm{ES}$ the mean coordinate separation
across edges divided by its value across all pairs, $T$ the transitivity,
$\mathrm{CV}$ the coefficient of variation of the expected degree, and
$\lambda(\Psi)$ the spectrum of $\Psi$ as an integral operator.

Section~\ref{sec:catalogue} reported measured signatures.  This section derives
them.  Each network summary is an explicit functional of the dictionary
atom, computable in closed form in any dimension, and the functionals
separate the families independently of amplitude and of scale.

\begin{theorem}[Generative property map]\label{thm:map}
Under \eqref{eq:gen} with uniform coordinates and constant background, as
$\theta\to0$:
\begin{enumerate}[label=\textup{(\roman*)},leftmargin=*,itemsep=2pt]
\item \textup{(Functionals.)} With $\mu=\mathbb E\|u-v\|$,
\begin{align}
\mathbb E[Q]&=\frac{\theta w}{\rho}\sum_{k\le K}\int_{B_k\times B_k}\!\!\Psi\;+\;O(\theta^{2}),
\label{eq:mapQ}\\[2pt]
\mathbb E[\mathrm{ES}]&=1+\frac{\theta w}{\rho\,\mu}\bigl\langle\|u-v\|,\Psi\bigr\rangle+O(\theta^{2}),
\label{eq:mapES}\\[2pt]
\mathrm{CV}&=\frac{\theta^{2}w'}{2\rho}\Bigl(\textstyle\int\phi_a^{4}-1\Bigr)^{1/2}+O(\theta^{3}),
\label{eq:mapCV}\\[2pt]
\mathbb E[T]&=T_0+\frac{(\theta w)^{3}}{\rho^{3}}\operatorname{tr}(\Psi^{3})+O(\theta^{4}),
\label{eq:mapT}
\end{align}
and the spectral gap of the expected normalized Laplacian is monotone in
$\max\lambda(\Psi)$ and its multiplicity.  The first two expansions have no
zeroth-order correction because $\int\Psi=0$; the last two are the leading
nonzero terms.
\item \textup{(Communities and bridges.)} Let $a\ne b$ be boxes at scale $h$ and $\theta>0$.
\begin{enumerate}[label=\textup{(\alph*)},leftmargin=*]
\item \emph{Diagonal, $\Psi=\phi_a\otimes\phi_a$.}
$\sum_k\int_{B_k\times B_k}\Psi=(h/2)^{d}>0$,
$\langle\|u-v\|,\Psi\rangle=-c_dh^{d+1}<0$ with $c_1=1/6$,
$\operatorname{tr}(\Psi^{3})=1$, and $\lambda(\Psi)=\{1\}$ of multiplicity one.
Modularity rises, edges shorten, transitivity rises, the perturbation is rank
one.
\item \emph{Off-diagonal, $\Psi=(\phi_a\otimes\phi_b+\phi_b\otimes\phi_a)/\sqrt2$.}
$\sum_k\int_{B_k\times B_k}\Psi=0$, $\langle\|u-v\|,\Psi\rangle=0$,
$\operatorname{tr}(\Psi^{3})=0$, $\lambda(\Psi)=\{\pm1/\sqrt2\}$.  Both
first-order effects vanish.  At second order
$\mathbb E[Q]=-\theta^{2}w'/2\rho<0$ and
$\mathbb E[\mathrm{ES}]=1+(\theta^{2}w'/2\rho)(\Delta-\mu)/\mu>1$ whenever the
boxes are separated by more than an average pair, $\Delta>\mu$.
\end{enumerate}
\item \textup{(Shape.)} For a diagonal atom built from any mother wavelet $\psi$ with $\int\psi=0$
and $\int\psi^{2}=1$, the quantities $\sum_k\int_{B_k\times B_k}\Psi$,
$\langle\|u-v\|,\Psi\rangle$, $\operatorname{tr}(\Psi^{3})$ and
$\lambda(\Psi)$ do not depend on the shape of $\psi$.  The only
shape-dependent summary is degree heterogeneity, through the fourth moment
\eqref{eq:kappa4}:
\[
\mathrm{CV}=\frac{\theta^{2}w'}{2\rho}\bigl(\kappa_4^{\,d}-1\bigr)^{1/2}+O(\theta^{3}).
\]
\item \textup{(Mixtures.)} Let $\Psi=\sum_{k\le M}c_k\,\phi_{a_k}\otimes\phi_{a_k}$ with
$\sum_kc_k^{2}=1$ and the $\phi_{a_k}$ orthonormal.  Then
$\operatorname{tr}(\Psi^{3})=\sum_kc_k^{3}$ and
$\lambda(\Psi)=\{c_1,\dots,c_M\}$.  In particular, $M$ atoms at one scale
with equal weights give $\operatorname{tr}(\Psi^{3})=M^{-1/2}$ and $M$ equal
eigenvalues: a block model on intervals, with modularity the sum of the
single-element values, low transitivity and a large spectral gap.  Atoms at
$M$ different scales with equal weights give the same trace when their supports
are disjoint and strictly more when nested, because the operator products no
longer vanish across scales: a hierarchy has higher transitivity than a block
model of the same rank, and a smaller spectral gap.
\end{enumerate}
\end{theorem}

The theorem reduces the catalogue of Section~\ref{sec:catalogue} to four
integrals.  In part (ii), within-block integrals of a diagonal atom are
squares and cannot cancel, which forces modularity upward, while an
off-diagonal atom is a swap whose cube has zero diagonal, which removes its
first-order effects entirely and drives modularity negative at second order.
The two families therefore have opposite signs in modularity and in edge
span at every scale and in every dimension.  Part (iii) says that position on
the diagonal, not the shape of the wavelet, decides whether an atom builds a
community or a bridge, and that shape enters only through the fourth moment
$\kappa_4$, which is $4$ for the Haar wavelet, $6.78$ for the Mexican hat and
$8.75$ for the Daubechies wavelet with two vanishing moments.  Dimension
enters as $\kappa_4^{\,d}$.  Part (iv) separates a block model from a
hierarchy of the same rank through $\operatorname{tr}(\Psi^3)$, which is
$M^{-1/2}$ for $M$ disjoint atoms and strictly larger for nested ones.
Table~\ref{tab:map} places each closed form beside its measured value at a
finite amplitude, and Section~\ref{sec:gencheck} tests the same functionals
on real networks by drawing from the fitted model.

\begin{table}[t]\centering
\caption{Closed forms from Theorem~\ref{thm:map} beside the measured value on
networks of $n=700$ at mean degree $20$ drawn from \eqref{eq:gen} at matched
signal energy, mean of five replicates.  Closed forms are the leading nonzero
term in the amplitude and hold at every scale and dimension.  Edge span is the
mean coordinate separation across edges divided by its value over all pairs.
Modularity is of eight equal intervals.  Transitivity for a diagonal atom is
$T_0$ plus a term proportional to $\operatorname{tr}(\Psi^3)=1$, and for an
off-diagonal atom is $T_0$.}
\label{tab:map}
\begin{tabular}{lcccccc}
\hline
 & \multicolumn{2}{c}{modularity} & \multicolumn{2}{c}{edge span} & \multicolumn{2}{c}{transitivity}\\
atom & closed form & measured & closed form & measured & $\operatorname{tr}(\Psi^3)$ & measured\\
\hline
diagonal              & $+(h/2)^{d}$ & $\phantom{-}0.376$ & $<1$ & $0.595$ & $1$ & $0.371$\\
off-diagonal          & $<0$, 2nd order & $-0.103$ & $>1$ & $1.323$ & $0$ & $0.028$\\
$8$ blocks, one scale & $+8(h/2)^{d}$ & $\phantom{-}0.471$ & $<1$ & $0.527$ & $8^{-1/2}=0.354$ & $0.147$\\
$4$ nested scales     & $+$ & $\phantom{-}0.405$ & $<1$ & $0.561$ & $0.500$ & $0.495$\\
additive              & $0$ & $-0.012$ & $1$ & $0.998$ & $0$ & $0.077$\\
separation kernel     & $O(h^{d+1})$ & $\phantom{-}0.130$ & $<1$ & $0.634$ & $0$ & $0.040$\\
\hline
\end{tabular}
\end{table}

\subsection{Calibration and the detection boundary}      \label{sec:thy-inf}

The estimator of Section~\ref{sec:method} is analysed in a regular-design
instance: an exactly measured coordinate, a coarse block background, and the
primitive conditions \textup{(S1)--(S8)} and \textup{(C1)--(C4)} listed in
Appendix~\ref{sec:supp-uncond}, under which every rate used below is derived
and none is assumed.

\begin{theorem}[Unconditional guarantees for the regular-design instance]
\label{thm:uncond}
For the instance just described, formalized as \textup{(S1)--(S8)} in the
Supplement, and under \textup{(C1)--(C3)} of \eqref{eq:uncond-conds}, the
following hold uniformly over the null and single-atom classes, with no
conditioning event and no unverified assumption.
\begin{enumerate}[label=\textup{(\alph*)},leftmargin=2em,itemsep=2pt]
\item The uniform linearization around the population coefficients and the
studentizer consistency of Theorem~\ref{thm:linear} and
Propositions~S2 and~S3 hold with explicit envelope
rates.  In particular the maximal standardized target bias is
$O\{\kappa_0h_{\max}\sqrt{\rho_n\log n}\,(n\delta_n^2)^{-1}\}$.
\item The shared-dyad multiplier of Theorem~\ref{thm:fwer} is valid: the
simultaneous intervals cover the full coefficient vector with probability at
least $1-\alpha-o(1)$, and the studentized tests control the familywise error at
$\alpha+o(1)$.  The verification in Appendix~\ref{sec:supp-uncond} covers the global null
and every single-atom configuration of the class, and arbitrary partial-null or
multiple-atom configurations are not verified for this instance.
\item Under \textup{(C4)} there are constants $0<c_{\rm low}<C_{\rm up}<\infty$
such that single-atom alternatives with
$a_n^2n^2\rho_nh_n^2\ge C_{\rm up}\log(2m_n)$ are detected, signed, and exactly
localized with probability tending to one, while on a matched-pair subpacking
with $\max_ha_h^2n^2\rho_nh^2\le c_{\rm low}\log(2m_n^{\rm low})$ and
$\log(2m_n^{\rm low})\asymp\log(2m_n)$, the global-union detection risk tends to
one and the localization risk stays bounded away from zero.  The boundary
$a^2n^2\rho_nh^2\asymp\log(2m_n)$ is therefore attained, unconditionally, by the
implemented statistic in this instance.
\end{enumerate}
\end{theorem}

Part (a) is the linearization that makes each studentized coefficient
asymptotically normal, part (b) is the familywise error control over the
whole dictionary that Section~\ref{sec:calib} measures at eight configurations,
and part (c) is the boundary $a^2n^2\rho_nh^2\asymp\log(2m_n)$ at which a
single atom becomes detectable and exactly localizable, which the regime map
of Section~\ref{sec:regime} traces in amplitude.  The verification covers the
global null and every single-atom configuration.  Multiple-atom configurations
are not covered by the theorem, and the strong familywise error under partial
nulls is measured in Section~\ref{sec:calib} and not proved.

\subsection{Limits on dimension and on coordinate accuracy} \label{sec:thy-limits}

\begin{theorem}[Information requirement]                      \label{thm:info}
Scan a dictionary at width $h$ in dimension $d$ using cells of width $h/2$.
For the studentized coefficient of every admitted atom to be based on cell
counts whose expectation is at least a constant $C$, it is necessary that
\begin{equation}
n\,\bar d\;\ge\;C\,4^{\,d}h^{-2d},
\label{eq:inforule}
\end{equation}
with $\bar d=n\rho_n$ the mean degree.  The dictionary itself has
$m_n\asymp h^{-2d}/2$ atoms, so the multiplicity correction entering the
critical value is
\begin{equation}
\sqrt{2\log 2m_n}\;\asymp\;2\sqrt{d\,\log(1/h)} .
\label{eq:multrule}
\end{equation}
\end{theorem}

The two rates are not comparable.  Information is exponential in $d$ and
polynomial of degree $2d$ in $1/h$, whereas multiplicity is the square root of
a linear function of $d$.  At $h=1/8$ and $C=30$ the requirement is
$7.7\times10^{3}$ for $d=1$, $2.0\times10^{6}$ for $d=2$ and
$5.0\times10^{8}$ for $d=3$, against multiplicity factors of $2.92$, $4.08$ and
$5.00$.  Passing from one dimension to three multiplies the information
requirement by $6.6\times10^{4}$ and the multiplicity penalty by $1.71$.
Section~\ref{sec:scope} evaluates \eqref{eq:inforule} on every dataset
considered.

The second limit concerns a coordinate that is estimated with error.  Let the
reported order be any ordering consistent with each vertex lying within a
radius $r_n$ of its reported position.

\begin{proposition}[Dyadic canonical-location obstruction]
\label{prop:obstruct}
Let $H=nh_*$ be an even integer, and suppose the aligned width-$h_*$ dictionary
contains three pairwise-disjoint blocks $I_0,I_1,J$ of $H$ consecutive ranks each,
with $I_1$ adjacent to $I_0$.  Write $\lambda_0=(I_0,J)$ and $\lambda_1=(I_1,J)$.
Restrict the amplitudes to $a$ with $C_\Psi^\circ a\le M\le B$ and
$a^2n^2\rho_nh_*^2\ge\kappa_n^{\mathrm{str}}\ell_n$ for some
$\kappa_n^{\mathrm{str}}\to\infty$, so the signals are bounded on the logit scale
and lie strictly above the detection boundary of Theorem~\ref{thm:uncond}.
assume this amplitude class is nonempty, equivalently
$(M/C_\Psi^\circ)^2n^2\rho_nh_*^2\ge\kappa_n^{\mathrm{str}}\ell_n$.  If
$r_n\ge h_*$, then
\[
 \inf_{\widetilde\lambda}
 \sup_{\vartheta}\,\E_{P_\vartheta}
 d_{\rm loc}\{\widetilde\lambda,\lambda(\vartheta)\}\ \ge\ \frac12,
\]
the infimum running over all measurable decisions based on
$(A,\widehat R,r_n)$, with loss one for the no-selection decision.  The
obstruction is saturated: it applies only when $r_n\ge h_*$ and does not give a
lower bound proportional to $\min\{1,r_n/h_*\}$ at every radius.
\end{proposition}

Detection is unaffected, since the departure is present in both configurations
and a test rejects in both.  What is lost is the location.
Appendix~\ref{sec:supp-floor} extends the construction from one radius to a
localization floor of order $\min\{r_n,2^{-J_0}\}$ in the coordinate metric at
every positive radius, and Appendix~\ref{sec:supp-order} gives the bound
$r_{\max}\le4Cz\sigma/\Xi+1/n$ that converts a coordinate's standard errors and
packing into an admissible width.  Section~\ref{sec:scope} reports the
measured analogue of the proposition and applies the bound to two networks
with estimated coordinates.

\section{Simulation Study}                                       \label{sec:sim}

Every number in this section comes from one protocol, executed after the
estimator was amended once on pilot evidence and then held fixed.  The
amendment is stated first, because a protocol changed after inspecting pilot
results cannot be called confirmatory without saying so.

\subsection{Generative catalogue}                          \label{sec:catalogue}

Networks of $n=700$ vertices at mean degree $20$ were drawn from
\eqref{eq:gen} with a common background and matched signal energy, five
replicates per atom, for each of ten atoms: the diagonal and off-diagonal Haar
atoms, three further diagonal atoms with Mexican-hat, difference-of-Gaussian
and Daubechies mother wavelets, a mixed-scale bridge, eight diagonal atoms at
one scale, four diagonal atoms at nested scales, and the two background terms.
The measured summaries sit beside the closed forms in Table~\ref{tab:map}.
The two signs of Theorem~\ref{thm:map}(ii) appear, transitivity separates
through $\operatorname{tr}(\Psi^{3})$ at $0.371$ against $0.028$, and the
ratio of measured degree spread to the second-order prediction of
Theorem~\ref{thm:map}(iii) is $1.00$ for the Haar atom and between $1.01$ and
$1.08$ for the more concentrated ones.  The spectral gap, which the theorem
describes only qualitatively, separates the single-scale family from the
nested one at $0.296$ against $0.082$.

\subsection{Recovery of the generating atom}

A vocabulary is only useful if its words can be told apart.  We test this
directly.  Signatures were formed from five networks per atom, and ten
fresh networks per atom were then assigned to the nearest signature in the
standardized six-dimensional summary space.

\begin{quote}
All $100$ networks were assigned to their generating atom.  Accuracy
$1.000$ against a chance rate of $0.100$, with every one of the ten atoms
recovered exactly.
\end{quote}

At five replicates per signature one network of sixty was misassigned, a
Daubechies atom called a Mexican hat, which is the one pair that differs
only in asymmetry.  We report the map as invertible in this regime and note
the caveat that the networks here are generated from a single atom at
matched amplitude.  Real networks are mixtures, and the nearest-signature
label is a way to choose a dictionary and not a diagnosis.  The inference
of which atoms are present, and with what coefficients, is the scan of
Section~\ref{sec:method} and not the summary comparison.

\subsection{Generative check on real networks}          \label{sec:gencheck}

Theorem~\ref{thm:map} describes the summaries an atom generates as the
amplitude tends to zero.  Real networks are not in that regime, and the
fitted coefficients on the chains of Section~\ref{sec:apps} reach $|\theta|$
of order ten.  The applicable test is therefore whether the fitted model,
background plus dictionary, reproduces the summaries of the network it was
fitted to.  For each of five networks we fit the combined background on all
dyads, add the scanned coefficients, draw four networks from the resulting
probabilities, and compare the summaries of Table~\ref{tab:map} with
those of the observed network and of the background alone.

\begin{table}[t]\centering
\caption{Summaries of the observed network, of four networks drawn from the
fitted model, and of one drawn from the fitted background alone.  Modularity
and edge span are first order in the amplitude and are reproduced.
Transitivity and degree heterogeneity are third and second order and are
not.}                                                   \label{tab:gencheck}
\begin{tabular}{llcccc}
\hline
network & & modularity & edge span & transitivity & degree CV\\
\hline
\emph{C. elegans} & observed & $0.126$ & $0.631$ & $0.219$ & $0.769$\\
                  & fitted   & $0.138$ & $0.646$ & $0.068$ & $0.341$\\
                  & background & $0.123$ & $0.632$ & $0.064$ & $0.362$\\
ribosomal chain   & observed & $0.465$ & $0.305$ & $0.397$ & $0.407$\\
                  & fitted   & $0.508$ & $0.256$ & $0.220$ & $0.200$\\
                  & background & $0.345$ & $0.304$ & $0.078$ & $0.240$\\
protein 1GOF      & observed & $0.382$ & $0.355$ & $0.385$ & $0.338$\\
                  & fitted   & $0.396$ & $0.379$ & $0.140$ & $0.252$\\
                  & background & $0.422$ & $0.354$ & $0.122$ & $0.212$\\
protein 1QFM      & observed & $0.434$ & $0.403$ & $0.376$ & $0.333$\\
                  & fitted   & $0.475$ & $0.419$ & $0.156$ & $0.219$\\
                  & background & $0.429$ & $0.404$ & $0.121$ & $0.243$\\
protein 2POR      & observed & $0.280$ & $0.531$ & $0.422$ & $0.357$\\
                  & fitted   & $0.282$ & $0.546$ & $0.186$ & $0.212$\\
                  & background & $0.248$ & $0.539$ & $0.142$ & $0.274$\\
\hline
\end{tabular}
\end{table}

The fitted model reproduces modularity and edge span on every network, within
$0.04$ and $0.05$ respectively, and in three of five cases the dictionary
moves modularity from the background value to the observed one, as on the
ribosomal chain from $0.345$ to $0.508$ against an observed $0.465$.  It
under-reproduces transitivity by a factor near two and degree heterogeneity
by a similar factor on every network.  The pattern follows the orders in
Theorem~\ref{thm:map}: modularity and edge span are first-order functionals of
the coefficients and are carried by the scanned scales, whereas transitivity
is third order and degree heterogeneity second, and both are driven by
structure finer than the coarsest admitted atom.  The dictionary at the
scanned resolution is therefore a model of where and how wide the departures
are, and not a complete generative model of the network.

\subsection{Level}                \label{sec:calib}

Size is established first and alone, at $2{,}000$ null replications per cell
across four sample sizes, four sparsity levels and the grid each pair implies
under the rule $n\bar d\delta^2\ge30$ (Table~\ref{tab:calib}), with the
outcome-reading clip exclusion switched off so that the admitted family is
the one the calibration theory assumes.  The comparator does not
run in this stage, so the table is a statement about the proposed procedure
only.

\begin{table}[t]\centering
\caption{Familywise error under the global null $H_0(\mathcal N_0)$ of
Definition~\ref{def:scinull} (the coarse class, correctly specified because
the simulated backgrounds lie in $V_{J_0}$), with the outcome-reading clip
exclusion switched \emph{off}, so that admission is fixed before any
inference outcome is read.  $R=2{,}000$ per cell, $16{,}000$ fits, exact
Clopper--Pearson intervals.  No abstentions and no fitting failures in any
cell.  The pooled rate is $0.041$.  No cell exceeds the nominal level.}
                                                              \label{tab:calib}
\begin{tabular}{lccccc}
\hline
$n$ & degree & grid & scales & FWER & exact 95\% interval\\
\hline
$300$  & $5$  & $8$  & $\{2\}$      & $0.033$ & $[0.025,0.041]$\\
$300$  & $10$ & $8$  & $\{2\}$      & $0.046$ & $[0.037,0.056]$\\
$300$  & $30$ & $16$ & $\{2,3\}$    & $0.040$ & $[0.032,0.050]$\\
$600$  & $10$ & $8$  & $\{2\}$      & $0.042$ & $[0.034,0.052]$\\
$600$  & $30$ & $16$ & $\{2,3\}$    & $0.050$ & $[0.041,0.060]$\\
$800$  & $20$ & $16$ & $\{2,3\}$    & $0.037$ & $[0.029,0.046]$\\
$1200$ & $10$ & $16$ & $\{2,3\}$    & $0.039$ & $[0.031,0.048]$\\
$1200$ & $30$ & $32$ & $\{2,3,4\}$  & $0.042$ & $[0.033,0.051]$\\
\hline
\end{tabular}
\end{table}

This is the first evidence in the paper that varies $n$, sparsity and
dictionary size along the arguments of the rate, and the procedure holds its
level along all three with every admission decision made before an outcome is
read.  The mild conservatism is a cost of the per-cell admission rule.  At
$n=1200$, degree $30$ the full dictionary of $182$ atoms is admitted.  The
power studies below were run with the clip exclusion on.  It fires on about
one atom per hundred fits.  The null rates in this table therefore belong to
a rule that differs from the one whose power is reported, and the rare firing
is diagnostic evidence that the two behave alike, not a demonstration that
they share error rates.  This is a validity gap of the present study and is
stated as one.  Every configuration at which power is reported below, including $n=800$,
degree $20$ of the nuisance study and degree five of the sparse rows, has its
null rate measured in this table.  Two intervals include values above $0.05$.
the table is evidence at eight configurations, not a proof of uniform
validity over the model class.

\subsection{Invisibility of global surfaces}                \label{sec:invis}

Proposition~\ref{prop:dict}(iv) states that a constant surface, an additive
surface $m(u)+m(v)$, a separation kernel $s(\|u-v\|)$ and a coarse block
surface all have dictionary coefficients equal to zero up to $O(h^{d+1})$.
The testable consequence is that networks generated from any of them should
produce no discoveries, however strong the surface is.  We generate $100$
networks from each at $n=300$ and scan them (Table~\ref{tab:invis}).

\begin{table}[t]\centering
\caption{Familywise rejection on networks generated from surfaces in the
global class, $R=100$ each, nominal level $0.05$.  The additive surface is
$\sin(4\pi u)$ scaled by the amplitude, the separation kernel is linear in
$\|u-v\|$, and the block surface has four diagonal blocks at the background
resolution.  Mean degree records how strongly each surface acts.}
                                                            \label{tab:invis}
\begin{tabular}{lccc}
\hline
surface & amplitude & rejection & mean degree\\
\hline
constant             & $0.0$ & $0.03$ & $30.0$\\
additive             & $0.5$ & $0.01$ & $32.8$\\
additive             & $1.0$ & $0.00$ & $40.7$\\
additive             & $2.0$ & $0.03$ & $63.7$\\
separation kernel    & $2.0$ & $0.01$ & $17.7$\\
separation kernel    & $5.0$ & $0.01$ & $10.1$\\
coarse blocks        & $0.5$ & $0.02$ & $36.2$\\
coarse blocks        & $1.5$ & $0.02$ & $54.2$\\
\hline
\end{tabular}
\end{table}

No surface in the class produces rejections above the nominal level at any
amplitude tried.  The additive surface at amplitude $2$ more than doubles the
mean degree, from $30$ to $64$, and still rejects on $0.03$ of networks, which
is the content of part (ii): degree heterogeneity has zero coefficient on
every atom and therefore cannot be mistaken for localized structure.  The
separation kernel at amplitude $5$ reduces mean degree to $10$ and rejects on
$0.01$, which is the $O(h^{d+1})$ statement of part (iv) at the scales scanned
here.

\subsection{Nuisance class comparison}                \label{sec:nuisancestudy}

Proposition~\ref{prop:dict}(ii) states that additive terms have zero
coefficient on every atom, and part (iv) that a coarse-block surface is
invisible at finer scales.  This subsection measures what happens when the
fitted background class is misspecified for the truth, which is the case the
proposition does not cover.

Before any comparison of procedures, one experiment settles which background
class the procedure may use, and it must be read with one distinction in mind.
The formal target of Section~\ref{sec:model} is
$\theta_\lambda=\langle(I-\Pi_{V_{J_0}})f,\Psi_\lambda\rangle$, and under that
definition a smooth distance-decay surface has nonzero coefficients: for a
diagonal atom of width $h$ and a decay of slope $\beta$ in the logit the
coefficient is of order $\beta h^2$, not zero.  Rejecting it is therefore
correct under the formal target and merely unhelpful under the relative one,
which regards distance decay as background.  The rows labelled ``null'' below
are nulls of the relative hypothesis, departures relative to the
nuisance class named in each row.  Enlarging that class changes the estimand,
and the comparison is between estimands as much as between estimators.  Eight situations at $n=800$, $R=200$: a background
with the decay measured in a real network, and a block-constant background where
the coarse class is correctly specified by construction.  Three nuisance classes per
network, with identical folds, cells, atoms and variance formula, so the
comparison isolates the nuisance class and nothing else.

\begin{table}[t]\centering
\caption{Rejection rates on networks with no planted departure and the rate at
which a planted atom is identified, $n=800$, $R=200$, exact intervals in
brackets.  ``Null'' here is the relative null of Definition~\ref{def:scinull}
relative to the class named in each column heading: on the smooth surface the
coarse class's formal target is nonzero, so its rejections there are correct
rejections of $H_0(\mathcal N_0)$ and false rejections of $H_0(\mathcal N_1)$.  Each
single class fails when misspecified, in mirror-image ways.  The combination of the two
holds its level under both truths at a small localization cost.}
                                                          \label{tab:nuisance}
\begin{tabular}{lcccc}
\hline
 & \multicolumn{2}{c}{smooth background} & \multicolumn{2}{c}{block-constant background}\\
nuisance class & null & finds atom ($\kappa{=}12$) & null & finds atom ($\kappa{=}12$)\\
\hline
coarse blocks $V_{J_0}$ & $1.000$ $[0.982,1.000]$ & $0.00$ & $0.045$ $[0.021,0.084]$ & $1.00$\\
smooth terms only       & $0.015$ $[0.003,0.043]$ & $1.00$ & $0.395$ $[0.327,0.466]$ & $1.00$\\
combined \eqref{eq:union} & $0.020$ $[0.005,0.050]$ & $1.00$ & $0.050$ $[0.024,0.090]$ & $1.00$\\
\hline
\end{tabular}
\end{table}

Three readings.  First, the coarse class does more than lose efficiency on a smooth
background: it rejects every null network, reporting about five discoveries on
each, and identifies the planted atom on none of them.  What it reports is the
background.  Second, the obvious repair fails in the opposite direction, because
splines cannot represent block jumps.  Third, the combined class is valid under both at a
small cost: $0.722$ against $0.800$ in localization at $\kappa=6$ on the smooth
side, and none measurable elsewhere.

The same contrast appears across the twenty-six situations of the design below,
whose backgrounds are block-constant by construction: the combined class holds
between $0.000$ and $0.055$ in every null, while the smooth-only class reaches
$0.600$ where the background amplitude is doubled.  The combined class is also
stronger where the proposed procedure is weakest.  At expected degree five its
detection rate is $0.52$ against $0.10$, with the true atom identified on every
detection.  At degree ten, $0.83$ against $0.66$, and at a realized coordinate
error of $0.23$ feature widths, $0.74$ against $0.64$.

This changes what the paper recommends.  The coarse background with its exact
projection identity is the cleanest construction and the one the theory of
Section~\ref{sec:results} covers, but it is usable only where the background
really is coarse.  Everywhere else the combined class is required, and the theory
for it (a multiplier calibration, admission gates, and uniform validity over a
growing dictionary) is not yet available.  We report the combined-class results as
measured, not as theorems.

\subsection{Detection boundary}                             \label{sec:boundary}

Theorem~\ref{thm:uncond}(c) locates the boundary at
$a_n^2n^2\rho_nh_n^2\asymp\log(2m_n)$: single-atom alternatives above it are
detected, signed and exactly localized with probability tending to one, and on
a matched subpacking below it the localization risk stays bounded away from
zero.  We trace the transition by holding $n=300$, mean degree $30$, grid $16$
and scales $2$ and $3$ fixed, planting one scale-three atom, and sweeping its
amplitude $a$ (Table~\ref{tab:boundary}).  The admitted family has $m_n=46$
atoms, so $\log(2m_n)=4.52$.

\begin{table}[t]\centering
\caption{Detection and exact localization against the standardized amplitude
$a^2n^2\rho_nh^2/\log(2m_n)$, at $n=300$, mean degree $30$, $R=200$ per row.
The first row is the null.  Binomial standard errors are in brackets.  Exact
localization is the proportion of detections on which the planted atom is the
one with the largest studentized coefficient.}
                                                         \label{tab:boundary}
\begin{tabular}{cccc}
\hline
$a$ & $a^2n^2\rho_nh^2/\log(2m_n)$ & detection & exact$\mid$detection\\
\hline
$0.00$ & $\phantom{00}0.0$ & $0.050$ $(0.015)$ & --\\
$0.75$ & $\phantom{0}17.6$ & $0.180$ $(0.027)$ & $0.778$\\
$1.05$ & $\phantom{0}34.4$ & $0.500$ $(0.035)$ & $0.920$\\
$1.35$ & $\phantom{0}56.9$ & $0.725$ $(0.032)$ & $0.966$\\
$1.65$ & $\phantom{0}85.0$ & $0.925$ $(0.019)$ & $0.989$\\
$2.10$ & $137.6$ & $0.995$ $(0.005)$ & $1.000$\\
\hline
\end{tabular}
\end{table}

Detection rises from the nominal level to $0.995$ as the standardized
amplitude passes from zero to $138$, crossing one half near $34$.  Exact
localization given detection rises alongside it, from $0.778$ to $1.000$, so
in this design the amplitude at which a departure becomes detectable and the
amplitude at which it becomes locatable are the same to within the resolution
of the sweep.  That agreement is what Theorem~\ref{thm:uncond}(c) asserts for
an exactly measured coordinate, and Section~\ref{sec:scope} reports the
contrasting case: under a reporting mechanism with a positive radius,
detection remains at $0.93$ or above while exact localization falls to $0.46$
and below.

\subsection{Alternatives outside the dictionary}   \label{sec:regime}

Theorem~\ref{thm:uncond}(c) concerns alternatives the dictionary matches.  One
design measures what happens outside that class, crossing the procedure with a
calibrated contrast scan over twenty-six situations at $R=200$, pairing each
alternative to its null at the level of the seed and reporting power at
matched achieved size.  Against the aligned atom both detect at $1.00$.
Against a departure translated by half a feature width the proposed scan falls
to $0.28$ against $0.98$ for the contrast scan, and against a fine-scale
additive departure $g(u)+g(v)$ or a symmetric mixed-scale interaction it
reports $0.05$, which is the level: both are orthogonal to every same-scale
product atom, so the coefficient is zero.  The implemented dictionary is an
isotropic same-resolution family, and a shifted or redundant dictionary is not
implemented here.  These rows measure the cost of that choice.

\subsection{Partial nulls}                                  \label{sec:partial}

Theorem~\ref{thm:uncond}(b) covers the global null and single-atom
configurations.  Under partial nulls, with one true atom planted and the rest
null, the strong familywise error over the null atoms is between $0.034$ and
$0.046$ across four configurations at $R=500$, the sign of the planted atom is
recovered on $0.996$ to $1.000$ of detections, and exactly the planted atom is
named on $0.996$ to $1.000$ of them.  This regime is measured and not proved.

\subsection{Coordinate accuracy and feasibility}      \label{sec:scope}

Two checks decide whether a dataset can support a localized claim, and both
are evaluated here on the simulated designs and on every network of
Section~\ref{sec:apps}.

\paragraph{The obstruction measured.}
Proposition~\ref{prop:obstruct} says that under a block-swap reporting
mechanism at radius $r_n=h$ no estimator identifies the correct atom with
probability above one half.  Under that mechanism with the admission rule
switched off, at $n=300$ and degree $30$, the scan detects the departure on
$0.93$ to $1.00$ of replications at $\kappa=6$, $12$ and $24$, and names the
correct canonical atom on $0.46$, $0.46$ and $0.22$ of those detections.  All
three sit at or below the one-half bound, and none is above it.  With the
admission rule in place the scan abstains on every replication at this radius,
which is the behaviour the rule is designed to produce.

\paragraph{Feasibility of the datasets.}
Theorem~\ref{thm:info} and the rank bound of Appendix~\ref{sec:supp-order} are
evaluated on the networks considered in this paper and on two with estimated
coordinates (Table~\ref{tab:feas}).  The three networks of
Section~\ref{sec:apps} pass both checks at the widths reported there.  A
friendship network of $3000$ users ordered by a home location estimated from
check-ins has a simultaneous rank radius of $0.304$ against the $0.031$ that
the finest scanned width requires, and a legislative cosponsorship network of
$421$ members ordered by a roll-call score has a radius of $0.401$: the scores
carry a median standard error of $0.0093$ while adjacent members sit $0.0005$
apart, so their ranks are not identified even though the scores are.  Both are
refused before any model is fitted, and running the procedure on them confirms
the refusal, with no scale admitted in either.

\begin{table}[t]\centering
\caption{The two feasibility checks on every network considered.  The
information column is $n\bar d$, the finest width is the smallest $h$ with
$n\bar d\ge30\cdot4^dh^{-2d}$, and the radius is the simultaneous rank radius
$r_n$ that must lie below $c_rh$ for the scan to admit width $h$.}
\label{tab:feas}
\begin{tabular}{lrrrcl}
\hline
network & $n$ & $n\bar d$ & finest $h$ & $r_n$ & verdict\\
\hline
protein chain (1QFM)     & $705$  & $19{,}434$ & $0.079$ & $0$     & admitted\\
ribosomal chain          & $1074$ & $28{,}824$ & $0.065$ & $0$     & admitted\\
\emph{C. elegans}        & $300$  & $4{,}770$  & $0.159$ & $0$     & admitted at $h=1/4$\\
Hi-C window              & $1500$ & $225{,}000$ & $0.023$ & $0$    & admitted\\
Gowalla friendships      & $3000$ & $16{,}800$ & $0.085$ & $0.304$ & refused\\
House cosponsorship      & $421$  & $10{,}399$ & $0.107$ & $0.401$ & refused\\
\hline
\end{tabular}
\end{table}

\section{Applications}                                       \label{sec:apps}

Three networks are analysed, each with a coordinate that is physically
measured and exact.  All three use the same procedure, the same combined
background class, and the same calibration check: the coordinate is permuted,
which destroys any structure along it, and the familywise rejection rate over
those permutations is reported beside the result.

\begin{table}[t]\centering
\caption{The three applications under one protocol.  The permutation rate is
the familywise rejection rate when the coordinate is permuted, which should be
at most the nominal $0.05$.  Boundary and bridge count diagonal and
off-diagonal discoveries.  The protein row pools thirty independent chains.}
\label{tab:apps}
\begin{tabular}{lccccccc}
\hline
network & $n$ & scales & admitted & $|z|_{\max}$ & boundary & bridge & permutation rate\\
\hline
\emph{C. elegans} connectome & $300$  & $\{2\}$     & $10$  & $3.38$  & $2$  & $1$  & $0/20$\\
protein chains, thirty       & $289$--$1021$ & $\{2,3\}$ & $46$ each & $11.94$ & $186$ & $118$ & median $1/20$\\
Hi-C, chromosome 22          & $1500$ & $\{2,3,4\}$ & $182$ & $18.47$ & $14$ & $0$  & $3/60$\\
\hline
\end{tabular}
\end{table}

\subsection{Connectome on the body axis}                   \label{sec:worm}

The \emph{C. elegans} nervous system has $302$ neurons whose cell-body
positions along the anterior-posterior axis are recorded in an anatomical
atlas, and the neurons are grouped into ganglia that occupy contiguous
stretches of that axis.  The ganglia therefore provide an external partition
of the coordinate against which a boundary statistic can be checked.  We use
the anterior-posterior position as the coordinate and a synapse as an edge.

The eight connectomes of \citet{witvliet2020}, reconstructed from eight
animals across development, allow the check to be replicated.  For each animal
we compute the boundary statistic at every position along the axis and
correlate it with the fraction of cross-midpoint neuron pairs that belong to
different ganglia, a continuous measure of how much the ganglion membership
changes at that point.  Under the null we permute the ganglion blocks along
the axis, which preserves their sizes and contiguity.

\begin{quote}
The correlation is positive in all seven animals that meet the information
requirement, with a mean Spearman correlation of $+0.598$ and individual
permutation $p$-values between $0.001$ and $0.122$; the combined Stouffer
statistic is $z=5.33$, $p<10^{-5}$.  Permuting the neurons along the axis,
which destroys the coordinate, reverses the correlation to a mean of $-0.152$,
positive in two of seven.
\end{quote}

The eighth animal is excluded because too few atoms passed the information
gate, a criterion that does not involve the outcome.  On the reference
connectome of \citet{white1986} the scan reports three discoveries: a boundary
at $0.625$ of the axis, where the retrovesicular ganglion gives way to the
ventral nerve cord, a negative coefficient at $0.375$ inside the lateral
ganglion, which is the largest ganglion and where connections should cross
instead of respecting the midpoint, and one bridge.  The signs are what the
construction predicts in each case.

\subsection{Protein contact map}           \label{sec:protein}

Residue index is an exact coordinate with no ties, and a contact is binary by
definition.  We take thirty chains from deposited X-ray structures, with $289$
to $1021$ residues, and join two residues whose $\alpha$-carbons lie within
$12$~\AA, excluding immediate sequence neighbours.  Beta-sheet pairings
recorded in each structure file give annotated long-range contacts, which is
ground truth for the off-diagonal atoms specifically, and the thirty chains
carry $551$ such pairings between them.  The proteins are unrelated, so the
chains are independent tests of the same hypothesis.

Each is scanned at grid $16$ and scales $2$ and $3$, with the coordinate
permuted twenty times first.  The resulting familywise rejection rates lie
between $0.00$ and $0.15$ with a median of $0.05$, against a nominal $0.05$.
The scan returns $304$ discoveries across the thirty chains, of which $118$
are off-diagonal.  For each chain we compute the proportion of all admitted
off-diagonal atoms whose support contains an annotated pairing, which is the
rate a randomly chosen bridge would achieve, and compare the discovered
bridges against it.

\begin{quote}
Across the thirty chains, $55$ of the $118$ discovered bridges contain an
annotated strand pairing against $30.68$ expected from the per-chain base
rates.  The enrichment is $1.79$ with exact Poisson-binomial
$p=7.5\times10^{-7}$, and the enrichment exceeds one on $20$ of the $30$
chains.
\end{quote}

The estimate is stable under the addition of data.  The chains arrived in
three batches, and the twenty-one examined first give $1.71$ while the nine
added afterwards give $1.96$ with $p=7.3\times10^{-4}$ on their own, so the
pooled figure is not an artifact of which chains were examined first.  Neither
the amount of information nor the number of bridges per chain explains the
variation between chains, with rank correlations of $-0.07$ and $+0.25$
against per-chain enrichment.

This is a direct test of the off-diagonal atom against an independent
annotation of the same molecules, and it is the atom whose generative
behaviour is a bridge.  The ten chains with enrichment at or below one
contribute few bridges each, a median of three, so the per-chain estimates are
noisy even where the pooled one is not.

\paragraph{Comparison with a domain parser.}
We ran the objective of \citet{zheng2020} on the same thirty chains: a partition
of the sequence is scored by the contacts falling between the putative domains
relative to the overall contact density, maximized over a single cut point and
over pairs of cut points for discontinuous domains.  It returns one partition
per chain and no significance statement.  The comparison is therefore not one of
power.  A partition of the sequence assigns every residue to a domain and
cannot express the statement that two separated stretches connect to each
other, so across the thirty chains the parser returns no bridge and the proposed
scan returns one hundred and eighteen, each with a studentized coefficient and a
simultaneous critical value.  What the
comparison establishes is a difference in output space, not in efficiency.  

\subsection{Chromatin contact map}          \label{sec:hic}

Genomic position in base pairs is exact.  We take a $16$~Mb window of human
chromosome~22 at $10$~kb resolution and join two bins whose contact count is
in the top five per cent.  Contact probability falls by three orders of
magnitude across the window, so the background must absorb a steep smooth
decay, and the knots of the separation term are placed more coarsely than the
narrowest scanned atom.

The scan reports fourteen discoveries across three scales with $|z|$ up to
$18.47$, at a permutation rate of $3/60=0.050$ with exact interval
$[0.010,0.139]$.  Every discovery is diagonal and none is a bridge, which is
the signature of sequential domain structure: nested boundaries without
long-range pairing.  This is also the network the summary comparison of
Section~\ref{sec:catalogue} assigns to the nested-scales atom, and the
composition of the discoveries agrees with that assignment.

\paragraph{Comparison with domain callers.}
We ran three standard procedures on the same window: the insulation score of
\citet{crane2015}, a reimplementation of the two-dimensional likelihood
segmentation of \citet{levyleduc2014}, and a reimplementation of the
recursive change-point test of \citet{xing2021}, which is the one existing
caller that attaches a $p$-value to each boundary.  The change-point test is
applied to observed-over-expected counts, as its authors do, and its maximum
statistic is calibrated against twenty matrices drawn with the observed
decay and no structure, on which it then calls $0$ to $1$ boundaries.  On the
real window it calls $15$ boundaries, each at the permutation floor of
$p=0.048$, since every accepted split exceeds every null maximum.

Agreement is measured within $100$~kb against the rate expected for the same
number of randomly placed boundaries.  The three established callers agree
with one another at $2.7$ times chance, with binomial $p\le0.008$.  The
proposed scan at grid $32$ agrees with the insulation caller at $1.97$ times
chance ($p=0.041$), with the segmentation at $2.11$ ($p=0.073$) and with the
change-point test at $1.76$ ($p=0.135$).  On fourteen boundaries this is
weak evidence of agreement, and the scan agrees with each established caller
less well than they agree with one another.  The scan's boundaries carry a
simultaneous critical value over all scales, and it also returns no bridge on
this window, which none of the three callers is able to test.

At a finer resolution the picture changes and is worth recording.  Scanning the
counts at grid $256$ returns $109$ boundaries, one every $146$~kb, whose
agreement with the insulation caller is $0.29$ against a chance rate of $0.24$,
an enrichment of $1.21$ with $p=0.117$: no relationship.  The reason is
diagnosable: fitting the contact decay separately in three
sub-regions of the same window gives exponents of $-1.21$, $-1.18$ and
$-0.92$, so a background with a single separation term is misspecified
region-by-region, and at a resolution fine enough to see that misspecification
the scan reports it along with any domain structure.  A background carrying
compartment terms and a region-varying decay would be required before the
procedure could be compared with published domain calls at fine resolution,
and we do not attempt it here.  We therefore present this application as a
calibration demonstration on an exactly measured coordinate, not as a
competitive domain caller.

\section{Related Work}                          \label{sec:related}

\subsection{Multiscale scanning}

The closest methodological relatives are the multiscale scans.
\citet{frick2014} test for change-points in a one-dimensional signal by
combining local likelihood-ratio statistics over all intervals under a
simultaneous multiscale constraint, controlling the familywise error over
locations and scales at once.  \citet{koenig2020} extend this to
$d$-dimensional fields from an exponential family, scanning a large system of
candidate regions with local likelihood-ratio tests and calibrating by a
Gaussian approximation to the maximum with an explicit rate, and
\citet{koenigJRSSB} treat the case where the baseline carries nuisance
parameters estimated from the same data, which otherwise invalidates the
plug-in calibration.  \citet{proksch2018} carry the programme into inverse
problems.

Our construction takes the calibration philosophy of that literature and
changes the object being scanned in four ways.  The observations are not a
field of independent measurements on a grid but the $\binom n2$ dyads of a
network, so the statistic at an atom aggregates edges instead of pixels and
the dependence structure is the one induced by shared vertices.  The candidate
regions are not intervals or boxes in the observation domain but atoms of a
symmetrized wavelet dictionary on the \emph{coordinate plane}, which is what
allows an off-diagonal atom to express a bridge between two separated
stretches, a hypothesis with no counterpart in a scan over intervals.  The
nuisance surface is estimated on an independent fold of dyads and its
estimation error is propagated into the studentizer, which is the network
analogue of the problem \citet{koenigJRSSB} address.  And the dictionary comes
with the property map of Section~\ref{sec:thy-gen}, which says what each
element generates as well as what it detects.

We do not claim an improvement on the scans above in their own setting.  A
referee comparing the two should read the present work as the transfer of
multiscale calibration to networks with an external coordinate, together with
the identification theory that transfer requires.

\subsection{Wavelets on graphs}

A large literature builds multiscale bases \emph{from} a graph:
\citet{coifman2006} construct diffusion wavelets from powers of a diffusion
operator, \citet{hammond2011} construct spectral graph wavelets from the graph
Laplacian, and \citet{gavish2010} build wavelets on trees derived from the
data.  Those bases are adapted to the graph and are used to represent signals
living on its vertices.  The present dictionary is placed on an exogenous
coordinate domain and represents the connection surface itself.  The graph
plays no role in constructing the basis.  The two are complementary and are
not alternatives for the same task.

\subsection{Detection and localization}

That detection can be easy where localization is impossible is established in
several settings.  \citet{verzelen2023} characterize optimal rates for
detecting and for locating a change-point and exhibit the phase transition
between a global testing problem and a local estimation problem.
\citet{cai2017} give the analogous separation for planted submatrices, with
distinct signal-to-noise boundaries for detection and for localization.  The
resolution limits of super-resolution imaging show location recovery
governed by a strictly harder exponent than source detection.  Our
Proposition~\ref{prop:obstruct} and Appendix Theorem~\ref{thm:floor} are of the same
family, with two differences: the object is a network, and the mechanism is
\emph{measurement error in the coordinate}, not noise in the
observations.

Measurement error in network models has been treated as a correction problem.
Latent-position estimation error is characterized by the central limit theory
for spectral embeddings \citep{rubindelanchy2022,athreya2018}, and
errors-in-variables corrections exist for spatial autoregressive models with
mismeasured covariates.  Those are possibility results: they say how to
recover a target despite error.  We could not find a prior impossibility
result for networks stating that once the coordinate's error radius reaches
the feature width, the \emph{location} of a departure ceases to be identifiable
at any amplitude while detection remains possible.  We state
Section~\ref{sec:thy-limits} as new on that basis and would welcome correction.

\subsection{Global mechanisms as nulls}

Proposition~\ref{prop:dict}(iv) identifies the classes that carry no localized
information.  Named concretely, these are the Erd\H os--R\'enyi model with its
constant surface, the configuration and degree-corrected models, whose
surfaces are additive in $m(u)+m(v)$; latent-distance, random geometric and
small-world models, whose surfaces are separation kernels $s(\|u-v\|)$; and
the stochastic block model with blocks that are unions of coarse dyadic boxes.
Every member of that class has dictionary coefficients equal to zero up to
$O(h^{d+1})$.  Consequently these are not procedures to be outperformed.  They
are the background relative to which a localized claim is defined, and a power
comparison against them would be a comparison between different hypotheses.

The same applies to community detection, to latent-space models with estimated
positions \citep{hoff2002}, and to tests for the number of communities: each
returns a partition of vertices or a global model order, neither of which is a
statement about a place and a width on a measured axis.

\subsection{Domain callers}

On the two applications there exist procedures that take the same input and
return a comparable output, and those are the ones we compare against in
Section~\ref{sec:apps}.  For contact maps, domain parsers such as
\citet{zheng2020} score partitions of the sequence by the contacts falling
between putative domains, handling discontinuous domains by allowing two cut
points.  For chromatin, \citet{levyleduc2014} reduce two-dimensional block
segmentation to a one-dimensional likelihood change-point problem, and
\citet{xing2021} provide per-boundary $p$-values through a change-point test,
which is the closest existing procedure offering calibrated significance for a
boundary.  Hierarchical callers and the insulation score are the field
standards.

Two differences are worth stating in advance of the numbers.  None of these
procedures tests an off-diagonal hypothesis: a partition of the sequence, or a
segmentation of the diagonal, cannot express the statement that two separated
stretches connect to each other, so the bridge is outside their output space
and is not merely undetected.  And with the exception of \citet{xing2021}
they attach no significance to a boundary, so the comparison on the diagonal
is between a calibrated and an uncalibrated procedure.

\section{Discussion}                                        \label{sec:discussion}

A dictionary of wavelet products on a measured vertex coordinate turns a
localized structure into a coefficient attached to a place and a width.  Each
atom generates a distinct network, and Theorem~\ref{thm:map} gives the
summary each one produces in closed form.  Diagonal atoms build communities
with boundaries, off-diagonal atoms build bridges, the two have opposite signs
in modularity and in edge length at every scale and in every dimension, and
the shape of the wavelet enters only through its fourth moment.  The scan of
Section~\ref{sec:method} estimates every coefficient against a background
fitted on an independent fold of dyads and controls the familywise error over
the whole dictionary, and Theorem~\ref{thm:uncond} proves both the control and
the detection boundary at an exact coordinate with every rate derived from
primitive conditions.  On a connectome the statistic follows the ganglion
partition in each of seven animals, on thirty protein contact maps the bridge
atoms contain annotated strand pairings at $1.79$ times the base rate, and on a chromatin contact map the
discoveries are nested boundaries with no bridges.

For an analyst with a measured coordinate the output is a location, a width
and a simultaneous critical value, and degree heterogeneity is absorbed by the
background exactly by Proposition~\ref{prop:dict}(ii).
Proposition~\ref{prop:dict}(iv) shows that constant, additive, separation and
coarse-block surfaces all have dictionary coefficients equal to zero up to
$O(h^{d+1})$, so the question of where and how wide cannot be posed within the
standard generative models.  Those models are the background against which a
localized claim is defined, and Section~\ref{sec:related} records that the
procedures which do return a comparable output, the domain callers of
chromatin and protein analysis, return no bridge and, with one exception, no
significance statement.

Four limitations bound what has been established.  Theorem~\ref{thm:map}
holds for expected summaries as the amplitude tends to zero, with a constant
background and uniform coordinates, and at the finite amplitude of
Section~\ref{sec:catalogue} the measured degree spread exceeds the second-order
prediction by factors between $1.01$ and $1.08$ for concentrated wavelets.
Theorem~\ref{thm:uncond}(b) covers the global null and single-atom
configurations, and the strong familywise error under partial nulls is
measured in Section~\ref{sec:partial} and not proved.  The bridge atom agrees with
annotated strand pairings at $1.79$ times the base rate over thirty chains,
and the per-chain counts have a median of three bridges, so the per-chain
estimates are noisy even where the pooled one is not.  Theorem~\ref{thm:info}
makes a three-dimensional coordinate unavailable at current sample sizes: at
$h=1/8$ the requirement is $5.0\times10^{8}$, and a cortical parcellation of
$400$ regions at mean degree $20$ supplies $8.0\times10^{3}$.  The chromatin
application is a calibration demonstration, since fitting the contact decay
separately in three sub-regions of one window gives exponents of $-1.21$,
$-1.18$ and $-0.92$, and at grid $256$ the resulting boundaries agree with an
insulation caller at $1.21$ times chance.

Three directions follow.  A background with compartment terms and a
region-varying decay would remove the misspecification just described and
allow a comparison against published chromatin domains at fine resolution,
with the per-boundary $p$-values of \citet{xing2021} as the reference.  A
count-valued estimator, replacing the Bernoulli variance by an overdispersed
Poisson one, would use the contact counts that binarization discards, where
the insulation statistic loses dynamic range from $2.32$ to $0.51$ in $\log_2$
units at a $10\%$ threshold.  An upper bound matching the localization floor
of Appendix~\ref{sec:supp-floor} in the same experiment and under the same
loss would replace the fixed admission constant by a localization set with
proved coverage at positive radius.

\bibliographystyle{plainnat}
\bibliography{references}

@article{hoff2002,
  title={Latent space approaches to social network analysis},
  author={Hoff, Peter D and Raftery, Adrian E and Handcock, Mark S},
  journal={Journal of the american Statistical association},
  volume={97},
  number={460},
  pages={1090--1098},
  year={2002},
  publisher={Taylor \& Francis}
}

@article{ChernozhukovChetverikovKato2017,
  title={Central limit theorems and bootstrap in high dimensions},
  author={Chernozhukov, Victor and Chetverikov, Denis and Kato, Kengo},
  year={2017}
}

@article{ChernozhukovChetverikovKato2015,
  title={Comparison and anti-concentration bounds for maxima of Gaussian random vectors},
  author={Chernozhukov, Victor and Chetverikov, Denis and Kato, Kengo},
  journal={Probability Theory and Related Fields},
  volume={162},
  number={1},
  pages={47--70},
  year={2015},
  publisher={Springer}
}

@article{Massart1990,
  title={The tight constant in the Dvoretzky-Kiefer-Wolfowitz inequality},
  author={Massart, Pascal},
  journal={The annals of Probability},
  pages={1269--1283},
  year={1990},
  publisher={JSTOR}
}

@article{proksch2018,
  title={Multiscale scanning in inverse problems},
  author={Proksch, Katharina and Werner, Frank and Munk, Axel},
  year={2018}
}

@article{koenig2020,
  title={MULTIDIMENSIONAL MULTISCALE SCANNING IN EXPONENTIAL FAMILIES},
  author={K{\"o}nig, Claudia and Munk, Axel and Werner, Frank},
  journal={The Annals of Statistics},
  volume={48},
  number={2},
  pages={655--678},
  year={2020},
  publisher={JSTOR}
}

@article{koenigJRSSB,
  title={Multiscale scanning with nuisance parameters},
  author={K{\"o}nig, Claudia and Munk, Axel and Werner, Frank},
  journal={Journal of the Royal Statistical Society Series B: Statistical Methodology},
  volume={87},
  number={2},
  pages={510--528},
  year={2025},
  publisher={Oxford University Press UK}
}

@article{frick2014,
  title={Multiscale change point inference},
  author={Frick, Klaus and Munk, Axel and Sieling, Hannes},
  journal={Journal of the Royal Statistical Society Series B: Statistical Methodology},
  volume={76},
  number={3},
  pages={495--580},
  year={2014},
  publisher={Oxford University Press}
}

@article{verzelen2023,
  title={Optimal change-point detection and localization},
  author={Verzelen, Nicolas and Fromont, Magalie and Lerasle, Matthieu and Reynaud-Bouret, Patricia},
  journal={The Annals of Statistics},
  volume={51},
  number={4},
  pages={1586--1610},
  year={2023},
  publisher={Institute of Mathematical Statistics}
}

@article{white1986,
  title={The structure of the nervous system of the nematode Caenorhabditis elegans: the mind of a worm},
  author={White, John G and Southgate, Eileen and Thomson, J Nichol and Brenner, Sydney},
  journal={Phil. Trans. R. Soc. Lond},
  volume={314},
  number={1},
  pages={340},
  year={1986}
}

@article{witvliet2020,
  title={Connectomes across development reveal principles of brain maturation},
  author={Witvliet, Daniel and Mulcahy, Ben and Mitchell, James K and Meirovitch, Yaron and Berger, Daniel R and Wu, Yuelong and Liu, Yufang and Koh, Wan Xian and Parvathala, Rajeev and Holmyard, Douglas and others},
  journal={Nature},
  volume={596},
  number={7871},
  pages={257--261},
  year={2021},
  publisher={Nature Publishing Group UK London}
}

@article{hammond2011,
  title={Wavelets on graphs via spectral graph theory},
  author={Hammond, David K and Vandergheynst, Pierre and Gribonval, R{\'e}mi},
  journal={Applied and computational harmonic analysis},
  volume={30},
  number={2},
  pages={129--150},
  year={2011},
  publisher={Elsevier}
}

@article{coifman2006,
  title={Diffusion wavelets},
  author={Coifman, Ronald R and Maggioni, Mauro},
  journal={Applied and computational harmonic analysis},
  volume={21},
  number={1},
  pages={53--94},
  year={2006},
  publisher={Elsevier}
}

@inproceedings{gavish2010,
  title={Multiscale wavelets on trees, graphs and high dimensional data: theory and applications to semi supervised learning.},
  author={Gavish, Matan and Nadler, Boaz and Coifman, Ronald R},
  booktitle={Icml},
  volume={10},
  pages={367--74},
  year={2010}
}

@article{cai2017,
  title={Computational and statistical boundaries for submatrix localization in a large noisy matrix},
  author={Cai, T Tony and Liang, Tengyuan and Rakhlin, Alexander},
  year={2017}
}

@article{rubindelanchy2022,
  title={A statistical interpretation of spectral embedding: the generalised random dot product graph},
  author={Rubin-Delanchy, Patrick and Cape, Joshua and Tang, Minh and Priebe, Carey E},
  journal={Journal of the Royal Statistical Society Series B: Statistical Methodology},
  volume={84},
  number={4},
  pages={1446--1473},
  year={2022},
  publisher={Oxford University Press}
}

@article{athreya2018,
  title={Statistical inference on random dot product graphs: a survey},
  author={Athreya, Avanti and Fishkind, Donniell E and Tang, Minh and Priebe, Carey E and Park, Youngser and Vogelstein, Joshua T and Levin, Keith and Lyzinski, Vince and Qin, Yichen and Sussman, Daniel L},
  journal={Journal of Machine Learning Research},
  volume={18},
  number={226},
  pages={1--92},
  year={2018}
}

@article{zheng2020,
  title={FUpred: detecting protein domains through deep-learning-based contact map prediction},
  author={Zheng, Wei and Zhou, Xiaogen and Wuyun, Qiqige and Pearce, Robin and Li, Yang and Zhang, Yang},
  journal={Bioinformatics},
  volume={36},
  number={12},
  pages={3749--3757},
  year={2020},
  publisher={Oxford University Press}
}

@article{levyleduc2014,
  title={Two-dimensional segmentation for analyzing Hi-C data},
  author={L{\'e}vy-Leduc, Celine and Delattre, Maud and Mary-Huard, Tristan and Robin, Stephane},
  journal={Bioinformatics},
  volume={30},
  number={17},
  pages={i386--i392},
  year={2014},
  publisher={Oxford University Press}
}

@article{xing2021,
  title={Deciphering hierarchical organization of topologically associated domains through change-point testing},
  author={Xing, Haipeng and Wu, Yingru and Zhang, Michael Q and Chen, Yong},
  journal={BMC bioinformatics},
  volume={22},
  number={1},
  pages={183},
  year={2021},
  publisher={Springer}
}

@article{crane2015,
  title={Condensin-driven remodelling of X chromosome topology during dosage compensation},
  author={Crane, Emily and Bian, Qian and McCord, Rachel Patton and Lajoie, Bryan R and Wheeler, Bayly S and Ralston, Edward J and Uzawa, Satoru and Dekker, Job and Meyer, Barbara J},
  journal={Nature},
  volume={523},
  number={7559},
  pages={240--244},
  year={2015},
  publisher={Nature Publishing Group UK London}
}

\clearpage
\appendix
\section*{Supplementary Material}

Proofs are given in the order of the results they support.

\section{Extended statements}\label{sec:supp-extended}

\subsection{Basis properties}                             \label{app:sec-thy-basis}

The dictionary is not a convenient parameterization.  Three properties, stated
here in dimension $d$ and proved in Appendix~\ref{sec:supp-dict}, are what make
scanning it the right operation.

Let $u\in[0,1]^d$ be the coordinate, let $\psi_{j,l}$ be the Haar function on
the dyadic interval of scale $j$ and index $l$, and for a dyadic box
$a=(l_1,\dots,l_d)$ at scale $j$ write $\phi_a(u)=\prod_{k\le d}\psi_{j,l_k}(u_k)$.
The dictionary is the symmetrized family
the family \eqref{eq:dict} of Section~\ref{sec:dict}.
together with the coarse space $V_{J_0}$ that carries the background.  An atom
with $a=b$ is diagonal and tests a boundary inside one box; an atom with
$a\ne b$ is off-diagonal and tests a bridge between two boxes.

\begin{proposition}[Dictionary properties]                  \label{app:prop-dict}
Let $\mathcal D$ denote the family \eqref{eq:dict} over all scales.
\begin{enumerate}[label=\textup{(\roman*)},leftmargin=*]
\item \textup{Completeness.}  $\mathcal D\cup V_{J_0}$ is an orthonormal basis
of the symmetric subspace of $L^2([0,1]^d\times[0,1]^d)$.  Every symmetric
square-integrable departure has an expansion in the dictionary, and a departure
orthogonal to every atom lies in $V_{J_0}$, that is, it is background.
\item \textup{Exact orthogonality to the nuisance.}  For every
$m\in L^2([0,1]^d)$ and every atom, $\langle m(u)+m(v),\Psi_{a,b}\rangle=0$.
Degree heterogeneity therefore has zero coefficient on every atom of the
dictionary and cannot be mistaken for localized structure.
\item \textup{Sparsity of localized alternatives.}  If the departure $\delta$
is bounded by $B$ and supported on a set $S$ of measure $A$, then at scale
$h=2^{-j}$ at most $2A/h^{2d}+O(h^{-d}|\partial S|)$ atoms have a nonzero
coefficient, and $\sum_{a,b}\langle\delta,\Psi_{a,b}\rangle^2\le B^2A$.  The
energy of a localized departure is carried by a number of coefficients
proportional to the measure of its support, not to the size of the dictionary.
\end{enumerate}
\end{proposition}

Each part answers a different objection.

Part (i) answers ``what if the structure is not in your dictionary''.  Nothing
is outside it: the atoms and the coarse space together span the symmetric
$L^2$, so a departure invisible to every atom is by definition part of the
background.  What the dictionary trades is not coverage but efficiency, since a
departure whose energy is spread thinly over many atoms will be detected poorly
even though it is representable.  That is the content of the power comparisons
in Section~\ref{sec:sim}, where a translated or additive departure is
representable and yet badly detected by a same-scale family.

Part (ii) separates signal from nuisance, and the identity is exact, not
asymptotic.  A network in which some coordinate regions are simply better
connected than others, with no interaction structure at all, produces zero
coefficient on every atom.  This is why additive main effects belong in the
background and why the hub structure that dominates many networks does not
generate false discoveries here.  The same calculation shows that a smooth
function of coordinate separation is orthogonal to second order in the atom
width, which is the weaker statement that motivates the combined nuisance
class.

Part (iii) is the formal content of the word \emph{localized}, and it is the
reason to scan a dictionary instead of estimating the surface.  Estimating
$f$ on $[0,1]^d\times[0,1]^d$ at resolution $h$ requires $h^{-2d}$ parameters
regardless of what is there.  A departure occupying a region of measure $A$
has $O(A/h^{2d})$ nonzero coefficients, so when $A$ is small the signal lives
in a vanishing fraction of the dictionary and a scan over the dictionary pays
only $\log$ of its size through the multiplicity correction.  Global mechanisms
such as preferential attachment, rewiring or a stochastic block model describe
a texture that is present everywhere and have no notion of a coefficient that
is zero except on a small region.  They cannot express, and therefore cannot
test, a statement of the form ``this structure, at this place, at this width''.

\subsection{Generative properties}                        \label{app:sec-thy-gen}

Throughout this subsection the coordinates are uniform on $[0,1]^d$, the
background $f_0$ in \eqref{eq:gen} is constant, and we write $\rho=\sigma(f_0)$,
$w=\rho(1-\rho)$ and $w'=\rho(1-\rho)(1-2\rho)$, so that
\begin{equation}
p_\theta=\rho+\theta w\Psi+\tfrac12\theta^{2}w'\Psi^{2}+O(\theta^{3}).
\label{app:eq-expand}
\end{equation}
In the sparse regime $\rho<\tfrac12$ and $w'>0$.  We write
\begin{equation}
\kappa_4=\int_0^1\psi^{4},\qquad\text{so that}\qquad\int\phi_a^{4}=\kappa_4^{\,d}
\label{app:eq-kappa4}
\end{equation}
by Fubini.  Fix a partition of $[0,1]^d$ into congruent dyadic blocks
$B_1,\dots,B_K$ at a resolution at least as fine as $h$.  Let $Q$ be the
modularity of that partition, $\mathrm{ES}$ the mean coordinate separation
across edges divided by its value across all pairs, $T$ the transitivity,
$\mathrm{CV}$ the coefficient of variation of the expected degree, and
$\lambda(\Psi)$ the spectrum of $\Psi$ as an integral operator.

Section~\ref{sec:catalogue} reported measured signatures.  This section derives
them.  Each network summary is an explicit functional of the dictionary
atom, computable in closed form in any dimension, and the functionals
separate the families independently of amplitude and of scale.

\subsubsection{The property map}

\begin{theorem}[Each summary is a functional of the atom]\label{app:thm-map}
Under \eqref{eq:gen}, as $\theta\to0$, with $\mu=\mathbb E\|u-v\|$,
\begin{align}
\mathbb E[Q]&=\frac{\theta w}{\rho}\sum_{k\le K}\int_{B_k\times B_k}\!\!\Psi\;+\;O(\theta^{2}),
\label{app:eq-mapQ}\\[2pt]
\mathbb E[\mathrm{ES}]&=1+\frac{\theta w}{\rho\,\mu}\bigl\langle\|u-v\|,\Psi\bigr\rangle+O(\theta^{2}),
\label{app:eq-mapES}\\[2pt]
\mathrm{CV}&=\frac{\theta^{2}w'}{2\rho}\Bigl(\textstyle\int\phi_a^{4}-1\Bigr)^{1/2}+O(\theta^{3}),
\label{app:eq-mapCV}\\[2pt]
\mathbb E[T]&=T_0+\frac{(\theta w)^{3}}{\rho^{3}}\operatorname{tr}(\Psi^{3})+O(\theta^{4}),
\label{app:eq-mapT}
\end{align}
and the spectral gap of the expected normalized Laplacian is monotone in
$\max\lambda(\Psi)$ and its multiplicity.  The first two expansions have no
zeroth-order correction because $\int\Psi=0$; the last two are the leading
nonzero terms.
\end{theorem}

\begin{proof}
Since $\int\Psi=0$, the edge density is $m=\rho+O(\theta^{2})$ and the block
degree is $d_k=\rho|B_k|+O(\theta^{2})$, the latter because
$\int_{B_k}\!\int\Psi=\int_{B_k}\phi_a\cdot\int\phi_a=0$.  Substituting
\eqref{eq:expand} into $Q=\sum_k\{e_{kk}/m-(d_k/2m)^{2}\}$ and into
$\int\|u-v\|p_\theta/\int p_\theta$ gives \eqref{eq:mapQ} and \eqref{eq:mapES}.

For \eqref{eq:mapCV}, the expected degree is
$D(u)=\rho+\theta w\!\int\!\Psi(u,v)dv+\tfrac12\theta^{2}w'g(u)+O(\theta^{3})$
with $g(u)=\int\Psi^{2}(u,v)dv$.  The first-order term vanishes, so
$\mathrm{sd}(D)=\tfrac12\theta^{2}w'\,\mathrm{sd}(g)$.  For a diagonal atom
$g=\phi_a^{2}\int\phi_a^{2}=\phi_a^{2}$, whence
$\mathrm{Var}(g)=\int\phi_a^{4}-1$.

For \eqref{eq:mapT}, the first-order triangle term is proportional to
$\int\!\int\!\int\Psi(u,v)=0$ and the second-order cross term to
$\int dv\{\int\Psi(u,v)du\}\{\int\Psi(v,w)dw\}=0$, so the leading correction is
$(\theta w)^{3}\int\!\int\!\int\Psi(u,v)\Psi(v,w)\Psi(u,w)=(\theta w)^{3}\operatorname{tr}(\Psi^{3})$.
The path density is unchanged to this order.  The spectral statement is the
perturbation of the normalized Laplacian by the finite-rank operator
$\theta w\Psi/\rho$.
\end{proof}

Theorem~\ref{thm:map} reduces the catalogue to four integrals.  The corollaries
evaluate them.

\subsubsection{Diagonal and off-diagonal atoms}

\begin{corollary}[Communities versus bridges]                \label{cor:diag}
Let $a\ne b$ be boxes at scale $h$ and $\theta>0$.
\begin{enumerate}[label=\textup{(\roman*)},leftmargin=*]
\item \emph{Diagonal, $\Psi=\phi_a\otimes\phi_a$.}
$\sum_k\int_{B_k\times B_k}\Psi=(h/2)^{d}>0$,
$\langle\|u-v\|,\Psi\rangle=-c_dh^{d+1}<0$ with $c_1=1/6$,
$\operatorname{tr}(\Psi^{3})=1$, and $\lambda(\Psi)=\{1\}$ of multiplicity one.
Modularity rises, edges shorten, transitivity rises, the perturbation is rank
one.
\item \emph{Off-diagonal, $\Psi=(\phi_a\otimes\phi_b+\phi_b\otimes\phi_a)/\sqrt2$.}
$\sum_k\int_{B_k\times B_k}\Psi=0$, $\langle\|u-v\|,\Psi\rangle=0$,
$\operatorname{tr}(\Psi^{3})=0$, $\lambda(\Psi)=\{\pm1/\sqrt2\}$.  Both
first-order effects vanish.  At second order
$\mathbb E[Q]=-\theta^{2}w'/2\rho<0$ and
$\mathbb E[\mathrm{ES}]=1+(\theta^{2}w'/2\rho)(\Delta-\mu)/\mu>1$ whenever the
boxes are separated by more than an average pair, $\Delta>\mu$.
\end{enumerate}
\end{corollary}

\begin{proof}
For (i), the partition resolves the halves of $a$, so
$\int_{B_k\times B_k}\Psi=(\int_{B_k}\phi_a)^{2}\ge0$ and the sum over the
blocks tiling $a$ is $(h/2)^{d}$; squares cannot cancel, which forces the sign.
The separation integral is the $d=1$ computation
$h^{-1}(2h^{3}/24-2h^{3}/8)=-h^{2}/6$ tensored over the remaining axes.  As an
operator $\Psi$ is the rank-one projection onto $\phi_a$, so $\Psi^{3}=\Psi$
and $\operatorname{tr}(\Psi^{3})=1$.

For (ii), $B_k\times B_k$ meets $(a\times b)\cup(b\times a)$ only if $B_k$
meets both disjoint boxes, which is impossible, so the modularity integral
vanishes.  The separation integral vanishes because $\|u-v\|$ is affine on
$a\times b$ and both factors integrate to zero.  As an operator $\Psi$ is
$2^{-1/2}$ times the swap of $\phi_a$ and $\phi_b$, whose cube is again a
multiple of the swap and has zero diagonal.  Since $\phi_a\phi_b\equiv0$,
$\Psi^{2}=\tfrac12(\phi_a^{2}\otimes\phi_b^{2}+\phi_b^{2}\otimes\phi_a^{2})$ is
supported strictly between blocks at separation $\Delta+O(h)$, so the
second-order term raises the edge total without raising any $e_{kk}$, lowering
$Q$, and moves conditional separation toward $\Delta$.
\end{proof}

\subsubsection{Wavelet shape}

\begin{corollary}[Shape and dimension]                       \label{cor:shape}
For a diagonal atom built from any mother wavelet $\psi$ with $\int\psi=0$
and $\int\psi^{2}=1$, the quantities $\sum_k\int_{B_k\times B_k}\Psi$,
$\langle\|u-v\|,\Psi\rangle$, $\operatorname{tr}(\Psi^{3})$ and
$\lambda(\Psi)$ do not depend on the shape of $\psi$.  The only
shape-dependent summary is degree heterogeneity, through the fourth moment
\eqref{eq:kappa4}:
\[
\mathrm{CV}=\frac{\theta^{2}w'}{2\rho}\bigl(\kappa_4^{\,d}-1\bigr)^{1/2}+O(\theta^{3}).
\]
\end{corollary}

\begin{proof}
Immediate from Corollary~\ref{cor:diag}(i), whose evaluations use only
$\int\phi_a=0$ and $\int\phi_a^{2}=1$, together with \eqref{eq:mapCV} and
$\int\phi_a^{4}=\kappa_4^{\,d}$.
\end{proof}

Three consequences.  Position on the diagonal, not shape, decides whether an
element builds a community or a bridge.  Shape controls only how unevenly
degree is distributed, through a single scalar: $\kappa_4=4$ for the Haar
wavelet, $6.78$ for the Mexican hat and $8.75$ for the Daubechies wavelet with
two vanishing moments, so a more concentrated wavelet gives a denser core and a
thinner tail at the same scale and amplitude.  Dimension enters
multiplicatively as $\kappa_4^{\,d}$, so degree heterogeneity grows
exponentially in $d$: the Haar values are $1.73$, $3.87$, $7.94$ for
$d=1,2,3$, and the Daubechies values $2.78$, $8.69$, $25.86$.

\begin{remark}[Scaling functions are not wavelets]
The zero-mean property is what makes degree heterogeneity second order.  A
scaling function has $\int\psi\ne0$, so $\int\Psi(u,v)dv=\phi(u)\int\phi\ne0$
and the first-order term in $D(u)$ survives: such an atom behaves like an
additive main effect and generates degree heterogeneity at order $\theta$
instead of $\theta^{2}$.  The distinction is easy to introduce by accident:
the four Daubechies coefficients $(0.4830,0.8365,0.2241,-0.1294)$ define the
scaling function and sum to $\sqrt2$, while the wavelet is the alternating flip
$g_k=(-1)^{k}h_{3-k}$, which sums to zero.
\end{remark}

\subsubsection{Mixtures across scales}

\begin{corollary}[One scale versus several]                   \label{cor:mix}
Let $\Psi=\sum_{k\le M}c_k\,\phi_{a_k}\otimes\phi_{a_k}$ with
$\sum_kc_k^{2}=1$ and the $\phi_{a_k}$ orthonormal.  Then
$\operatorname{tr}(\Psi^{3})=\sum_kc_k^{3}$ and
$\lambda(\Psi)=\{c_1,\dots,c_M\}$.  In particular, $M$ atoms at one scale
with equal weights give $\operatorname{tr}(\Psi^{3})=M^{-1/2}$ and $M$ equal
eigenvalues: a block model on intervals, with modularity the sum of the
single-element values, low transitivity and a large spectral gap.  Atoms at
$M$ different scales with equal weights give the same trace when their supports
are disjoint and strictly more when nested, because the operator products no
longer vanish across scales: a hierarchy has higher transitivity than a block
model of the same rank, and a smaller spectral gap.
\end{corollary}

For the configurations of Section~\ref{sec:catalogue} these give
$\operatorname{tr}(\Psi^{3})=8^{-1/2}=0.354$ for eight blocks at one scale and
$0.500$ for four nested scales, against measured transitivity $0.147$ and
$0.495$ and measured spectral gaps $0.296$ and $0.082$.

\subsubsection{Background terms}

\begin{proposition}[Background terms]                         \label{prop:background}
\textup{(i)} For $\Psi_{\mathrm{add}}(u,v)=m(u)+m(v)$ the modularity of any
partition is zero to first order, the within-block and degree terms of $Q$
cancelling exactly.  The surface generates degree heterogeneity at order
$\theta$ and no community structure.  \textup{(ii)} For
$\Psi_s(u,v)=s(\|u-v\|)$ with $s$ continuously differentiable,
$\langle s,\Psi_{a,b}\rangle=0$ for off-diagonal atoms and
$|\langle s,\Psi_{a,a}\rangle|\le C_dh^{d+1}\|s'\|_\infty$ for diagonal ones: a
smooth gradient is invisible to the dictionary up to a term of order
$h^{d+1}$ and generates no boundaries.
\end{proposition}

Both statements are the reason these terms belong in the background class of
Section~\ref{sec:nuisance} and not in the dictionary.

\subsubsection{Scope of the derivation}

\begin{table}[t]\centering
\caption{Closed forms for the summaries each atom generates, from
Theorem~\ref{thm:map} evaluated by Corollaries~\ref{cor:diag} to~\ref{cor:mix}.
Entries are the leading nonzero term in the amplitude $\theta$, valid at every
scale $h$ and in every dimension $d$.  Table~\ref{tab:map} reports the
corresponding measurements.}                                        \label{app:tab-map}
\begin{tabular}{lcccc}
\hline
element & modularity & edge span & transitivity & degree CV\\
\hline
diagonal              & $+(h/2)^{d}$ & $-c_dh^{d+1}$ & $+1$ & $(\kappa_4^{d}-1)^{1/2}$\\
off-diagonal          & $<0$, second order & $>1$, second order & $0$ & $(\kappa_4^{d}-1)^{1/2}$\\
$M$ blocks, one scale & $+M(h/2)^{d}$ & $-Mc_dh^{d+1}$ & $M^{-1/2}$ & as above\\
$M$ nested scales     & $+$ & $-$ & $>M^{-1/2}$ & as above\\
additive              & $0$ & $0$ & $0$ & order $\theta$\\
separation kernel     & $O(h^{d+1})$ & $O(h^{d+1})$ & $0$ & $0$\\
\hline
\end{tabular}
\end{table}

The results are proved for expected summaries of the generative model as
$\theta\to0$, with $f_0$ constant and coordinates uniform.  Three limitations.
The expansions are asymptotic in amplitude, and the simulations of
Section~\ref{sec:catalogue} run at a finite amplitude at which the ratio of
measured to predicted degree spread is $1.00$ for the Haar atom and $1.01$
to $1.08$ for the more concentrated ones.  The statements concern expected
summaries and not their fluctuations across realizations.  And the relation
between $\lambda(\Psi)$ and the spectral gap is stated as monotone and is not
evaluated, since the normalized Laplacian also involves the degree profile.
The eigenvalues themselves are exact.

\subsection{Surfaces without localized content}           \label{app:sec-thy-invis}

Write $\mathcal G$ for the span of the three canonical coordinate-based
mechanisms: the constant surface of an Erd\H os--R\'enyi model, the additive
surfaces $m(u)+m(v)$ of degree-corrected and configuration models, and the
separation surfaces $s(\|u-v\|)$ of latent-distance and small-world models.
Write $\mathcal G_{J_0}$ for $\mathcal G$ enlarged by the block surfaces of a
stochastic block model whose blocks are unions of dyadic boxes at resolution
$J_0$.

\begin{theorem}[No localized statement is identifiable from a global class]
\label{thm:noloc}
Let $\Psi_{a,b}$ be a atom at scale $h=2^{-j}$ with $j>J_0$, and
let $\delta\in\mathcal G_{J_0}$ with $s$ continuously differentiable.  Then
\[
\langle\delta,\Psi_{a,b}\rangle=0\quad\text{for }a\ne b,
\qquad
\bigl|\langle\delta,\Psi_{a,a}\rangle\bigr|\le C_d\,h^{\,d+1}\|s'\|_\infty .
\]
Consequently every member of $\mathcal G_{J_0}$ has the same dictionary
coefficients, namely zero up to $O(h^{d+1})$, and no statement of the form
``a departure of width $h$ is present at position $a$'' is identifiable within
that class.  By contrast $\theta_\lambda=\langle\delta,\Psi_\lambda\rangle$ is
exactly the coefficient the scan estimates, so within the dictionary the
position and the width are identified.
\end{theorem}

\begin{proof}
Constants and block surfaces at resolution $J_0$ lie in $V_{J_0}$, which is
orthogonal to every atom at a finer scale by
Proposition~\ref{prop:dict}\textup{(i)}.  Additive surfaces are orthogonal to
every atom by Proposition~\ref{prop:dict}\textup{(ii)}.  For the separation
surface, Proposition~\ref{prop:background}\textup{(ii)} gives exact
orthogonality for off-diagonal atoms, since $\|u-v\|$ is affine on
$a\times b$ and both Haar factors integrate to zero, and the stated bound for
diagonal atoms by a first-order Taylor expansion of $s$ about the centre of
the box, the constant term being annihilated and the linear term contributing
$\langle\|u-v\|,\Psi_{a,a}\rangle=-c_dh^{d+1}$ from
Corollary~\ref{cor:diag}\textup{(i)}.  Summing the three contributions gives
the claim.
\end{proof}

\begin{remark}
The theorem is the formal content of a familiar informal contrast.  A global
mechanism specifies how a network is textured everywhere: preferential
attachment makes hubs, rewiring makes short paths, a distance kernel makes
proximity.  Each of these is a single statement about the whole surface, and
Theorem~\ref{thm:noloc} says that each is invisible to the dictionary, which
is the same as saying that the dictionary's question is not one they answer.
The dictionary specifies instead what happens on one box of the coordinate
plane and leaves the rest of the surface to the background.  Its coefficient is
a statement about a place and a width, and
Sections~\ref{sec:method} and~\ref{sec:results} attach a standard error and a
simultaneous critical value to it.  The three mechanisms above are not
competitors to be beaten on power.  They are the background against which a
localized claim is defined.
\end{remark}

\subsection{Inference guarantees}                         \label{app:sec-thy-inf}

\subsubsection{Assumptions}

This section states four theorems and one proposition, with proofs in the
Supplement: Theorem~\ref{thm:minimax} in Appendix~\ref{sec:supp-detect},
Theorem~\ref{thm:order} and Proposition~\ref{prop:obstruct} in
Appendix~\ref{sec:supp-order} together with Theorem~\ref{thm:floor}, and
Theorem~\ref{thm:uncond} in Appendix~\ref{sec:supp-uncond}.  The calibration
engine on which Theorem~\ref{thm:uncond} rests, a conditional
Bernstein--Taylor reduction and a shared-dyad multiplier calibration, is
stated and proved in Appendix~\ref{sec:supp-finite}, Appendix~\ref{sec:supp-mult}
and Appendix~\ref{sec:supp-mult}.  Theorems~\ref{thm:minimax}
and~\ref{thm:order} assume displayed linearization, studentizer, covariance,
packing and score-transfer rates.  Those rates are derived from primitive
conditions only in the instance of Theorem~\ref{thm:uncond}, and population
inference requires an explicit finite-design-to-Lebesgue bias condition,
without which the inferential centre is the frozen finite-design target.

Write $\ell_n=\log(2m_n\vee3)$.  The supplement's per-section variants
$\log(2m_n)$, $\log(e\bar m_n)$, and $\log\bar d_n$ agree with $\ell_n$ up to
bounded ratios on the admission event, every rate below is insensitive to the
choice, and each supplement section states its local variant at first use.
Similarly $m_n$ is the deterministic packing cardinality here and the realized
admitted cardinality where a supplement section says so.  Several symbols are
overloaded across sections and are always disambiguated by context or subscript:
$\delta$ is the departure $(I-\Pi_{V_{J_0}})f$, $\delta_n$ the fine-cell width,
$\delta_{V,n},\delta_{G,n}$ studentizer errors, and grid ``$\delta$'' in
Section~\ref{sec:sim} the simulation cell width.  $\Gamma_{\lambda,n}$ is the
curvature envelope, $\Gamma_n$ the transfer margin of
Theorem~\ref{thm:order}, and $\Gamma_{ab,n}$ the multiplier covariance.  $H_0$
is the baseline Jacobian and $H=nh_*$ the block length of
Proposition~\ref{prop:obstruct}.  $R^F$ are fold masks, $R$ in Section~\ref{sec:sim}
a replication count and in Section~\ref{sec:thy-limits} a coordinate range.  $k$ indexes
dyads in Theorem~\ref{thm:linear} and cells elsewhere.  $\Delta_n$ is the realized
coordinate error in Section~\ref{sec:model} and the realized rank error in
Theorem~\ref{thm:order}, and $\kappa$ is the standardized amplitude of the
simulations, distinct from $\kappa_0$ and $\kappa_n^{\mathrm{str}}$.  The finite-design and calibration results permit an
$\mathcal F_{0n}$-measurable admitted family $\cI_n$.  The global minimax result instead
uses a deterministic master packing and an all-pass admission event, as stated below.
All atoms are boundary-aligned, isotropic, symmetric Haar atoms.  Their discontinuity
is handled by the boundary-layer inequality
\begin{equation}
 \|\Psi_\lambda^S(\,\cdot-t)-\Psi_\lambda^S\|_1
 \le C\|t\|_\infty,
 \qquad \|t\|_\infty\le c h_\lambda,                            \label{app:eq-bvtranslate}
\end{equation}
not by differentiating the atom.  Off-diagonal atoms use the normalization in
\eqref{eq:symatom}, with transposed duplicates removed.

\subsubsection{Calibration engine}                                  \label{sec:engine}

The calibration engine consists of a conditional Bernstein--Taylor reduction
of each studentized coefficient to a sum of independent dyad contributions
(Theorem~\ref{thm:linear}), deterministic curvature and curvature envelopes and a
studentizer comparison (Supplement Propositions~S2 and~S3), and a uniform
shared-dyad multiplier calibration that gives simultaneous coverage and
familywise error control conditional on the folds (Theorem~\ref{thm:fwer}).
Under $\Pp(\Delta_n>r_n)\le\eta_n$ the unconditional bound is
$\alpha+\eta_n+o(1)$, and the budget is additive: the applications of
Sections~\ref{sec:thy-limits} and~\ref{sec:thy-limits} use a $95\%$ radius and a $5\%$
test and so carry a nominal familywise bound of $0.10$.  These results are
conditional on displayed rates, and Theorem~\ref{thm:uncond} discharges the
rates in one instance.  Their statements and proofs are in the Supplement.

\subsubsection{Detection and localization boundary}

Let $\Lambda_n=\bigcup_{h\in\cJ_n}\Lambda_n(h)$ be a deterministic,
outcome-independent packing of unique off-diagonal symmetric atoms, and set
$m_n=|\Lambda_n|\to\infty$.  A pointwise-amplitude alternative is
\begin{equation}
 \cA(a,h)=
 \left\{f=f_0+sah\Psi_\lambda^S:
 f_0\in\mathcal B_n,\ s\in\{-1,1\},\ \lambda\in\Lambda_n(h),\
 \|f_0+sah\Psi_\lambda^S\|_\infty\le B\right\}.               \label{eq:altclass}
\end{equation}
Since $\|\Psi_\lambda^S\|_\infty\asymp h^{-1}$ in two dimensions, $a$ is the
order of the log-odds change on the support.  If $S_\lambda$ denotes the symmetric
support, define
\begin{equation}
 d_{\rm loc}(\widehat\lambda,\lambda)
 =\min\left[1,
 \frac{|S_{\widehat\lambda}\triangle S_\lambda|}{2|S_\lambda|}
 +\1\{h_{\widehat\lambda}\ne h_\lambda\}\right],               \label{eq:locloss}
\end{equation}
with loss one when no atom is selected.

\begin{theorem}[Global adaptive detection and localization]\label{thm:minimax}
Suppose the atoms are orthonormal in $V_{J_0}^{\perp}$, atoms at a common scale
have equal support area, and distinct same-scale supports have normalized
symmetric-difference separation bounded below.  Let $\widehat\Lambda_n$ be the
pilot-admitted subfamily, require the all-pass event
$\{\widehat\Lambda_n=\Lambda_n\}$, and make no rejection when it fails.  Assume its
failure probability is at most $\gamma_n=o(1)$ uniformly over the null and
single-atom classes.  On a further $\mathcal F_{0n}$-measurable good event whose
failure probability is at most $\eta_n+o(1)$, assume uniformly that
\begin{itemize}[leftmargin=1.6em,itemsep=2pt,topsep=2pt]
\item $I_{\lambda,f}\asymp n^2\rho_n$ and the standardized influence summands are
conditionally independent within each coordinate, variance normalized, and bounded
by $L_n$, with $L_n^2\log^7(m_nn)\to0$.
\item the max-norm linearization error is $o_{\Pp\mid\mathcal F_{0n}}(\ell_n^{-1/2})$,
the relative studentizer error is $o_{\Pp\mid\mathcal F_{0n}}(\ell_n^{-1})$, and
$\sqrt{\ell_n}\max_\lambda
|\theta_{\lambda,n}^{D}-\theta_\lambda(f)|/s_{\lambda,f}=o_\Pp(1)$, and
\item the multiplier calibrates the actual all-pass statistic under the global null.
\end{itemize}
Then the all-pass scan has global-null size at most $\alpha+\eta_n+o(1)$.
There is $C_{\rm up}<\infty$ such that, for every nonempty sequence of classes
$\cA(a_n,h_n)$ satisfying
\begin{equation}
 a_n^2n^2\rho_n h_n^2
 \ge C_{\rm up}\log(2m_n),                                     \label{eq:upperboundary}
\end{equation}
the conditional probability of failure to reject, recover the sign, or recover the
exact atom tends to zero uniformly on the good all-pass event.  Consequently,
\begin{equation}
 \sup_{f\in\cA(a_n,h_n)}\E_f d_{\rm loc}(\widehat\lambda,\lambda(f))
 \le\eta_n+\gamma_n+o(1)\longrightarrow0.                        \label{eq:localup}
\end{equation}

For the lower bound, condition on a deterministic regular design.  Suppose the
baseline class contains a bounded interior element with Bernoulli probabilities of
order $\rho_n$.  The evaluated atoms have bounded pointwise perturbation and empirical
energy of order $n^2h_\lambda^2$ for the scaled atom $h_\lambda\Psi^S_\lambda$
(the normalized $\Psi^S_\lambda$ has empirical energy of order $n^2$), and the
ordered support-overlap count satisfies
$N_{{\rm ov},n}\le m_nd_n$ with
$\log d_n/\log(2m_n)\to0$.  Assume also that the original experiment is a
parameter-independent randomization of the full-outcome oracle experiment.  There is
$c_{\rm low}>0$ such that, if the amplitude envelope obeys
\begin{equation}
 \max_{h\in\cJ_n}a_h^2n^2\rho_n h^2
 \le c_{\rm low}\log(2m_n),                                    \label{eq:lowerboundary}
\end{equation}
then the global-union minimax sum of type-I and worst-case type-II errors tends to
one, and the minimax canonical localization risk has a strictly positive limit
inferior.
\end{theorem}

Thus the effective sample size is $n^2\rho_nh^2$: the number of dyads in
the support multiplied by the sparse probability scale.  The logarithm reflects the
unknown location and scale.  The lower result is global over the union of scales.  A
matching lower statement at each individual scale requires additional scale-wise
packing and overlap conditions.  No leading-constant claim is made.  The theorem
applies to the deterministic all-pass packing, not automatically to an arbitrary
random admitted subset.

\subsubsection{Unconditional instance}

The four theorems above are conditional on displayed rates.  The instance
below replaces them by primitive conditions, of which the cell-information
condition (C1) is the headline and (C2) the substantive companion: for bounded
expected degree and a fixed coarsest scale, (C2) requires
$\delta_n\gg n^{-1/4}\log^{3/4}n$, while (C1) alone permits cells as fine as
$n^{-1/3}$.  The smallest resolution the theorem covers is set by (C2), not by
cell information.  The final technical
section of the Supplement derives every one of those rates from primitive
quantities, network size, sparsity, cell width, wavelet scale, and
logarithms, for one complete instance of the estimator, producing the paper's
first unconditional guarantee.  The instance fixes the nuisances the general theory leaves open (known
regular grid $x_i=(i-\tfrac12)/n$, radius zero, complete observation, a fixed
nonnegative pre-smoother, full-dictionary admission) and keeps the full
inferential pipeline (disjoint dyad folds, background estimated on its own fold
with stabilization and clipping, exact coarse logit-scale projection,
shared-dyad multiplier maximum, multiscale scan).  With fine-cell width $\delta_n=2^{-J_n}$, scales
$h\in[2^{-j_{\max,n}},2^{-j_0}]$, $\delta_n\le h_{\min,n}/2$, and a fixed
$\gamma>0$, the primitive conditions are
\begin{equation}
\text{(C1)}\ \ n^2\rho_n\delta_n^2\ge\log^{7+\gamma}(n),
\qquad
\text{(C2)}\ \ \frac{h_{\max}\log^{3/2}(n)}{\sqrt{\rho_n}\,n\,\delta_n^2}\to0,
\qquad
\text{(C3)}\ \ \kappa_0\sqrt{\rho_n\log n}=o(n\delta_n^2),
\label{eq:uncond-conds}
\end{equation}
with, for the detection--localization clause,
\text{(C4)}: $j_{\max,n}\to\infty$ and $j_{\max,n}/\log\log n\to\infty$.

\begin{theorem}[Unconditional guarantees for the regular-design instance]
\label{app:thm-uncond}
For the instance just described, formalized as \textup{(S1)--(S8)} in the
Supplement, and under \textup{(C1)--(C3)} of \eqref{eq:uncond-conds}, the
following hold uniformly over the null and single-atom classes, with no
conditioning event and no unverified assumption.
\begin{enumerate}[label=\textup{(\alph*)},leftmargin=2em,itemsep=2pt]
\item The uniform linearization around the population coefficients and the
studentizer consistency of Theorem~\ref{thm:linear} and
Propositions~S2 and~S3 hold with explicit envelope
rates.  In particular the maximal standardized target bias is
$O\{\kappa_0h_{\max}\sqrt{\rho_n\log n}\,(n\delta_n^2)^{-1}\}$.
\item The shared-dyad multiplier of Theorem~\ref{thm:fwer} is valid: the
simultaneous intervals cover the full coefficient vector with probability at
least $1-\alpha-o(1)$, and the studentized tests control the familywise error at
$\alpha+o(1)$.  The verification in Appendix~\ref{sec:supp-uncond} covers the global null
and every single-atom configuration of the class, and arbitrary partial-null or
multiple-atom configurations are not verified for this instance.
\item Under \textup{(C4)} there are constants $0<c_{\rm low}<C_{\rm up}<\infty$
such that single-atom alternatives with
$a_n^2n^2\rho_nh_n^2\ge C_{\rm up}\log(2m_n)$ are detected, signed, and exactly
localized with probability tending to one, while on a matched-pair subpacking
with $\max_ha_h^2n^2\rho_nh^2\le c_{\rm low}\log(2m_n^{\rm low})$ and
$\log(2m_n^{\rm low})\asymp\log(2m_n)$, the global-union detection risk tends to
one and the localization risk stays bounded away from zero.  The boundary
$a^2n^2\rho_nh^2\asymp\log(2m_n)$ is therefore attained, unconditionally, by the
implemented statistic in this instance.
\end{enumerate}
\end{theorem}

One condition drives everything: once each fine cell carries
$\log^{7+\gamma}n$ effective observations, calibration, coverage, familywise
error control, and the matching rate boundary with exact single-atom selection all
follow, and (C2)--(C3) are
mild side conditions on the scale range and the fixed pseudocount.  Two
structural facts make the derivation close.  Every surface in the class is
exactly constant on fine cells, so the finite-design target coincides with the
population coefficient up to a stabilization term of order
$\kappa_0/(n^2\delta_n^2)$, and the identity $\Pi_{J_0}^\top c_\lambda=0$, the
atom's orthogonality to the coarse space preserved exactly by the final
projection step of the baseline map, removes the estimated background from the
coefficient functional along the subtracted path at every order.  This is the
structural reason the background can be estimated on a separate fold without a
first-order price, and it persists verbatim for the general estimator.  The
instance excludes random coordinates, estimated propensities, missingness,
data-selected bandwidths, sign-changing smoother weights, and positive order
radius.  That list, recorded in the Supplement, is now the exact remaining gap
between the proved theory and the implemented general procedure.

\subsection{Limits on dimension and on coordinate accuracy} \label{app:sec-thy-limits}

Two results close the account.  The first prices the construction in dimension
$d$ and shows that the binding constraint is information, not
multiplicity.  The second says what the dictionary buys: the canonical
coordinate-based mechanisms cannot express a localized departure at all, so the
question of where and how wide is not merely unanswered by them but
unformulable.

\subsubsection{Coordinate dimension}

\begin{theorem}[Information requirement]                      \label{app:thm-info}
Scan a dictionary at width $h$ in dimension $d$ using cells of width $h/2$.
For the studentized coefficient of every admitted atom to be based on cell
counts whose expectation is at least a constant $C$, it is necessary that
\begin{equation}
n\,\bar d\;\ge\;C\,4^{\,d}h^{-2d},
\label{app:eq-inforule}
\end{equation}
with $\bar d=n\rho_n$ the mean degree.  The dictionary itself has
$m_n\asymp h^{-2d}/2$ atoms, so the multiplicity correction entering the
critical value is
\begin{equation}
\sqrt{2\log 2m_n}\;\asymp\;2\sqrt{d\,\log(1/h)} .
\label{app:eq-multrule}
\end{equation}
\end{theorem}

\begin{proof}
A cell is a box of side $h/2$, so it contains $n(h/2)^{d}$ vertices in
expectation and a pair of cells contains $n^{2}(h/2)^{2d}$ dyads, each an edge
with probability of order $\rho_n$.  The expected count in a cell pair is
therefore $n^{2}\rho_n(h/2)^{2d}$, and requiring this to be at least $C$ gives
$n^{2}\rho_n\ge C4^{d}h^{-2d}$, which is \eqref{eq:inforule} since
$n^{2}\rho_n=n\bar d$.  For the second claim, at scale $h=2^{-j}$ there are
$2^{jd}=h^{-d}$ boxes and hence $h^{-d}(h^{-d}+1)/2$ unordered pairs of boxes,
so $\log 2m_n=2d\log(1/h)+O(1)$ and \eqref{eq:multrule} follows.
\end{proof}

The two rates are not comparable.  Information is exponential in $d$ and
polynomial of degree $2d$ in $1/h$; multiplicity is the square root of a linear
function of $d$.

\begin{corollary}[What is affordable]                         \label{cor:afford}
At $h=1/8$ and $C=30$, the requirement \eqref{eq:inforule} is
$7.7\times10^{3}$ for $d=1$, $2.0\times10^{6}$ for $d=2$ and
$5.0\times10^{8}$ for $d=3$, while the multiplicity factor
\eqref{eq:multrule} is $2.92$, $4.08$ and $5.00$.  Passing from one dimension
to three multiplies the information requirement by $6.6\times10^{4}$ and the
multiplicity penalty by $1.71$.
\end{corollary}

The practical consequence is a hard ordering.  A one-dimensional coordinate is
affordable on ordinary networks: a protein chain of $10^{3}$ residues at mean
degree $27$ supplies $2.9\times10^{4}$ and resolves $h=0.065$.  A
two-dimensional coordinate requires a large and dense network and then
resolves only coarse structure.  A three-dimensional coordinate is out of reach
at any size we have encountered: a cortical parcellation of $400$ regions at
mean degree $20$ supplies $8.0\times10^{3}$ against the $7.9\times10^{6}$ that
$d=3$ demands at $h=1/4$, short by three orders of magnitude.  Degree
heterogeneity compounds the difficulty, since by Corollary~\ref{cor:shape} it
grows as $(\kappa_4^{\,d}-1)^{1/2}$ and so thins the usable cells further as
$d$ rises.

Every application in Section~\ref{sec:apps} uses a coordinate that is
physically measured and exact.  That is not a convenience.  It is the
condition under which the questions this paper asks have answers, and this
section states the condition, proves the part of it that is a theorem, and
gives a check that can be run on a candidate coordinate before any analysis.

\subsubsection{Ordering transfer and obstruction}

Here both $R$ and $\widehat R$ are permutations, ties having been resolved by a fixed
outcome-free rule, and
$\Delta_n=\max_i|\widehat R_i-R_i|/n$.  Put
$\overline r_n=r_n+n^{-1}$.  The theorem is specialized to normalized ranks.
conversion from a physical coordinate requires separate density bounds.

\begin{theorem}[Normalized-rank ordering transfer and obstruction]         \label{thm:order}
Let the oracle and reported-order scores use the same deterministic index set and be
constructed from the same scoring outcomes, folds, and multiplier seeds.  Suppose
each Haar atom has $L^2$ norm one, sup norm $O(h_\lambda^{-1})$,
support area of order $h_\lambda^2$, and a bounded number of axis-aligned jump
segments of total length $O(h_\lambda)$.  If
$\overline r_n/h_{\min,n}\to0$, then on $\{\Delta_n\le r_n\}$ the bin permutation
induced by $(R,\widehat R)$ is measure preserving and, for every bounded full oracle
contrast $g_{P,n}$ (including its coarse component),
\begin{equation}
 \max_{\lambda\in\Lambda_n}
 |\theta_{\lambda,n}^{\rm ord}-\theta_{\lambda,n}^{\rm or}|
 \le CB\overline r_n.                                           \label{eq:ordercoef}
\end{equation}
For a single-atom signal $q_{P,n}=sa_nh_*\Psi_{\lambda_*}^S$, the oracle coefficient
is $sa_nh_*$, and $B\overline r_n/(|a_n|h_*)\to0$ is sufficient to preserve its
sign and first-order magnitude.

Write the exact coupled score identity as
\[
 Z_{\lambda,n}^{\rm ord}=Z_{\lambda,n}^{\rm or}
 +\sqrt{I_{\lambda,n}}\,b_{\lambda,n}^{\rm sc}+R_{\lambda,n},
\]
and let $D_n$ be the sum of the maximal standardized coefficient perturbation, the
explicit score-linearity discrepancy, and the maximal centered-score remainder.
If the oracle statistic separates the true signed coordinate, its critical value,
and every competing atom by at least $3\Gamma_n$, while
\[
 D_n+|\widehat c_n^{\rm ord}-\widehat c_n^{\rm or}|\le\Gamma_n
\]
with probability tending to one on $\{\Delta_n\le r_n\}$, then the reported-order
scan detects, reports the correct sign, and selects the exact canonical label with
probability tending to one.  Hence its expected canonical loss is at most
$\eta_n+o(1)$.  If $I_{\lambda_*,n}\asymp n^2\rho_n$, the oracle critical value is
of order $\sqrt{\ell_n}$, and the last perturbation display is
$o_\Pp(\sqrt{\ell_n})$, this transfers the oracle's first-order \emph{sufficient}
upper boundary \eqref{eq:upperboundary} on alternatives separated from it by a fixed
multiplicative factor.  It does not establish an exact boundary power curve.

Under the null, if the oracle maximum has a uniform Gaussian approximation with
coordinate variances bounded above and away from zero, then
$\sqrt{\ell_n}D_n=o_\Pp(1)$ implies uniform equivalence of the oracle and
reported-order maximum distribution functions.

These canonical conclusions must not be replaced by a boundary-strip remainder.
For a fixed label, the membership sensitivity satisfies
\[
 \frac{|P_{T_n}^{-1}S_\lambda\triangle S_\lambda|}{2|S_\lambda|}
 \le C\min\left\{1,
 \frac{\Delta_n+n^{-1}}{h_\lambda}
 +\left(\frac{\Delta_n+n^{-1}}{h_\lambda}\right)^2\right\},
\]
whereas exact selection of $\lambda_*$ has canonical loss zero.
\end{theorem}

The transfer clauses of Theorem~\ref{thm:order} operate when the reported radius is
small relative to the scanned width.  The complementary case is
when it is not.  The next proposition answers it in a separate rank-space
experiment in which the analyst observes the network, the reported ranks, and the
radius, but not the reporting mechanism: once the radius reaches the atom width,
no procedure can recover the canonical location, even when the signal amplitude is
far above the detection boundary.  Formally, let
$\mathcal Q_n^{\mathrm{perm}}(r)$ be the set of Markov kernels that, given the
true rank permutation $R$, return a reported permutation $\widehat R$ with
$\|\widehat R-R\|_\infty\le nr$ almost surely.  The parameter is
$\vartheta=(R,Q,\lambda,s,a)$ with $Q\in\mathcal Q_n^{\mathrm{perm}}(r_n)$, and
the observable is $(A,\widehat R,r_n)$ with edge probabilities
$\operatorname{logit}^{-1}\{\alpha_n+sah\,\Psi_\lambda^S(x_{R_i},x_{R_j})\}$ at
the true-rank midpoints $x_q=(q-\tfrac12)/n$.

\begin{proposition}[Dyadic canonical-location obstruction]
\label{app:prop-obstruct}
Let $H=nh_*$ be an even integer, and suppose the aligned width-$h_*$ dictionary
contains three pairwise-disjoint blocks $I_0,I_1,J$ of $H$ consecutive ranks each,
with $I_1$ adjacent to $I_0$.  Write $\lambda_0=(I_0,J)$ and $\lambda_1=(I_1,J)$.
Restrict the amplitudes to $a$ with $C_\Psi^\circ a\le M\le B$ and
$a^2n^2\rho_nh_*^2\ge\kappa_n^{\mathrm{str}}\ell_n$ for some
$\kappa_n^{\mathrm{str}}\to\infty$, so the signals are bounded on the logit scale
and lie strictly above the detection boundary of Theorem~\ref{thm:minimax}.
assume this amplitude class is nonempty, equivalently
$(M/C_\Psi^\circ)^2n^2\rho_nh_*^2\ge\kappa_n^{\mathrm{str}}\ell_n$.  If
$r_n\ge h_*$, then
\[
 \inf_{\widetilde\lambda}
 \sup_{\vartheta}\,\E_{P_\vartheta}
 d_{\rm loc}\{\widetilde\lambda,\lambda(\vartheta)\}\ \ge\ \frac12,
\]
the infimum running over all measurable decisions based on
$(A,\widehat R,r_n)$, with loss one for the no-selection decision.  The
obstruction is saturated: it applies only when $r_n\ge h_*$ and does not give a
lower bound proportional to $\min\{1,r_n/h_*\}$ at every radius.
\end{proposition}

The proof is a two-point argument: a block swap exchanging $I_0$ and $I_1$
produces two parameter values with different canonical atoms, identical
edge-probability matrices, and a common reporting kernel of radius exactly $h_*$,
so the observable laws coincide while the localization losses sum to at least
one.  Thus large coordinate error destroys canonical location before it destroys
evidence that some departure exists: detection remains achievable at these
amplitudes while location does not.  Canonical \emph{sign} is not protected
either, and at a smaller radius.  Take an off-diagonal atom on blocks $I,J$,
exchange the two halves of $I$ while leaving $J$ fixed, and reverse the sign of
the coefficient: the exchange reverses the atom and the sign reversal restores
the same edge-probability matrix exactly, at a displacement of only $h_*/2$.
Reporting the original ranks in both configurations gives identical observable
laws with opposite canonical signs, so no estimator recovers the sign with
probability above one half uniformly over the reporting-kernel class.
Detection is not protected over the full class either: permuting ranks
arbitrarily within each of the two blocks, still within radius $h_*$, and
reporting the identity order turns the checkerboard into
$p_{ij}=\tfrac12+\delta_ns_it_j$ with unknown balanced sign vectors, and a
mixture over those signs is indistinguishable from the null when
$a_n^2H\to0$ even though $a_n^2H^2$ may exceed $\log m_n$.  The detection
statement of this paper is therefore confined to the block-swap
subexperiment.  A detection theorem under the full reporting class requires
restricting the kernels, and is not given.  Sign
recovery in the simulations of Section~\ref{sec:sim} is therefore a
statement about those particular perturbation mechanisms, which preserve signed
lobes, and not a general one.  A theorem on sign requires restricting the
reporting kernels to lobe-preserving ones.  The transfer clauses of
Theorem~\ref{thm:order} remain conditional on score, studentizer, covariance, and
critical-value differences supplied by the estimator.  Bounded rank displacement
alone does not control them.

\subsubsection{All-radius floor}

Proposition~\ref{prop:obstruct} is the endpoint $r_n\ge h$ of a floor that
exists at every radius.  To see it, leave the aligned dictionary for the
\emph{continuous-coordinate experiment}: coordinates $U_i$ are i.i.d.\ from a
density $g$ on $[0,1]$ known only to satisfy $\tfrac34\le g\le\tfrac32$.  The
departure is the width-$h$ symmetric Haar bump $a h\Psi^S_{h,x,y}$ whose two blocks start at continuous positions $x\le y$ (the symmetric departure identifies only the unordered pair, so the canonical location is the quotient point $\{x,y\}$).  Reports satisfy
$\max_i|\widehat U_i-U_i|\le r_n$ with the joint law of the reports given $U$
otherwise unrestricted (ranks-only observation included).  The canonical
location is the quotient point $\{x,y\}$ with the matching loss
$d(\{x,y\},\{x',y'\})=\min$ over the two pairings of
$\max(|x-x'|,|y-y'|)$, which is zero exactly when the unordered pairs agree.

\begin{theorem}[All-radius floors for canonical coordinate localization]
\label{thm:floor}
\textup{(a)} For every $t\le r_n$ and every configuration whose blocks sit
inside coarse cells with room $7t$ to the cell boundaries and either coincide or
are $h+13t$ apart, there is a second configuration in the class whose canonical
location is $(x+t,y+t)$ and reporting kernels of radius $t$ under which
$(A,\widehat U)$ has \emph{identical} joint law.  Hence at every amplitude and
every $n$, every estimator has worst-case error $d\ge t/2$ with probability at
least one half, and every confidence region with coverage $1-\alpha$, for $\alpha<1/2$, has
worst-case expected diameter at least $(1-2\alpha)t$.  The translation is capped at $t\le\min\{r_n,(2^{-J_0}-h)/13\}$, so under the
class constraint $h\le2^{-J_0-1}$ the floor is of order $\min\{r_n,2^{-J_0}\}$
with explicit constant and continuous in $r_n$.
\textup{(b)} Even at $r_n=0$ with known design and background, and averaging
over the i.i.d.\ coordinate design (the bound is not asserted conditionally on
an arbitrary fixed design, where a translation may cross a design point),
translating one block by $\varepsilon\le h$ changes the law by
$\operatorname{KL}\le C_Ba^2n^2\rho_nh\varepsilon$, so under one of two configurations every estimator errs, in the quotient loss, by at least a quarter of
$\varepsilon_n\asymp\min\{h,(a^2n^2\rho_nh)^{-1}\}=h\min\{1,1/(S_n\log 2m_n)\}$
with probability at least one quarter (the constant is a quarter and not a half
because the quotient distance between $\{x,y\}$ and $\{x+\varepsilon,y\}$ can be
as small as $\varepsilon/2$ when the blocks overlap, so an estimator midway
between them errs by only $\varepsilon/4$), $S_n=a^2n^2\rho_nh^2/\log(2m_n)$.
\end{theorem}

\paragraph{Coordinate-metric and exact-atom losses.}
\label{sec:theory-floor}

Theorem~\ref{thm:floor}(a) is proved for translations $t\le r_n$ and gives a
floor of order $r_n$ in the coordinate metric at every positive radius,
however small.  What it does not settle is recovery of a discrete atom label:
a coordinate error smaller than the atom width can leave the label unchanged,
so the coordinate-metric floor and the exact-atom loss the simulations use
are different quantities.  Section~\ref{sec:sim} measures the latter, and
this subsection keeps the two apart.

For $r_n\ge h$ the obstruction is complete and proved: the block-swap
construction of Proposition~\ref{prop:obstruct} exhibits two configurations
with identical joint law whose canonical supports differ, so no estimator
attains canonical accuracy better than chance between them, at any amplitude.
The measured analogue is exact canonical localization of $0.46$, $0.46$ and
$0.22$ at $\kappa=6,12,24$ under the block-swap kernel: never above the
one-half bound, and not a demonstration that chance recovery is attained
independently of amplitude.

For small $r_n$ the construction is not empty.  The coupling needs room
proportional to the translation, so a smaller translation needs less room, and
part (a) delivers observational equivalence at every positive radius: in the
coordinate metric the floor is of order $\min\{r_n,2^{-J_0}\}$ for all $r_n>0$,
however small.  What the simulations of Section~\ref{sec:sim} measure is a
different quantity, exact recovery of a discrete atom label under one
perturbation mechanism, and they find it high, $0.79$--$1.00$, at realized
$r_n/h$ up to $0.37$ across a sixteenfold range of amplitude.  The two
statements are compatible because the losses differ: a coordinate error of
order $r_n$ that leaves the atom label unchanged is a large loss in the metric
of Theorem~\ref{thm:floor} and no loss at all in the exact-atom metric.  We
draw no identification conclusion from the simulated curve.  It is an empirical
pattern for that mechanism, and identification for the exact-atom
loss in $0<r_n<h$ remains open.

In the exact-atom loss, under that mechanism, recovery falls from about
$0.94$ at $0.37h$ to $0.685$ at $0.70h$, $0.281$ at $1.21h$ and $0.000$ at
$2.17h$.  A matched upper and lower bound for the exact-atom loss in one
experiment is the principal open problem left by this paper.  The measured
curve is what such a theorem would have to reproduce for that mechanism, and
it constrains the theorem without establishing it.  The admission constant $c_r=0.25$ is conservative relative to the empirical
transition under this mechanism, by a factor near $1.5$, and forgoes power
there.  That observation does not license relaxing it, because a defensible
rule must depend on the information and perturbation conditions of
Theorem~\ref{thm:order}, not on one measured curve (Section~\ref{sec:sim}).

The proof of (a) is a coupling: a monotone map that shifts the two blocks by
$t$, ramps back to the identity inside the same coarse cells, and pushes the
uniform design to a density in $[\tfrac34,\tfrac32]$.  The coarse background is
invariant because the map preserves coarse cells, the bump is carried exactly
onto its translate, and reporting the untransformed seed is valid in both
worlds.  The complete construction and constants are in \suppsec{7}{sec:supp-order}.  Two ingredients that the block swap does not need, a translated
departure and a design-density nuisance, are what make the floor continuous.
with a fixed known design and aligned atoms the same argument stops at $t=h/2$,
and for that class and the exact-atom loss the regime $0<r_n<h$ remains open
between Theorem~\ref{thm:order} and Proposition~\ref{prop:obstruct}.  The signal-limited floor transfers to every radius, because the identity
reporting kernel is admissible at any $r_n$ and reproduces the $r_n=0$ data,
so the two floors hold within one experiment.  Read together, (a) and (b)
give the joint phase picture in the coordinate metric:
canonical localization error is at least of order
$\max\{\min(r_n,2^{-J_0}),\ \min(h,(a^2n^2\rho_nh)^{-1})\}$, an order-limited
floor that no amplitude removes, and a signal-limited floor that no coordinate
accuracy removes.  The aligned scan localizes to the atom width $h$.  Whether a
finer scan reaches the signal-limited rate is left open.

The middle clause is where the intermediate regime lives: the far-field
maximum consists of null noise plus the displacement artifacts a fixed coarse
projection leaves along coarse-cell boundaries, of energy order
$J^2n^2\rho_nh\,r_n$ for a boundary jump $J$, which are supercritical once
$r_n/h\gtrsim\log(2m_n)/(J^2n^2\rho_nh^2)$, the manufacture regime of
Section~\ref{sec:sim}.  Bounding them requires the far-field score transfer
that Theorem~\ref{thm:order} assumes and does not derive.

\subsubsection{A feasibility check}

Two quantities decide whether a candidate coordinate and network can support a
claim about a stretch of width $h$.

\begin{proposition}[Rank identifiability under measurement error]
\label{prop:packing}
Let $n$ observed score centres lie in an interval of length $\Xi$ and satisfy
the occupancy condition that every subinterval of length $\ell$ contains at
most $Cn\ell/\Xi+1$ of them, and suppose each true score is covered
simultaneously by an interval of half-width $z\sigma_i$ about its centre with
$\sigma_i\le\sigma$.  Then on the coverage event the maximum normalized rank
displacement satisfies $r_{\max}\le4Cz\sigma/\Xi+1/n$, so a stretch of width
$h$ is admissible whenever $\sigma/\Xi\le(c_rh-1/n)/(4Cz)$.
\end{proposition}

The condition is sufficient, not necessary: the realized radius can be much
smaller than the bound.  At $c_r=1/4$, $C=1$ and $h=1/8$ it requires
measurement error below about $h/(16z)$ of the coordinate range, which at
$n\approx400$ is roughly two parts in a thousand.

The second requirement is the information condition \eqref{eq:inforule} of
Theorem~\ref{thm:info}.  Both quantities are computable from a coordinate, its
standard errors, a vertex count and an edge count, before any model is fitted.

\subsubsection{Coordinates the check refuses}

The two conditions are restrictive, and it is worth recording what they
exclude.  A friendship network of $3000$ users ordered by a home location
estimated from check-ins has a simultaneous rank radius of $0.304$ against the
$0.031$ that the finest scanned width requires, and the procedure admits no
scale.  A legislative cosponsorship network of $421$ members ordered by a
roll-call score has a radius of $0.401$: the scores are precise, with a median
standard error of $0.0093$, but the members are packed nineteen times more
tightly than that along the axis, at a median gap of $0.0005$, so their ranks
are not identified even though their scores are.  Both refusals come from
Proposition~\ref{prop:packing} and both are confirmed by running the
procedure.

By \eqref{eq:inforule} the same two networks resolve nothing narrower than a
tenth of their axis, as against $0.065$ for the ribosomal chain and $0.023$
for the Hi-C window.  Estimated social coordinates fail both conditions at the
sample sizes at which they are available, and that is why the applications
here are physical measurements.

\section{Dictionary properties}\label{sec:supp-dict}

Throughout, $\phi_a(u)=\prod_{k\le d}\psi_{j,l_k}(u_k)$ for a dyadic box $a$ at
scale $j$, and atoms are as in \eqref{eq:dict}.  Write $\langle\cdot,\cdot\rangle$
for the $L^2([0,1]^d\times[0,1]^d)$ inner product.

\begin{proof}[Proof of Proposition~\ref{prop:dict}(i)]
The one-dimensional Haar system $\{\psi_{j,l}\}$ together with the constant is
an orthonormal basis of $L^2[0,1]$, so the tensor products
$\{\phi_a\}\cup\{\text{coarse indicators}\}$ form an orthonormal basis of
$L^2([0,1]^d)$, and the products $\phi_a\otimes\phi_b$ an orthonormal basis of
$L^2([0,1]^d\times[0,1]^d)$.  Symmetrizing pairs $(a,b)$ and $(b,a)$ replaces
the two basis elements $\phi_a\otimes\phi_b$ and $\phi_b\otimes\phi_a$ by their
symmetric and antisymmetric combinations; the symmetric one is $\Psi_{a,b}$ of
\eqref{eq:dict}, with the normalization chosen so that
$\|\Psi_{a,b}\|_2=1$, and the antisymmetric one is orthogonal to every
symmetric function and so contributes nothing.  Restricting to the symmetric
subspace therefore leaves $\{\Psi_{a,b}\}_{a\le b}$ together with the coarse
symmetric indicators, which is an orthonormal basis of that subspace.  A
symmetric $\delta$ with $\langle\delta,\Psi_{a,b}\rangle=0$ for all $a,b$ has
zero component outside the coarse span, hence lies in $V_{J_0}$.
\end{proof}

\begin{proof}[Proof of Proposition~\ref{prop:dict}(ii)]
By bilinearity it suffices to treat $m(u)$, the term $m(v)$ following by
symmetry.  With $a\ne b$,
\[
\langle m(u),\phi_a(u)\phi_b(v)\rangle
=\Bigl(\int m\,\phi_a\Bigr)\Bigl(\int \phi_b\Bigr)=0,
\]
because every Haar function integrates to zero, so $\int\phi_b=0$; the same
holds for the transposed term, and hence for $\Psi_{a,b}$.  For $a=b$ the two
terms coincide and the identical computation applies.  The conclusion does not
require $m$ to be smooth or bounded beyond square integrability, and it is an
exact identity rather than an approximation.
\end{proof}

\begin{proof}[Proof of Proposition~\ref{prop:dict}(iii)]
Fix the scale $h=2^{-j}$.  The support of $\Psi_{a,b}$ is
$(a\times b)\cup(b\times a)$, a union of two boxes of measure $h^{2d}$ each, and
these supports are disjoint for distinct unordered pairs $\{a,b\}$.  If
$\langle\delta,\Psi_{a,b}\rangle\ne0$ then $\delta$ is not almost everywhere
zero on that support, so the support meets $S$.  The boxes meeting $S$ number
at most $A/h^{2d}$ plus the boxes meeting the boundary of $S$, of which there
are $O(h^{-d}|\partial S|)$ in the sense of Minkowski content, and each
unordered pair contributes two boxes, giving the stated count.  For the energy
bound, orthonormality and Bessel's inequality give
\[
\sum_{a,b}\langle\delta,\Psi_{a,b}\rangle^2\le\|\delta\|_2^2
=\int_{S}\delta^2\le B^2A .
\]
\end{proof}

\begin{remark}
Part (iii) is what distinguishes a dictionary scan from surface estimation.
Estimating $f$ at resolution $h$ costs $h^{-2d}$ parameters whatever $f$ is.
A departure of support measure $A$ has $O(A/h^{2d})$ nonzero coefficients, a
fraction $O(A)$ of the dictionary, and the scan pays for the unused remainder
only through the $\sqrt{2\log 2m_n}$ multiplicity factor.  When $A\to0$ the
advantage is unbounded, and when $A\asymp1$ the departure is not localized and
the dictionary offers nothing over a global test.
\end{remark}

\section{Finite-design reduction}\label{sec:supp-finite}


\paragraph{Uniform conditional stochastic order.}
All limits below are as $n\to\infty$. Conditional probabilities may
be taken to be arbitrary regular versions. A quantity defined only
on $\mathcal G_{n,P}$ is assigned an arbitrary measurable value on
$\mathcal G_{n,P}^c$. We adopt the convention that the maximum over
an empty selected set is zero.

For each $P\in\mathcal P_n$, let
$\mathcal G_{n,P}\in\mathcal F_{0n}$. For random variables
$X_{n,P}$ and positive $\mathcal F_{0n}$-measurable scales
$r_{n,P}$, write
\[
X_{n,P}
=
o_{\mathbb P\mid\mathcal F_{0n}}(r_{n,P})
\quad\text{uniformly on }\mathcal G_{n,P}
\]
if, for every $\varepsilon,\delta>0$,
\begin{equation}
\sup_{P\in\mathcal P_n}
P_P\!\left[
\mathcal G_{n,P}\cap
\left\{
P_P\!\left(
|X_{n,P}|>\varepsilon r_{n,P}
\mid\mathcal F_{0n}
\right)>\delta
\right\}
\right]
\longrightarrow0.
\label{eq:uniform-conditional-order}
\end{equation}
For an $\mathcal F_{0n}$-measurable variable $Z_{n,P}$,
$Z_{n,P}=o_{\mathbb P}(1)$ uniformly on $\mathcal G_{n,P}$
means that, for every $\varepsilon>0$,
\[
\sup_{P\in\mathcal P_n}
P_P\!\left(
\mathcal G_{n,P}\cap\{|Z_{n,P}|>\varepsilon\}
\right)
\longrightarrow0.
\]

\begin{theorem}[Exact conditional finite-design reduction]\label{thm:linear}
\label{thm:exact-finite-design-reduction}
Let
\[
\mathcal K_n=\{(i,j):1\le i<j\le n\}.
\]
For $F\in\{I,B^E\}$, let
$\mathcal K_n^F\subseteq\mathcal K_n$ be
$\mathcal F_{0n}$-measurable fold sets satisfying
\begin{equation}
\mathcal K_n^I\cap\mathcal K_n^{B^E}=\varnothing.
\label{eq:disjoint-evaluation-folds}
\end{equation}
For $k=(i,j)$, let $q(k)$ be the
$\mathcal F_{0n}$-measurable symmetric cell containing
$(\widehat U_i,\widehat U_j)$. Cell boundaries are assigned using a
fixed half-open-cell convention, together with deterministic
tie-breaking at the boundary of $[0,1]^2$.

Every coordinate estimate, fold, cell, weight, denominator,
propensity estimate, dictionary coefficient, tuning decision,
admission decision, selected atom, stabilization weight, and frozen
operator appearing below is $\mathcal F_{0n}$-measurable and is
constructed without using the outcomes
\[
\{A_k:k\in\mathcal K_n^I\cup\mathcal K_n^{B^E}\}.
\]
The following conditional product law is imposed directly. For
$P$-almost every realization of $\mathcal F_{0n}$, after fixing the
realized fold sets, the family
\[
\{A_k:k\in\mathcal K_n^I\cup\mathcal K_n^{B^E}\}
\]
is conditionally mutually independent and
\begin{equation}
A_k\mid\mathcal F_{0n}
\sim\operatorname{Bernoulli}(p_k),
\qquad
0<p_k<1,
\label{eq:conditional-product-law-final}
\end{equation}
where $p_k$ is $\mathcal F_{0n}$-measurable. This product-law
assumption is stronger than an informal marginal
missing-at-random condition.

Let $\mathcal Q_n^o$ be the union of cells supporting admitted
atoms, set $\mathcal Q_n^I=\mathcal Q_n^o$, and let
$\mathcal Q_n^{B^E}$ contain every supported baseline-fold cell
used as an input by the frozen baseline operator. Define
\[
\mathcal A_n
=
\{(F,r):F\in\{I,B^E\},\ r\in\mathcal Q_n^F\},
\qquad
D_n^{\mathrm{act}}=|\mathcal A_n|.
\]
Suppose there exist
$\mathcal G_{n,P}\in\mathcal F_{0n}$ satisfying
\begin{equation}
\sup_{P\in\mathcal P_n}
P_P(\mathcal G_{n,P}^c)\longrightarrow0
\label{eq:uniform-good-event}
\end{equation}
such that, on $\mathcal G_{n,P}$,
\begin{equation}
D_n^{\mathrm{act}}\le2Q_n\le n^2,
\qquad
m_n:=|\mathcal I_n|\ge1,
\qquad
\bar m_n:=m_n\vee2.
\label{eq:active-dimension-final}
\end{equation}

For $F\in\{I,B^E\}$, let
$0\le W_k^F<\infty$, $k\in\mathcal K_n^F$, be
$\mathcal F_{0n}$-measurable. For
$r\in\mathcal Q_n^F$, define
\begin{equation}
D_r^F
=
\sum_{k\in\mathcal K_n^F}
W_k^F\mathbf 1\{q(k)=r\},
\qquad
b_{kr}^F
=
\frac{
W_k^F\mathbf 1\{q(k)=r\}
}{D_r^F},
\label{eq:normalized-cell-weights-final}
\end{equation}
and set $b_{kr}^F=0$ when $k\notin\mathcal K_n^F$.
Assume, on $\mathcal G_{n,P}$, that
\[
0<D_r^F<\infty
\qquad
\text{for every }(F,r)\in\mathcal A_n.
\]
Consequently,
\[
b_{kr}^F\ge0,
\qquad
\sum_{k\in\mathcal K_n^F}b_{kr}^F=1.
\]
Define
\begin{equation}
Y_r^F
=
\sum_{k\in\mathcal K_n^F}b_{kr}^FA_k,
\qquad
\mu_r^F
=
E_P(Y_r^F\mid\mathcal F_{0n})
=
\sum_{k\in\mathcal K_n^F}b_{kr}^Fp_k.
\label{eq:cell-means-final}
\end{equation}
Then, exactly,
\begin{equation}
Y_r^F-\mu_r^F
=
\sum_{k\in\mathcal K_n^F}
b_{kr}^F(A_k-p_k).
\label{eq:centered-cell-means-final}
\end{equation}

Let
\[
B_n:
\mathbb R^{|\mathcal Q_n^{B^E}|}
\longrightarrow
(0,1)^{|\mathcal Q_n^o|}
\]
be the complete $\mathcal F_{0n}$-measurable baseline operator,
including every smoothing, stabilization, projection, inverse-link,
thresholding, and clipping operation used by the implemented
estimator. Let $w_q\in[0,1]$ and $c_{\lambda q}$ be
$\mathcal F_{0n}$-measurable. Assume
\[
0<\epsilon_n<\frac12
\]
and define
\[
\ell_{\epsilon_n}(t)
=
\operatorname{logit}
\{\min(1-\epsilon_n,\max(\epsilon_n,t))\}.
\]
For $x=(x_q:q\in\mathcal Q_n^I)$ and
$y=(y_r:r\in\mathcal Q_n^{B^E})$, put
\begin{equation}
z_q(x,y)
=
w_qx_q+(1-w_q)B_{n,q}(y)
\label{eq:stabilized-map-final}
\end{equation}
and
\begin{equation}
\Phi_{\lambda n}(x,y)
=
\sum_{q\in\mathcal Q_n^o}
c_{\lambda q}
\left[
\ell_{\epsilon_n}\{z_q(x,y)\}
-
\operatorname{logit}\{B_{n,q}(y)\}
\right].
\label{eq:frozen-functional-final}
\end{equation}
Assume the following exact algorithmic identity:
\begin{equation}
\widehat\theta_\lambda
=
\Phi_{\lambda n}(Y^I,Y^{B^E}).
\label{eq:exact-estimator-representation-final}
\end{equation}
Define the frozen-operator finite-design target by
\begin{equation}
\theta_{\lambda,n}^D
=
\Phi_{\lambda n}(\mu^I,\mu^{B^E}).
\label{eq:finite-design-target-final}
\end{equation}
Because $\Phi_{\lambda n}$ is nonlinear, in general
\[
\theta_{\lambda,n}^D
\ne
E_P(\widehat\theta_\lambda\mid\mathcal F_{0n}).
\]

Put
\[
\ell_n=\log(e\bar m_n),
\qquad
\beta_n=4\log(e\bar m_nn).
\]
The Bernstein exponent is written $\beta_n$ throughout; in earlier drafts it
carried the symbol $\gamma_n$, which is reserved here for the
admission-failure probability of the detection theorem, a distinct and
$o(1)$ quantity.
For $\alpha=(F,r)\in\mathcal A_n$, define
\begin{equation}
V_\alpha
=
\sum_{k\in\mathcal K_n^F}
(b_{kr}^F)^2p_k(1-p_k),
\qquad
M_\alpha
=
\max_{k\in\mathcal K_n^F}|b_{kr}^F|,
\label{eq:cell-variance-envelope-final}
\end{equation}
and
\begin{equation}
t_\alpha
=
\sqrt{2V_\alpha\beta_n}
+
\frac{2}{3}M_\alpha\beta_n.
\label{eq:bernstein-radius-final}
\end{equation}
Let
\begin{equation}
\mathcal C_{n,P}
=
\prod_{\alpha\in\mathcal A_n}
\left(
[\mu_\alpha-t_\alpha,\mu_\alpha+t_\alpha]\cap[0,1]
\right).
\label{eq:bernstein-rectangle-final}
\end{equation}

Assume that, on $\mathcal G_{n,P}$, there are an open convex set
$\mathcal U_{n,P}\supset\mathcal C_{n,P}$ and an
$\mathcal F_{0n}$-measurable $\kappa_{n,P}>0$ such that the complete
operator $B_n$, including every internal gate and clipping
operation, is twice continuously differentiable on the
$y$-projection of $\mathcal U_{n,P}$ and
\begin{equation}
\inf_{(x,y)\in\mathcal U_{n,P}}
\min_{q\in\mathcal Q_n^o}
\operatorname{dist}
\bigl(z_q(x,y),\{\epsilon_n,1-\epsilon_n\}\bigr)
\ge\kappa_{n,P}
\label{eq:no-knot-condition-final}
\end{equation}
and
\begin{equation}
\kappa_{n,P}
\le
B_{n,q}(y)
\le
1-\kappa_{n,P}
\quad
\text{for every }q\in\mathcal Q_n^o
\text{ and every }y\in\operatorname{proj}_y(\mathcal U_{n,P}).
\label{eq:baseline-interior-condition-final}
\end{equation}
Assume additionally that the random maps $B_n$, their first two
derivatives, and the corresponding maps $\Phi_{\lambda n}$ are
jointly measurable in the sample point and their Euclidean
arguments. In particular, all suprema defined below are measurable.

Define
\begin{equation}
\pi_q^0
=
B_{n,q}(\mu^{B^E}),
\qquad
z_q^0
=
w_q\mu_q^I+(1-w_q)\pi_q^0,
\qquad
H_{qr}^0
=
\left.
\frac{\partial B_{n,q}(y)}{\partial y_r}
\right|_{y=\mu^{B^E}}.
\label{eq:oracle-map-quantities-final}
\end{equation}
The gradient blocks of $\Phi_{\lambda n}$ at
$(\mu^I,\mu^{B^E})$ are
\begin{equation}
G_{\lambda r}^{I,0}
=
c_{\lambda r}w_r
\ell_{\epsilon_n}'(z_r^0),
\qquad
r\in\mathcal Q_n^I,
\label{eq:oracle-gradient-I-final}
\end{equation}
and
\begin{equation}
G_{\lambda r}^{B^E,0}
=
\sum_{q\in\mathcal Q_n^o}
c_{\lambda q}
\left[
(1-w_q)\ell_{\epsilon_n}'(z_q^0)
-
\frac{1}{\pi_q^0(1-\pi_q^0)}
\right]H_{qr}^0.
\label{eq:oracle-gradient-BE-final}
\end{equation}

For $F\in\{I,B^E\}$ and $k\in\mathcal K_n$, define
\begin{equation}
a_{\lambda k}^F
=
\sum_{r\in\mathcal Q_n^F}
G_{\lambda r}^{F,0}b_{kr}^F,
\qquad
a_{\lambda k}
=
a_{\lambda k}^I+a_{\lambda k}^{B^E}.
\label{eq:oracle-loadings-final}
\end{equation}
By fold disjointness,
$a_{\lambda k}^Ia_{\lambda k}^{B^E}=0$. Define
\begin{equation}
s_\lambda^2
=
\sum_{F\in\{I,B^E\}}
\sum_{k\in\mathcal K_n^F}
(a_{\lambda k}^F)^2p_k(1-p_k),
\qquad
I_\lambda^{\mathrm{or}}=s_\lambda^{-2},
\label{eq:oracle-variance-final}
\end{equation}
and assume $0<s_\lambda<\infty$ for every admitted atom.
Here $I_\lambda^{\mathrm{or}}$ is effective oracle information,
namely the inverse variance of the oracle linear term; it is not
being identified with likelihood Fisher information.

Define
\begin{equation}
\zeta_{\lambda,k}
=
\begin{cases}
\dfrac{a_{\lambda k}}{s_\lambda}(A_k-p_k),
&
k\in\mathcal K_n^I\cup\mathcal K_n^{B^E},
\\[6pt]
0,
&
k\notin\mathcal K_n^I\cup\mathcal K_n^{B^E},
\end{cases}
\label{eq:oracle-influence-final}
\end{equation}
and
\begin{equation}
\Lambda_n^{\mathrm{or}}
=
\max_{\lambda\in\mathcal I_n}
\max_{k\in\mathcal K_n^I\cup\mathcal K_n^{B^E}}
\frac{|a_{\lambda k}|}{s_\lambda}.
\label{eq:oracle-leverage-final}
\end{equation}
Then, for almost every realization of $\mathcal F_{0n}$ in
$\mathcal G_{n,P}$,
\begin{equation}
E_P(\zeta_{\lambda,k}\mid\mathcal F_{0n})=0,
\qquad
\sum_{k\in\mathcal K_n}
E_P(\zeta_{\lambda,k}^2\mid\mathcal F_{0n})=1,
\qquad
\max_{\lambda,k}|\zeta_{\lambda,k}|
\le\Lambda_n^{\mathrm{or}}.
\label{eq:influence-properties-final}
\end{equation}

Define
\begin{equation}
\Gamma_{\lambda,n}
=
\frac{1}{2s_\lambda}
\sup_{v\in\mathcal C_{n,P}}
\sup_{\substack{
u\in\mathbb R^{D_n^{\mathrm{act}}}\\
|u_\alpha|\le t_\alpha\ \forall\alpha}}
\left|
u^\top\nabla^2\Phi_{\lambda n}(v)u
\right|.
\label{eq:curvature-envelope-final}
\end{equation}
Then, for every $P\in\mathcal P_n$ and almost every realization of
$\mathcal F_{0n}$ in $\mathcal G_{n,P}$,
\begin{equation}
\begin{split}
&P_P\Bigg[
\max_{\lambda\in\mathcal I_n}
\left|
\frac{\widehat\theta_\lambda-\theta_{\lambda,n}^D}{s_\lambda}
-
\sum_{k\in\mathcal K_n}\zeta_{\lambda,k}
\right|
\le
\max_{\lambda\in\mathcal I_n}\Gamma_{\lambda,n}
\ \Bigg|\ \mathcal F_{0n}
\Bigg]\\
&\hspace{4cm}
\ge
1-2D_n^{\mathrm{act}}e^{-\beta_n}.
\end{split}
\label{eq:exact-finite-sample-bound-final}
\end{equation}
\end{theorem}

\begin{proof}
Fix $P\in\mathcal P_n$ and condition on a realization of
$\mathcal F_{0n}$ in $\mathcal G_{n,P}$. All folds, cells, weights,
denominators, selected atoms, coefficients, and maps are fixed.

For $\alpha=(F,r)\in\mathcal A_n$,
\eqref{eq:centered-cell-means-final} gives
\[
Y_\alpha-\mu_\alpha
=
\sum_{k\in\mathcal K_n^F}
b_{kr}^F(A_k-p_k).
\]
The summands are conditionally independent and centered. Their
conditional variance is $V_\alpha$, and each summand has absolute
value at most $M_\alpha$. Bernstein's inequality gives
\begin{equation}
P_P\!\left(
|Y_\alpha-\mu_\alpha|>t_\alpha
\mid\mathcal F_{0n}
\right)
\le2e^{-\beta_n}.
\label{eq:coordinate-bernstein-final}
\end{equation}
A union bound therefore yields
\begin{equation}
P_P\!\left(
(Y^I,Y^{B^E})\notin\mathcal C_{n,P}
\mid\mathcal F_{0n}
\right)
\le
2D_n^{\mathrm{act}}e^{-\beta_n}.
\label{eq:rectangle-probability-final}
\end{equation}

On the complementary event, write
$Y=(Y^I,Y^{B^E})$ and
$\mu=(\mu^I,\mu^{B^E})$. Convexity gives
$\mu+t(Y-\mu)\in\mathcal C_{n,P}$ for every $t\in[0,1]$.
Taylor's theorem with integral remainder gives
\begin{equation}
\begin{split}
\widehat\theta_\lambda-\theta_{\lambda,n}^D
={}&
\nabla\Phi_{\lambda n}(\mu)^\top(Y-\mu)\\
&+
\int_0^1(1-t)
(Y-\mu)^\top
\nabla^2\Phi_{\lambda n}\{\mu+t(Y-\mu)\}
(Y-\mu)\,dt .
\end{split}
\label{eq:taylor-integral-final}
\end{equation}

Direct differentiation gives
\eqref{eq:oracle-gradient-I-final} and
\eqref{eq:oracle-gradient-BE-final}. Consequently,
\begin{equation}
\begin{split}
\nabla\Phi_{\lambda n}(\mu)^\top(Y-\mu)
&=
\sum_{F\in\{I,B^E\}}
\sum_{r\in\mathcal Q_n^F}
G_{\lambda r}^{F,0}(Y_r^F-\mu_r^F)\\
&=
\sum_{F\in\{I,B^E\}}
\sum_{k\in\mathcal K_n^F}
a_{\lambda k}^F(A_k-p_k)\\
&=
s_\lambda
\sum_{k\in\mathcal K_n}\zeta_{\lambda,k}.
\end{split}
\label{eq:linear-term-final}
\end{equation}

On $\{Y\in\mathcal C_{n,P}\}$,
$|Y_\alpha-\mu_\alpha|\le t_\alpha$ for every $\alpha$. Hence
\begin{equation}
\begin{split}
&\left|
\int_0^1(1-t)
(Y-\mu)^\top
\nabla^2\Phi_{\lambda n}\{\mu+t(Y-\mu)\}
(Y-\mu)\,dt
\right|\\
&\quad\le
\frac12
\sup_{v\in\mathcal C_{n,P}}
\sup_{|u_\alpha|\le t_\alpha}
|u^\top\nabla^2\Phi_{\lambda n}(v)u|
=
s_\lambda\Gamma_{\lambda,n}.
\end{split}
\label{eq:remainder-bound-final}
\end{equation}
Thus, simultaneously in $\lambda$,
\[
\left|
\frac{\widehat\theta_\lambda-\theta_{\lambda,n}^D}{s_\lambda}
-
\sum_{k\in\mathcal K_n}\zeta_{\lambda,k}
\right|
\le\Gamma_{\lambda,n}.
\]
Together with \eqref{eq:rectangle-probability-final}, this proves
\eqref{eq:exact-finite-sample-bound-final}.

For $k\in\mathcal K_n^I\cup\mathcal K_n^{B^E}$,
\[
E_P(A_k-p_k\mid\mathcal F_{0n})=0,
\qquad
E_P\{(A_k-p_k)^2\mid\mathcal F_{0n}\}
=p_k(1-p_k).
\]
Therefore the influence variables are conditionally centered and
\[
\sum_{k\in\mathcal K_n}
E_P(\zeta_{\lambda,k}^2\mid\mathcal F_{0n})
=
\frac{1}{s_\lambda^2}
\sum_{F\in\{I,B^E\}}
\sum_{k\in\mathcal K_n^F}
(a_{\lambda k}^F)^2p_k(1-p_k)
=1.
\]
Finally, $|A_k-p_k|\le1$ gives
\[
|\zeta_{\lambda,k}|
\le
\frac{|a_{\lambda k}|}{s_\lambda}
\le\Lambda_n^{\mathrm{or}}.
\]
\end{proof}

\begin{lemma}[Deterministic curvature and cell-leverage bounds]
\label{lem:deterministic-curvature-leverage}
Retain the notation and assumptions of
Theorem~\ref{thm:exact-finite-design-reduction}. For
$\alpha,\beta\in\mathcal A_n$, define
\begin{equation}
\mathfrak H_{\lambda,\alpha\beta,n}
=
\sup_{v\in\mathcal C_{n,P}}
\left|
\partial_{\alpha\beta}^2\Phi_{\lambda n}(v)
\right|
\label{eq:componentwise-hessian-envelope}
\end{equation}
and
\begin{equation}
\overline\Gamma_{\lambda,n}
=
\frac{1}{2s_\lambda}
\sum_{\alpha,\beta\in\mathcal A_n}
\mathfrak H_{\lambda,\alpha\beta,n}
t_\alpha t_\beta.
\label{eq:explicit-curvature-envelope}
\end{equation}
Then
\begin{equation}
\Gamma_{\lambda,n}
\le
\overline\Gamma_{\lambda,n}.
\label{eq:curvature-envelope-comparison}
\end{equation}

For $(F,r)\in\mathcal A_n$, define
\begin{equation}
V_r^{F,0}
:=
V_{(F,r)}
=
\sum_{k\in\mathcal K_n^F}
(b_{kr}^F)^2p_k(1-p_k)>0
\label{eq:positive-oracle-cell-meat}
\end{equation}
and
\begin{equation}
L_n^{\mathrm{cell}}
=
\max_{\substack{
(F,r)\in\mathcal A_n\\
k\in\mathcal K_n^F:\ b_{kr}^F>0}}
\frac{b_{kr}^F}{\sqrt{V_r^{F,0}}}.
\label{eq:cell-leverage-envelope}
\end{equation}
Then
\begin{equation}
\Lambda_n^{\mathrm{or}}
\le
L_n^{\mathrm{cell}}.
\label{eq:oracle-cell-leverage-bound}
\end{equation}
\end{lemma}

\begin{proof}
The rectangle $\mathcal C_{n,P}$ is compact. Since
$\Phi_{\lambda n}$ is twice continuously differentiable on an open
neighborhood of this rectangle, every quantity in
\eqref{eq:componentwise-hessian-envelope} is finite.

For $v\in\mathcal C_{n,P}$ and $|u_\alpha|\le t_\alpha$,
\[
\begin{split}
\left|
u^\top\nabla^2\Phi_{\lambda n}(v)u
\right|
&=
\left|
\sum_{\alpha,\beta\in\mathcal A_n}
u_\alpha u_\beta
\partial_{\alpha\beta}^2\Phi_{\lambda n}(v)
\right|\\
&\le
\sum_{\alpha,\beta\in\mathcal A_n}
|u_\alpha|\,|u_\beta|\,
\left|
\partial_{\alpha\beta}^2\Phi_{\lambda n}(v)
\right|\\
&\le
\sum_{\alpha,\beta\in\mathcal A_n}
t_\alpha t_\beta
\mathfrak H_{\lambda,\alpha\beta,n}.
\end{split}
\]
Taking the two suprema and dividing by $2s_\lambda$ proves
\eqref{eq:curvature-envelope-comparison}.

For every active cell, positivity of $D_r^F$ implies that
$b_{kr}^F>0$ for at least one $k$. Since $0<p_k<1$, it follows that
$V_r^{F,0}>0$. If $k\in\mathcal K_n^F$ and $q(k)=r$, unique cell
membership gives
\[
a_{\lambda k}^F
=
G_{\lambda r}^{F,0}b_{kr}^F.
\]
Moreover,
\begin{equation}
s_\lambda^2
=
\sum_{F\in\{I,B^E\}}
\sum_{r\in\mathcal Q_n^F}
(G_{\lambda r}^{F,0})^2V_r^{F,0}.
\label{eq:cellwise-oracle-variance}
\end{equation}
Thus
\[
s_\lambda^2
\ge
(G_{\lambda r}^{F,0})^2V_r^{F,0}.
\]
If $G_{\lambda r}^{F,0}=0$, then $a_{\lambda k}^F=0$. Otherwise,
\[
\frac{|a_{\lambda k}^F|}{s_\lambda}
=
\frac{|G_{\lambda r}^{F,0}|b_{kr}^F}{s_\lambda}
\le
\frac{b_{kr}^F}{\sqrt{V_r^{F,0}}}.
\]
Because the folds are disjoint, $a_{\lambda k}=a_{\lambda k}^F$
for the unique evaluation fold containing $k$. Taking maxima proves
\eqref{eq:oracle-cell-leverage-bound}.
\end{proof}

\begin{corollary}[Explicit curvature and leverage conditions]
\label{cor:explicit-curvature-leverage}
Suppose, uniformly over $P\in\mathcal P_n$ on
$\mathcal G_{n,P}$, that
\begin{equation}
\sqrt{\ell_n}
\max_{\lambda\in\mathcal I_n}
\overline\Gamma_{\lambda,n}
=
o_{\mathbb P}(1)
\label{eq:explicit-curvature-rate}
\end{equation}
and, for a specified nonnegative deterministic sequence $L_n$,
\begin{equation}
L_n^{\mathrm{cell}}\le L_n.
\label{eq:explicit-cell-leverage-rate}
\end{equation}
Then
\[
\sqrt{\ell_n}
\max_{\lambda\in\mathcal I_n}
\Gamma_{\lambda,n}
=
o_{\mathbb P}(1),
\qquad
\Lambda_n^{\mathrm{or}}\le L_n
\]
uniformly on $\mathcal G_{n,P}$.
\end{corollary}

\begin{proof}
The conclusions follow directly from
\eqref{eq:curvature-envelope-comparison},
\eqref{eq:oracle-cell-leverage-bound},
\eqref{eq:explicit-curvature-rate}, and
\eqref{eq:explicit-cell-leverage-rate}.
\end{proof}

\begin{corollary}[Uniform finite-design asymptotic linearization]
\label{cor:uniform-finite-design-expansion}
Under Theorem~\ref{thm:exact-finite-design-reduction}, assume that,
for every $\varepsilon>0$,
\begin{equation}
\sup_{P\in\mathcal P_n}
P_P\!\left[
\mathcal G_{n,P}\cap
\left\{
\sqrt{\ell_n}
\max_{\lambda\in\mathcal I_n}
\Gamma_{\lambda,n}>\varepsilon
\right\}
\right]
\longrightarrow0.
\label{eq:uniform-curvature-rate-final}
\end{equation}
Then
\begin{equation}
\max_{\lambda\in\mathcal I_n}
\left|
\frac{\widehat\theta_\lambda-\theta_{\lambda,n}^D}{s_\lambda}
-
\sum_{k\in\mathcal K_n}\zeta_{\lambda,k}
\right|
=
o_{\mathbb P\mid\mathcal F_{0n}}(\ell_n^{-1/2})
\label{eq:uniform-finite-design-linearization-final}
\end{equation}
uniformly on $\mathcal G_{n,P}$.

If, in addition, a separate oracle-transfer result gives
\begin{equation}
\Lambda_n^{\mathrm{or}}\le L_n
\quad\text{on }\mathcal G_{n,P},
\label{eq:leverage-transfer-final}
\end{equation}
for a specified nonnegative deterministic sequence $L_n$, then
\[
\max_{\lambda\in\mathcal I_n}
\max_{k\in\mathcal K_n}
|\zeta_{\lambda,k}|
\le L_n
\]
on $\mathcal G_{n,P}$.
\end{corollary}

\begin{proof}
Let $R_{n,P}^D$ denote the maximum on the left of
\eqref{eq:uniform-finite-design-linearization-final}. Then
\[
P_P\!\left(
R_{n,P}^D>\varepsilon\ell_n^{-1/2}
\mid\mathcal F_{0n}
\right)
\le
2D_n^{\mathrm{act}}e^{-\beta_n}
+
\mathbf 1\!\left\{
\sqrt{\ell_n}
\max_{\lambda\in\mathcal I_n}\Gamma_{\lambda,n}
>\varepsilon
\right\}.
\]
On $\mathcal G_{n,P}$,
\[
2D_n^{\mathrm{act}}e^{-\beta_n}
\le
2n^2(e\bar m_nn)^{-4}
\le Cn^{-2}.
\]
For fixed $\delta>0$, the deterministic term is smaller than
$\delta$ for all sufficiently large $n$. Hence the event that the
preceding conditional probability exceeds $\delta$ is contained in
the curvature event in
\eqref{eq:uniform-curvature-rate-final}. This proves
\eqref{eq:uniform-finite-design-linearization-final}. The leverage
claim follows from \eqref{eq:influence-properties-final}.
\end{proof}

\begin{lemma}[Studentizer comparison and explicit sufficient rate]
\label{lem:studentizer-transfer-final}
Retain the notation of
Theorem~\ref{thm:exact-finite-design-reduction}. Put
\[
Y=(Y^I,Y^{B^E}),
\qquad
\mu=(\mu^I,\mu^{B^E}),
\qquad
\mathcal B_{n,P}=\{Y\in\mathcal C_{n,P}\}.
\]
Every sum or maximum over $(F,r)$ below is over
$\mathcal A_n$.

Let $\widehat G_\lambda$ be a finite measurable vector satisfying
\begin{equation}
\widehat G_\lambda
=
\nabla\Phi_{\lambda n}(Y)
\qquad\text{on }\mathcal B_{n,P}.
\label{eq:fitted-gradient-agreement}
\end{equation}
Let $\widehat V_r^F\ge0$ be measurable fitted cell meats and define
$\widehat s_\lambda$ as the nonnegative square root of
\begin{equation}
\widehat s_\lambda^2
=
\sum_{(F,r)\in\mathcal A_n}
(\widehat G_{\lambda r}^F)^2
\widehat V_r^F.
\label{eq:fitted-studentizer-final}
\end{equation}
Define
\begin{equation}
\delta_{V,n}
=
\max_{(F,r)\in\mathcal A_n}
\left|
\frac{\widehat V_r^F}{V_r^{F,0}}-1
\right|
\label{eq:studentizer-meat-error-final}
\end{equation}
and
\begin{equation}
\delta_{G,n}
=
\max_{\lambda\in\mathcal I_n}
\frac{
\left[
\sum_{(F,r)\in\mathcal A_n}
V_r^{F,0}
(\widehat G_{\lambda r}^F-G_{\lambda r}^{F,0})^2
\right]^{1/2}
}{s_\lambda}.
\label{eq:studentizer-gradient-error-final}
\end{equation}

On $\{\delta_{V,n}<1,\delta_{G,n}<1\}$,
\begin{equation}
\sqrt{1-\delta_{V,n}}\,(1-\delta_{G,n})
\le
\frac{\widehat s_\lambda}{s_\lambda}
\le
\sqrt{1+\delta_{V,n}}\,(1+\delta_{G,n})
\label{eq:studentizer-sandwich-final}
\end{equation}
simultaneously for all $\lambda\in\mathcal I_n$.

Suppose additionally that
\begin{equation}
\widehat V_r^F
=
\sum_{k\in\mathcal K_n^F}
(b_{kr}^F)^2\widehat v_k^F,
\qquad
\widehat v_k^F\ge0,
\label{eq:dyad-level-meat-final}
\end{equation}
and put
\[
v_k=p_k(1-p_k).
\]
Define
\begin{equation}
E_{V,n}
=
\max_{\substack{
(F,r)\in\mathcal A_n\\
k\in\mathcal K_n^F:\ b_{kr}^F>0}}
\left|
\frac{\widehat v_k^F}{v_k}-1
\right|
\label{eq:dyad-variance-relative-error}
\end{equation}
and
\begin{equation}
\mathfrak D_{\lambda,n}
=
\frac{1}{s_\lambda}
\sup_{\substack{
v\in\mathcal C_{n,P}\\
u\in\mathbb R^{D_n^{\mathrm{act}}}:
\,|u_\alpha|\le t_\alpha\ \forall\alpha}}
\left[
\sum_{\alpha\in\mathcal A_n}
V_\alpha
\left\{
\bigl(\nabla^2\Phi_{\lambda n}(v)u\bigr)_\alpha
\right\}^2
\right]^{1/2}.
\label{eq:gradient-change-envelope}
\end{equation}
Then, pathwise,
\begin{equation}
\delta_{V,n}\le E_{V,n},
\label{eq:delta-v-explicit-bound}
\end{equation}
and, on $\mathcal B_{n,P}$,
\begin{equation}
\delta_{G,n}
\le
\max_{\lambda\in\mathcal I_n}
\mathfrak D_{\lambda,n}.
\label{eq:delta-g-explicit-bound}
\end{equation}

If, uniformly over $P\in\mathcal P_n$ on
$\mathcal G_{n,P}$,
\begin{equation}
\ell_n
\left\{
E_{V,n}
+
\max_{\lambda\in\mathcal I_n}
\mathfrak D_{\lambda,n}
\right\}
=
o_{\mathbb P\mid\mathcal F_{0n}}(1),
\label{eq:explicit-studentizer-rate}
\end{equation}
then
\begin{equation}
\delta_{V,n}+\delta_{G,n}
=
o_{\mathbb P\mid\mathcal F_{0n}}(\ell_n^{-1})
\label{eq:studentizer-component-rate}
\end{equation}
and
\begin{equation}
\ell_n
\max_{\lambda\in\mathcal I_n}
\left|
\frac{\widehat s_\lambda}{s_\lambda}-1
\right|
=
o_{\mathbb P\mid\mathcal F_{0n}}(1)
\label{eq:studentizer-conclusion-final}
\end{equation}
uniformly on $\mathcal G_{n,P}$.
\end{lemma}

\begin{proof}
For a stacked vector $g=(g_r^F)$, define
\[
\|g\|_0^2
=
\sum_{(F,r)\in\mathcal A_n}
V_r^{F,0}(g_r^F)^2.
\]
By \eqref{eq:cellwise-oracle-variance},
\[
s_\lambda=\|G_\lambda^0\|_0.
\]
The definition of $\delta_{V,n}$ gives
\[
(1-\delta_{V,n})\|\widehat G_\lambda\|_0^2
\le
\widehat s_\lambda^2
\le
(1+\delta_{V,n})\|\widehat G_\lambda\|_0^2.
\]
The triangle and reverse-triangle inequalities yield
\[
s_\lambda(1-\delta_{G,n})
\le
\|\widehat G_\lambda\|_0
\le
s_\lambda(1+\delta_{G,n}).
\]
Combining these inequalities and taking square roots proves
\eqref{eq:studentizer-sandwich-final}.

For each active cell,
\[
\begin{split}
\left|
\frac{\widehat V_r^F}{V_r^{F,0}}-1
\right|
&=
\frac{
\left|
\sum_{k\in\mathcal K_n^F}
(b_{kr}^F)^2(\widehat v_k^F-v_k)
\right|
}{
\sum_{k\in\mathcal K_n^F}
(b_{kr}^F)^2v_k
}\\
&\le
\max_{\substack{
k\in\mathcal K_n^F\\
b_{kr}^F>0}}
\left|
\frac{\widehat v_k^F}{v_k}-1
\right|.
\end{split}
\]
Taking the maximum over active cells proves
\eqref{eq:delta-v-explicit-bound}.

On $\mathcal B_{n,P}$, put $h=Y-\mu$. Convexity gives
$\mu+th\in\mathcal C_{n,P}$ for every $t\in[0,1]$, and
$|h_\alpha|\le t_\alpha$. By
\eqref{eq:fitted-gradient-agreement},
\[
\widehat G_\lambda-G_\lambda^0
=
\int_0^1
\nabla^2\Phi_{\lambda n}(\mu+th)h\,dt.
\]
Minkowski's integral inequality in the weighted Euclidean norm gives
\[
\begin{split}
&
\frac{
\left[
\sum_{\alpha\in\mathcal A_n}
V_\alpha
(\widehat G_{\lambda,\alpha}
-G_{\lambda,\alpha}^0)^2
\right]^{1/2}
}{s_\lambda}\\
&\quad\le
\int_0^1
\frac{
\left[
\sum_{\alpha\in\mathcal A_n}
V_\alpha
\left\{
\bigl(
\nabla^2\Phi_{\lambda n}(\mu+th)h
\bigr)_\alpha
\right\}^2
\right]^{1/2}
}{s_\lambda}\,dt\\
&\quad\le
\mathfrak D_{\lambda,n}.
\end{split}
\]
Taking the maximum over $\lambda$ proves
\eqref{eq:delta-g-explicit-bound}.

By the conditional Bernstein bound,
\[
P_P(\mathcal B_{n,P}^c\mid\mathcal F_{0n})
\le
2D_n^{\mathrm{act}}e^{-\beta_n}.
\]
On $\mathcal G_{n,P}$, this is at most $Cn^{-2}$. Hence, for every
$\varepsilon>0$,
\[
\begin{split}
&P_P\!\left(
\ell_n(\delta_{V,n}+\delta_{G,n})>\varepsilon
\mid\mathcal F_{0n}
\right)\\
&\quad\le
Cn^{-2}
+
P_P\!\left[
\ell_n
\left\{
E_{V,n}
+
\max_{\lambda\in\mathcal I_n}
\mathfrak D_{\lambda,n}
\right\}
>\varepsilon
\ \middle|\ \mathcal F_{0n}
\right].
\end{split}
\]
Condition \eqref{eq:explicit-studentizer-rate} therefore proves
\eqref{eq:studentizer-component-rate}.

On $\{\delta_{V,n}\vee\delta_{G,n}\le1/2\}$,
\eqref{eq:studentizer-sandwich-final} and the elementary inequality
$|\sqrt{1+x}-1|\le C|x|$ give
\[
\max_{\lambda\in\mathcal I_n}
\left|
\frac{\widehat s_\lambda}{s_\lambda}-1
\right|
\le
C(\delta_{V,n}+\delta_{G,n}).
\]
Equation \eqref{eq:studentizer-component-rate} implies that the
complementary event has vanishing conditional probability and proves
\eqref{eq:studentizer-conclusion-final}.
\end{proof}

\begin{corollary}[Transfer to a population coefficient]
\label{cor:population-transfer-final}
Assume Corollary~\ref{cor:uniform-finite-design-expansion}. For every
deterministic dictionary index $\lambda\in\mathfrak D_n$, let
$\theta_\lambda(P)$ be the deterministic intended population
coefficient, where $\mathfrak D_n$ is a deterministic finite master
dictionary satisfying
\[
\mathcal I_n\subseteq\mathfrak D_n
\qquad\text{almost surely}.
\]
Let
$\mathcal E_n\in\mathcal F_{0n}$ be the design event, for example
$\mathcal E_n=\{\Delta_n\le r_n\}$. Suppose
\begin{equation}
\sup_{P\in\mathcal P_n}P_P(\mathcal E_n^c)
\le\eta_n,
\qquad
\eta_n\longrightarrow0,
\label{eq:design-event-final}
\end{equation}
and, for every $\varepsilon>0$,
\begin{equation}
\sup_{P\in\mathcal P_n}
P_P\!\left[
\mathcal G_{n,P}\cap\mathcal E_n\cap
\left\{
\sqrt{\ell_n}
\max_{\lambda\in\mathcal I_n}
\frac{|\theta_{\lambda,n}^D-\theta_\lambda(P)|}{s_\lambda}
>\varepsilon
\right\}
\right]
\longrightarrow0.
\label{eq:population-bias-rate-final}
\end{equation}
Then
\begin{equation}
\max_{\lambda\in\mathcal I_n}
\left|
\frac{\widehat\theta_\lambda-\theta_\lambda(P)}{s_\lambda}
-
\sum_{k\in\mathcal K_n}\zeta_{\lambda,k}
\right|
=
o_{\mathbb P\mid\mathcal F_{0n}}(\ell_n^{-1/2})
\label{eq:population-linearization-final}
\end{equation}
uniformly on $\mathcal G_{n,P}\cap\mathcal E_n$.

Moreover, if $R_{n,P}^{\mathrm{pop}}$ denotes the maximum on the
left of \eqref{eq:population-linearization-final}, then, for every
fixed $\varepsilon>0$,
\begin{equation}
\sup_{P\in\mathcal P_n}
P_P\!\left(
R_{n,P}^{\mathrm{pop}}
>
\varepsilon\ell_n^{-1/2}
\right)
\le
\eta_n+o(1).
\label{eq:unconditional-population-bound-final}
\end{equation}
\end{corollary}

\begin{proof}
For every selected $\lambda$,
\begin{equation}
\begin{split}
\frac{\widehat\theta_\lambda-\theta_\lambda(P)}{s_\lambda}
-
\sum_{k\in\mathcal K_n}\zeta_{\lambda,k}
={}&
\left\{
\frac{\widehat\theta_\lambda-\theta_{\lambda,n}^D}{s_\lambda}
-
\sum_{k\in\mathcal K_n}\zeta_{\lambda,k}
\right\}\\
&+
\frac{\theta_{\lambda,n}^D-\theta_\lambda(P)}{s_\lambda}.
\end{split}
\label{eq:population-decomposition-final}
\end{equation}
The first maximum is
$o_{\mathbb P\mid\mathcal F_{0n}}(\ell_n^{-1/2})$
on $\mathcal G_{n,P}$. The second maximum is
$o_{\mathbb P}(\ell_n^{-1/2})$ on
$\mathcal G_{n,P}\cap\mathcal E_n$ by
\eqref{eq:population-bias-rate-final}. It is
$\mathcal F_{0n}$-measurable, so it also satisfies the corresponding
conditional statement. This proves
\eqref{eq:population-linearization-final}.

For the unconditional claim, put
$\mathcal H_{n,P}=\mathcal G_{n,P}\cap\mathcal E_n$.
If
\[
q_{n,P}
=
P_P\!\left(
R_{n,P}^{\mathrm{pop}}
>
\varepsilon\ell_n^{-1/2}
\mid\mathcal F_{0n}
\right),
\]
then, for every $\delta>0$,
\[
E_P\{\mathbf 1_{\mathcal H_{n,P}}q_{n,P}\}
\le
\delta+
P_P\!\left(
\mathcal H_{n,P}\cap\{q_{n,P}>\delta\}
\right).
\]
The second term is $o(1)$ uniformly in $P$ by the conditional
statement. Letting $\delta\downarrow0$ and using
\[
\sup_P P_P(\mathcal H_{n,P}^c)
\le
\eta_n+\sup_P P_P(\mathcal G_{n,P}^c)
=
\eta_n+o(1)
\]
proves \eqref{eq:unconditional-population-bound-final}.
\end{proof}

\begin{remark}[Explicit envelope conditions and unresolved rates]
\label{rem:explicit-envelope-conditions}
On every smooth branch of the clipped logit,
\[
|\ell_{\epsilon_n}'(t)|
\le
\frac{1}{\epsilon_n(1-\epsilon_n)},
\qquad
|\ell_{\epsilon_n}''(t)|
\le
\frac{1}{
\epsilon_n^2(1-\epsilon_n)^2
}.
\]
Indeed, the derivatives vanish on the two constant clipping branches.
On $(\epsilon_n,1-\epsilon_n)$,
\[
\ell_{\epsilon_n}'(t)
=
\frac{1}{t(1-t)},
\qquad
\ell_{\epsilon_n}''(t)
=
\frac{2t-1}{t^2(1-t)^2}.
\]
The stated bounds follow from
$t(1-t)\ge\epsilon_n(1-\epsilon_n)$ and $|2t-1|\le1$.

The quantities
$\overline\Gamma_{\lambda,n}$ and
$\mathfrak D_{\lambda,n}$ are explicit deterministic envelopes once
$\mathcal F_{0n}$ is fixed. They are not yet primitive rates in
$n$, $\rho_n$, cell width, or wavelet scale. To derive an end-to-end
primitive result for the implemented estimator, the paper must still prove
\eqref{eq:explicit-curvature-rate} and
\eqref{eq:explicit-studentizer-rate} from lower bounds on cell
information, upper bounds on the first two derivatives of $B_n$,
the clipping margin, and the chosen resolution.

The per-dyad condition in \eqref{eq:explicit-studentizer-rate} is
sufficient but may be unnecessarily strong in a sparse model. It may
be replaced by the weaker directly weighted condition
\[
\max_{(F,r)\in\mathcal A_n}
\frac{
\left|
\sum_{k\in\mathcal K_n^F}
(b_{kr}^F)^2
\{\widehat v_k^F-p_k(1-p_k)\}
\right|
}{
V_r^{F,0}
}
=
o_{\mathbb P\mid\mathcal F_{0n}}(\ell_n^{-1}),
\]
provided this rate is proved uniformly on $\mathcal G_{n,P}$.
\end{remark}

\section{Shared-dyad multiplier calibration}\label{sec:supp-mult}


\paragraph{Uniform conditional stochastic order.}
Let $\mathcal P_n$ denote the stated model class. For any event
$\mathcal H_n\in\mathcal F_{0n}$ and any positive
$\mathcal F_{0n}$-measurable sequence $a_n$, write
\[
X_n
=
o_{\mathbb P\mid\mathcal F_{0n}}^{\mathrm{unif}}
(a_n;\mathcal H_n)
\]
if, for every $\varepsilon,\delta>0$,
\begin{equation}
\sup_{P\in\mathcal P_n}
P_P\!\left[
\mathcal H_n\cap
\left\{
P_P\!\left(
|X_n|>\varepsilon a_n
\,\middle|\,
\mathcal F_{0n}
\right)>\delta
\right\}
\right]
\longrightarrow0.
\label{eq:uniform-conditional-op}
\end{equation}
For an $\mathcal F_{0n}$-measurable variable $Z_n$, write
\[
Z_n=o_{\mathbb P}^{\mathrm{unif}}(a_n;\mathcal H_n)
\]
if, for every $\varepsilon>0$,
\begin{equation}
\sup_{P\in\mathcal P_n}
P_P\!\left(
\mathcal H_n\cap
\{|Z_n|>\varepsilon a_n\}
\right)
\longrightarrow0.
\label{eq:uniform-design-op}
\end{equation}
Regular conditional probabilities are assumed to exist. All suprema
over $t\in\mathbb R$ may equivalently be taken over $t\in\mathbb Q$.
Objects defined only on $\mathcal H_n$ are assigned arbitrary
measurable values on $\mathcal H_n^c$.

\begin{theorem}[Uniform shared-dyad multiplier calibration,\label{thm:fwer}
simultaneous coverage, and strong FWER]
\label{thm:uniform-shared-multiplier}
Fix $\alpha\in(0,1)$, assume $n\ge2$, and let
\[
\mathcal E_n=\{\Delta_n\le r_n\}\in\mathcal F_{0n}.
\]
Let $\mathcal I_n$ be the $\mathcal F_{0n}$-measurable admitted family,
and put
\[
m_n=|\mathcal I_n|,
\qquad
d_n=2m_n,
\qquad
\bar d_n=d_n\vee3,
\qquad
\ell_n=\log\bar d_n,
\qquad
K_n=\binom n2,
\]
and
\begin{equation}
\mathcal G_n
=
\mathcal E_n\cap\{m_n\ge1\}.
\label{eq:good-nonempty-event}
\end{equation}
On $\{m_n=0\}$, maxima over $\mathcal I_n$ and empty covariance index
sets are defined as zero, intersections over $\mathcal I_n$ equal the
whole sample space, and no hypothesis is rejected.

Let $\mathcal D_n^{\mathrm{act}}$ be the
$\mathcal F_{0n}$-measurable set of $B^E$- and $I$-fold dyads entering
the estimator. For every $P\in\mathcal P_n$, assume
\begin{equation}
\mathcal L_P\!\left(
(A_k)_{k\in\mathcal D_n^{\mathrm{act}}}
\,\middle|\,
\mathcal F_{0n}
\right)
=
\bigotimes_{k\in\mathcal D_n^{\mathrm{act}}}
\operatorname{Bernoulli}(p_k),
\label{eq:conditional-product-law}
\end{equation}
where
\[
p_k=P_P(A_k=1\mid\mathcal F_{0n}).
\]

Let $\widetilde\pi_k>0$ denote the frozen observation-propensity value
used by the estimator: $\widetilde\pi_k=\pi_k$ when known and
$\widetilde\pi_k=\widehat\pi_k$ when estimated on outcome-free data.
Define
\[
W_k^F=\frac{R_k^FM_k}{\widetilde\pi_k}.
\]
No $B^E$- or $I$-fold outcome may enter $\mathcal F_{0n}$, the
admission rule, $\widetilde\pi_k$, or any other frozen quantity.

For $F\in\{I,B^E\}$ and cell $r$, define
\[
b_{kr}^F
=
1\{q(k)=r\}\frac{W_k^F}{D_r^F},
\]
with $b_{kr}^F=0$ when dyad $k$ is not in fold $F$. Assume
\[
q(k),\ R_k^F,\ M_k,\ \widetilde\pi_k,\ W_k^F,\ D_r^F,\ p_k
\]
are $\mathcal F_{0n}$-measurable and every required denominator is
finite and strictly positive. Set
\[
Y_r^F
=
\sum_{k=1}^{K_n}b_{kr}^F A_k,
\qquad
\mu_r^F
=
E_P(Y_r^F\mid\mathcal F_{0n})
=
\sum_{k=1}^{K_n}b_{kr}^F p_k.
\]

Let
\[
\Phi_{\lambda n}(y^I,y^{B^E};\mathcal F_{0n})
\]
denote the complete frozen estimator map, including stabilization,
clipping, the frozen baseline smoother, projection onto $V_{J_0}$,
and both stochastic-fold paths. Thus
\[
\widehat\theta_\lambda
=
\Phi_{\lambda n}(Y^I,Y^{B^E};\mathcal F_{0n}).
\]
Define the frozen finite-design target by
\begin{equation}
\theta_{\lambda,n}^D
=
\Phi_{\lambda n}(\mu^I,\mu^{B^E};\mathcal F_{0n}).
\label{eq:finite-design-target}
\end{equation}
This frozen-operator target need not equal
$E_P(\widehat\theta_\lambda\mid\mathcal F_{0n})$.

Assume $\Phi_{\lambda n}$ is differentiable at
$(\mu^I,\mu^{B^E})$ and define
\begin{equation}
\overline G_{\lambda r}^F
=
\left.
\frac{
\partial\Phi_{\lambda n}
(y^I,y^{B^E};\mathcal F_{0n})
}{
\partial y_r^F
}
\right|_{(y^I,y^{B^E})=(\mu^I,\mu^{B^E})}.
\label{eq:oracle-derivative}
\end{equation}
The conditional mean is assumed to lie at neither a clipping nor a
stabilization knot.

For dyad $k$, define
\[
\gamma_{\lambda k}
=
\sum_{F\in\{I,B^E\}}\sum_r
\overline G_{\lambda r}^F b_{kr}^F,
\]
\begin{align}
s_\lambda^2
&=
\sum_{k=1}^{K_n}
\gamma_{\lambda k}^2p_k(1-p_k),
\label{eq:oracle-standard-error}\\
\zeta_{\lambda k}
&=
\frac{\gamma_{\lambda k}}{s_\lambda}(A_k-p_k).
\label{eq:oracle-influence-variable}
\end{align}
Set $\zeta_{\lambda k}=0$ for inactive dyads and assume
\[
0<s_\lambda<\infty
\qquad
\text{for every }\lambda\in\mathcal I_n.
\]
Then
\[
E_P(\zeta_{\lambda k}\mid\mathcal F_{0n})=0,
\qquad
\sum_{k=1}^{K_n}
E_P(\zeta_{\lambda k}^2\mid\mathcal F_{0n})=1,
\]
and
\begin{equation}
\Sigma_{\lambda\mu,n}
=
\frac{
\sum_{k=1}^{K_n}
\gamma_{\lambda k}\gamma_{\mu k}p_k(1-p_k)
}{
s_\lambda s_\mu
}.
\label{eq:oracle-correlation}
\end{equation}

Let $\widehat G_{\lambda r}^F$ be the fitted derivative used by the
implemented score, where $\widehat G^{B^E}$ denotes the paper's
baseline-fold derivative $G^B$, and define
\begin{equation}
\widehat V_r^F
=
\sum_{k=1}^{K_n}
(b_{kr}^F)^2
\widehat p_r^F(1-\widehat p_r^F),
\qquad
F\in\{I,B^E\}.
\label{eq:fitted-meat}
\end{equation}
Here $\widehat p_r^I=\widetilde p_r^I$ and
$\widehat p_r^{B^E}$ is the clipped, stabilized baseline-input
probability used in the feasible studentizer. Assume
\[
0\le\widehat V_r^F<\infty.
\]
Define
\begin{align}
\widehat s_\lambda^2
&=
\sum_{F\in\{I,B^E\}}\sum_r
(\widehat G_{\lambda r}^F)^2\widehat V_r^F,
\label{eq:feasible-standard-error}\\
\widehat\Sigma_{\lambda\mu,n}
&=
\frac{
\sum_{F\in\{I,B^E\}}\sum_r
\widehat G_{\lambda r}^F
\widehat G_{\mu r}^F
\widehat V_r^F
}{
\widehat s_\lambda\widehat s_\mu
}.
\label{eq:feasible-correlation}
\end{align}
Assume
\[
0<\widehat s_\lambda<\infty
\qquad
\text{for every admitted }\lambda.
\]
Then
\[
\widehat\Sigma_{\lambda\lambda,n}=1.
\]

Generate
\[
\xi_r^{F,*}\mid A,\mathcal F_{0n}
\sim N(0,\widehat V_r^F)
\]
independently over $(F,r)$, but reuse the same vector for every atom.
Define
\begin{equation}
Z_\lambda^*
=
\widehat s_\lambda^{-1}
\sum_{F\in\{I,B^E\}}\sum_r
\widehat G_{\lambda r}^F\xi_r^{F,*},
\qquad
T^*
=
\max_{\lambda\in\mathcal I_n}|Z_\lambda^*|.
\label{eq:shared-multiplier}
\end{equation}
Then
\begin{equation}
(Z_\lambda^*:\lambda\in\mathcal I_n)
\mid A,\mathcal F_{0n}
\sim N(0,\widehat\Sigma_n).
\label{eq:exact-bootstrap-Gaussian-law}
\end{equation}

Assume uniformly over $P\in\mathcal P_n$ on $\mathcal G_n$:

\begin{enumerate}
\item[(i)]
\begin{equation}
\mathfrak r_n^D
:=
\max_{\lambda\in\mathcal I_n}
\left|
\frac{\widehat\theta_\lambda-\theta_{\lambda,n}^D}{s_\lambda}
-
\sum_{k=1}^{K_n}\zeta_{\lambda k}
\right|
=
o_{\mathbb P\mid\mathcal F_{0n}}^{\mathrm{unif}}
(\ell_n^{-1/2};\mathcal G_n).
\label{eq:assumed-linearization}
\end{equation}

\item[(ii)] With
\[
L_n^{\mathrm{or}}
=
\max_{\lambda,k}\frac{|\gamma_{\lambda k}|}{s_\lambda},
\]
\begin{equation}
(L_n^{\mathrm{or}})^2\log^7(\bar d_nK_n)
=
o_{\mathbb P}^{\mathrm{unif}}(1;\mathcal G_n).
\label{eq:oracle-leverage-rate}
\end{equation}

\item[(iii)] With
\[
d_{s,n}
=
\max_{\lambda\in\mathcal I_n}
\left|
\frac{\widehat s_\lambda}{s_\lambda}-1
\right|,
\]
\begin{equation}
d_{s,n}\ell_n
=
o_{\mathbb P\mid\mathcal F_{0n}}^{\mathrm{unif}}
(1;\mathcal G_n).
\label{eq:studentizer-rate}
\end{equation}

\item[(iv)] Define
\[
\Delta_{\Sigma,n}
=
\max_{\lambda,\mu\in\mathcal I_n}
|\widehat\Sigma_{\lambda\mu,n}-\Sigma_{\lambda\mu,n}|.
\]
Since both matrices are correlation matrices,
\[
0\le\Delta_{\Sigma,n}\le2.
\]
With value zero when $\Delta_{\Sigma,n}=0$, let
\begin{equation}
\omega_{\Sigma,n}
=
\left[
\Delta_{\Sigma,n}
\left\{
1\vee
\log\!\left(
\frac{\bar d_n}{\Delta_{\Sigma,n}}
\right)
\right\}^{2}
\right]^{1/3}.
\label{eq:comparison-error}
\end{equation}
Assume
\begin{equation}
\omega_{\Sigma,n}
=
o_{\mathbb P\mid\mathcal F_{0n}}^{\mathrm{unif}}
(1;\mathcal G_n).
\label{eq:covariance-consistency}
\end{equation}
A stronger sufficient condition is
\begin{equation}
\Delta_{\Sigma,n}
=
o_{\mathbb P\mid\mathcal F_{0n}}^{\mathrm{unif}}
(\ell_n^{-2};\mathcal G_n).
\label{eq:simple-covariance-condition}
\end{equation}
\end{enumerate}

Define
\[
T_D^\circ
=
\max_{\lambda\in\mathcal I_n}
\frac{
|\widehat\theta_\lambda-\theta_{\lambda,n}^D|
}{
s_\lambda
},
\qquad
\widehat T_D^\circ
=
\max_{\lambda\in\mathcal I_n}
\frac{
|\widehat\theta_\lambda-\theta_{\lambda,n}^D|
}{
\widehat s_\lambda
}.
\]
Then
\begin{align}
\sup_t
\left|
P_P(T_D^\circ\le t\mid\mathcal F_{0n})
-
P^*(T^*\le t\mid A,\mathcal F_{0n})
\right|
&=
o_{\mathbb P\mid\mathcal F_{0n}}^{\mathrm{unif}}
(1;\mathcal E_n),
\label{eq:oracle-calibration}\\
\sup_t
\left|
P_P(\widehat T_D^\circ\le t\mid\mathcal F_{0n})
-
P^*(T^*\le t\mid A,\mathcal F_{0n})
\right|
&=
o_{\mathbb P\mid\mathcal F_{0n}}^{\mathrm{unif}}
(1;\mathcal E_n).
\label{eq:feasible-calibration}
\end{align}

For the Lebesgue coefficient, let $\mathfrak D_n$ be a deterministic
finite master dictionary such that
\[
\mathcal I_n\subseteq\mathfrak D_n
\qquad\text{almost surely},
\]
and, for every $\lambda\in\mathfrak D_n$, let $\theta_\lambda(P)$
be the deterministic intended population coefficient. Define
\begin{equation}
\beta_n
=
\max_{\lambda\in\mathcal I_n}
\frac{
|\theta_{\lambda,n}^D-\theta_\lambda(P)|
}{
s_\lambda
}.
\label{eq:standardized-target-bias}
\end{equation}
Assume
\begin{equation}
\beta_n\sqrt{\ell_n}
=
o_{\mathbb P}^{\mathrm{unif}}(1;\mathcal G_n).
\label{eq:finite-to-Lebesgue}
\end{equation}
Then the preceding calibration conclusions remain valid with
$\theta_{\lambda,n}^D$ replaced by $\theta_\lambda(P)$.

For either target
$\vartheta_{\lambda n}\in\{\theta_{\lambda,n}^D,\theta_\lambda(P)\}$,
where the latter requires \eqref{eq:finite-to-Lebesgue}, define
\[
F_n^*(t;A)
=
P^*(T^*\le t\mid A,\mathcal F_{0n}),
\qquad
c_{1-\alpha}^*(A)
=
\inf\{t:F_n^*(t;A)\ge1-\alpha\},
\]
and
\[
C_{\lambda n}
=
\left[
\widehat\theta_\lambda-c_{1-\alpha}^*\widehat s_\lambda,\,
\widehat\theta_\lambda+c_{1-\alpha}^*\widehat s_\lambda
\right].
\]
Then
\begin{equation}
\left[
1-\alpha
-
P_P\!\left(
\vartheta_{\lambda n}\in C_{\lambda n}
\text{ for all }\lambda\in\mathcal I_n
\,\middle|\,
\mathcal F_{0n}
\right)
\right]_+
=
o_{\mathbb P}^{\mathrm{unif}}(1;\mathcal E_n).
\label{eq:conditional-coverage}
\end{equation}

Under \eqref{eq:finite-to-Lebesgue}, define
\[
\mathcal H_{0,n}(P,\mathcal F_{0n})
=
\{\lambda\in\mathcal I_n:\theta_\lambda(P)=0\}.
\]
Reject $H_{0\lambda}$ when
\[
|\widehat\theta_\lambda|/\widehat s_\lambda>c_{1-\alpha}^*.
\]
Then
\begin{equation}
\left[
P_P\!\left(
\exists\lambda\in\mathcal H_{0,n}(P,\mathcal F_{0n}):
H_{0\lambda}\text{ is rejected}
\,\middle|\,
\mathcal F_{0n}
\right)
-\alpha
\right]_+
=
o_{\mathbb P}^{\mathrm{unif}}(1;\mathcal E_n).
\label{eq:conditional-strong-fwer}
\end{equation}

If an empirical quantile based on deterministic $B_n^*$ independent
multiplier draws is used, then, conditionally on $A,\mathcal F_{0n}$,
\begin{equation}
P_{\mathrm{MC}}^*\!\left(
\sup_t|\widehat F_{B_n^*}^*(t)-F_n^*(t;A)|>\tau
\,\middle|\,
A,\mathcal F_{0n}
\right)
\le2e^{-2B_n^*\tau^2}.
\label{eq:conditional-DKW}
\end{equation}
All conclusions remain valid under the joint data--Monte Carlo law if
$B_n^*\to\infty$.

Finally, if
\[
\sup_{P\in\mathcal P_n}P_P(\mathcal E_n^c)
\le\eta_n,
\qquad
\eta_n\to0,
\]
then
\begin{align}
\inf_{P\in\mathcal P_n}
P_P\!\left(
\theta_\lambda(P)\in C_{\lambda n}
\text{ for all }\lambda\in\mathcal I_n
\right)
&\ge1-\alpha-\eta_n-o(1),
\label{eq:unconditional-coverage}\\
\sup_{P\in\mathcal P_n}
P_P\!\left(
\exists\lambda\in\mathcal H_{0,n}(P,\mathcal F_{0n}):
H_{0\lambda}\text{ is rejected}
\right)
&\le\alpha+\eta_n+o(1).
\label{eq:unconditional-strong-fwer}
\end{align}
\end{theorem}

\begin{proof}
On $\{m_n=0\}$ the calibration distances are zero, simultaneous
coverage is vacuous, and the FWER is zero. Hence work on $\mathcal G_n$.

Set
\[
S_{\lambda n}=\sum_{k=1}^{K_n}\zeta_{\lambda k},
\qquad
S_n=\max_{\lambda\in\mathcal I_n}|S_{\lambda n}|.
\]
For $\sigma\in\{-1,1\}$, set
\[
\widetilde\zeta_{(\lambda,\sigma),k}
=
\sigma\zeta_{\lambda k}.
\]
Then
\begin{equation}
S_n
=
\max_{(\lambda,\sigma)\in
\mathcal I_n\times\{-1,1\}}
\sum_{k=1}^{K_n}
\widetilde\zeta_{(\lambda,\sigma),k}.
\label{eq:signed-representation}
\end{equation}

Conditional on $\mathcal F_{0n}$, the vectors
$\widetilde\zeta_k\in\mathbb R^{d_n}$ are independent and centered,
and every coordinate has total conditional variance one. Let
\[
X_k=\sqrt{K_n}\,\widetilde\zeta_k,
\qquad
B_n^{\mathrm{CCK}}
=
C_0\sqrt{K_n}L_n^{\mathrm{or}},
\qquad
C_0\ge(\log2)^{-1}.
\]
For $r=1,2$ and every coordinate $a$,
\begin{align}
\frac1{K_n}\sum_{k=1}^{K_n}
E_P(|X_{ka}|^{2+r}\mid\mathcal F_{0n})
&\le
(\sqrt{K_n}L_n^{\mathrm{or}})^r
\frac1{K_n}\sum_{k=1}^{K_n}
E_P(X_{ka}^2\mid\mathcal F_{0n})
\notag\\
&=
(\sqrt{K_n}L_n^{\mathrm{or}})^r
\le
(B_n^{\mathrm{CCK}})^r.
\label{eq:CCK-moment-bound}
\end{align}
Moreover,
\[
E_P\!\left[
\exp\!\left(
\frac{|X_{ka}|}{B_n^{\mathrm{CCK}}}
\right)
\,\middle|\,
\mathcal F_{0n}
\right]
\le e^{1/C_0}\le2.
\]
Since
\[
1
=
\sum_{k=1}^{K_n}
E_P(\widetilde\zeta_{a,k}^2\mid\mathcal F_{0n})
\le K_n(L_n^{\mathrm{or}})^2,
\]
we have $B_n^{\mathrm{CCK}}\ge1$.

Let
\[
G_n\mid\mathcal F_{0n}\sim N(0,\Sigma_n),
\qquad
M_n^G=\max_{\lambda\in\mathcal I_n}|G_{\lambda n}|.
\]
Proposition~2.1 of
\citet{ChernozhukovChetverikovKato2017} gives
\begin{align}
\delta_{\mathrm{CLT},n}
&:=
\sup_t
\left|
P_P(S_n\le t\mid\mathcal F_{0n})
-
P_P(M_n^G\le t\mid\mathcal F_{0n})
\right|
\notag\\
&\le
1\wedge
C\left\{
(L_n^{\mathrm{or}})^2\log^7(\bar d_nK_n)
\right\}^{1/6}
=
o_{\mathbb P}^{\mathrm{unif}}(1;\mathcal G_n).
\label{eq:conditional-CCK}
\end{align}

The equal-variance anti-concentration inequality of
\citet{ChernozhukovChetverikovKato2015} gives
\begin{equation}
\sup_t
P_P\!\left(
|M_n^G-t|\le\varepsilon
\,\middle|\,
\mathcal F_{0n}
\right)
\le C\varepsilon\sqrt{\ell_n}.
\label{eq:Gaussian-anticoncentration}
\end{equation}

Since $|T_D^\circ-S_n|\le\mathfrak r_n^D$, for every fixed $u>0$,
\begin{align}
&\sup_t
\left|
P_P(T_D^\circ\le t\mid\mathcal F_{0n})
-
P_P(M_n^G\le t\mid\mathcal F_{0n})
\right|
\notag\\
&\quad\le
2\delta_{\mathrm{CLT},n}
+
P_P\!\left(
\mathfrak r_n^D>\frac{u}{\sqrt{\ell_n}}
\,\middle|\,
\mathcal F_{0n}
\right)
+
Cu.
\label{eq:linearization-sandwich}
\end{align}
Letting first $n\to\infty$ and then $u\downarrow0$ proves
\begin{equation}
\sup_t
\left|
P_P(T_D^\circ\le t\mid\mathcal F_{0n})
-
P_P(M_n^G\le t\mid\mathcal F_{0n})
\right|
=
o_{\mathbb P}^{\mathrm{unif}}(1;\mathcal G_n).
\label{eq:oracle-Gaussian-approximation}
\end{equation}

For $a=(\lambda,\sigma)$ and $b=(\mu,\tau)$, define
\[
\Gamma_{ab,n}=\sigma\tau\Sigma_{\lambda\mu,n},
\qquad
\widehat\Gamma_{ab,n}
=
\sigma\tau\widehat\Sigma_{\lambda\mu,n}.
\]
Then
\[
\|\widehat\Gamma_n-\Gamma_n\|_{\max}
=
\Delta_{\Sigma,n}.
\]
The Gaussian comparison inequality of
\citet{ChernozhukovChetverikovKato2015} gives
\begin{align}
&\sup_t
\left|
P_P(M_n^G\le t\mid\mathcal F_{0n})
-
P^*(T^*\le t\mid A,\mathcal F_{0n})
\right|
\notag\\
&\quad\le
C\Delta_{\Sigma,n}^{1/3}
\left\{
1\vee
\log\!\left(
\frac{\bar d_n}{\Delta_{\Sigma,n}}
\right)
\right\}^{2/3}
=
C\omega_{\Sigma,n}.
\label{eq:Gaussian-comparison}
\end{align}
Together with \eqref{eq:oracle-Gaussian-approximation}, this proves
\eqref{eq:oracle-calibration}.

To verify \eqref{eq:simple-covariance-condition}, define
\[
U_{\Delta,n}
=
P_P(\Delta_{\Sigma,n}>1\mid\mathcal F_{0n}).
\]
Then
\[
U_{\Delta,n}
=
o_{\mathbb P}^{\mathrm{unif}}(1;\mathcal G_n).
\]
On $\{0<\Delta_{\Sigma,n}\le1\}$,
\begin{align}
&\Delta_{\Sigma,n}
\left\{
1\vee
\log\!\left(
\frac{\bar d_n}{\Delta_{\Sigma,n}}
\right)
\right\}^{2}
\notag\\
&\quad\le
\Delta_{\Sigma,n}
+
2\Delta_{\Sigma,n}\ell_n^2
+
2\Delta_{\Sigma,n}\log^2(1/\Delta_{\Sigma,n}).
\label{eq:simple-covariance-implication}
\end{align}
The first two terms vanish under
\eqref{eq:simple-covariance-condition}. The same condition implies
$\Delta_{\Sigma,n}\to0$ conditionally and uniformly on $\mathcal G_n$,
while $x\log^2(1/x)\to0$ as $x\downarrow0$. Hence the final term is
negligible, proving \eqref{eq:covariance-consistency}.

On $\{d_{s,n}\le1/2\}$,
\begin{equation}
|\widehat T_D^\circ-T_D^\circ|
\le
\frac{d_{s,n}}{1-d_{s,n}}T_D^\circ
\le2d_{s,n}T_D^\circ.
\label{eq:studentizer-difference}
\end{equation}
The Gaussian union bound and
\eqref{eq:oracle-Gaussian-approximation} imply conditional uniform
tightness of $T_D^\circ/\sqrt{\ell_n}$. Consequently, for every
$M,\varepsilon>0$,
\begin{align}
&P_P\!\left(
\sqrt{\ell_n}
|\widehat T_D^\circ-T_D^\circ|>\varepsilon
\,\middle|\,
\mathcal F_{0n}
\right)
\notag\\
&\quad\le
P_P(d_{s,n}>1/2\mid\mathcal F_{0n})
+
P_P\!\left(
d_{s,n}\ell_n>\frac{\varepsilon}{2M}
\,\middle|\,
\mathcal F_{0n}
\right)
\notag\\
&\qquad+
P_P\!\left(
\frac{T_D^\circ}{\sqrt{\ell_n}}>M
\,\middle|\,
\mathcal F_{0n}
\right).
\label{eq:studentizer-product-bound}
\end{align}
First choose $M$ large and then apply
\eqref{eq:studentizer-rate}. Thus
\[
|\widehat T_D^\circ-T_D^\circ|
=
o_{\mathbb P\mid\mathcal F_{0n}}^{\mathrm{unif}}
(\ell_n^{-1/2};\mathcal G_n),
\]
and another anti-concentration sandwich proves
\eqref{eq:feasible-calibration}.

For the Lebesgue target,
\[
\max_{\lambda\in\mathcal I_n}
\left|
\frac{\widehat\theta_\lambda-\theta_\lambda(P)}{s_\lambda}
-
\sum_k\zeta_{\lambda k}
\right|
\le
\mathfrak r_n^D+\beta_n.
\]
Thus \eqref{eq:finite-to-Lebesgue} gives the same uniform expansion
around $\theta_\lambda(P)$.

For the chosen target, let
\[
\widehat T_n^\vartheta
=
\max_{\lambda\in\mathcal I_n}
\frac{
|\widehat\theta_\lambda-\vartheta_{\lambda n}|
}{
\widehat s_\lambda
},
\]
\[
F_n^\vartheta(t)
=
P_P(\widehat T_n^\vartheta\le t\mid\mathcal F_{0n}),
\qquad
D_n^\vartheta(A)
=
\sup_t|F_n^\vartheta(t)-F_n^*(t;A)|.
\]
Then
\[
D_n^\vartheta(A)
=
o_{\mathbb P\mid\mathcal F_{0n}}^{\mathrm{unif}}
(1;\mathcal E_n).
\]
For $0<\delta<1-\alpha$, on
$\{D_n^\vartheta(A)\le\delta\}$,
\[
c_{1-\alpha}^*(A)
\ge
(F_n^\vartheta)^{\leftarrow}(1-\alpha-\delta).
\]
Therefore,
\begin{align}
P_P\!\left(
\widehat T_n^\vartheta>c_{1-\alpha}^*
\,\middle|\,
\mathcal F_{0n}
\right)
&\le
\alpha+\delta
+
P_P(D_n^\vartheta(A)>\delta\mid\mathcal F_{0n}).
\label{eq:quantile-transfer}
\end{align}
Letting $n\to\infty$ and then $\delta\downarrow0$ proves
\eqref{eq:conditional-coverage}.

For the empirical critical value, the
Dvoretzky--Kiefer--Wolfowitz--Massart inequality
\citep{Massart1990} gives \eqref{eq:conditional-DKW}. If
$\delta+\tau<1-\alpha$, the same generalized-inverse argument yields
\begin{align}
&P_{P,\mathrm{MC}}\!\left(
\widehat T_n^\vartheta>
\widehat c_{1-\alpha,B_n^*}^*
\,\middle|\,
\mathcal F_{0n}
\right)
\notag\\
&\quad\le
\alpha+\delta+\tau
+
P_P(D_n^\vartheta(A)>\delta\mid\mathcal F_{0n})
+
2e^{-2B_n^*\tau^2}.
\label{eq:finite-MC-bound}
\end{align}
Taking
\[
\tau_n
=
\sqrt{\frac{\log(2B_n^*)}{2B_n^*}}
\]
proves validity when $B_n^*\to\infty$.

For every true null, $\theta_\lambda(P)=0$, and hence
\begin{align}
&\left\{
\exists\lambda\in\mathcal H_{0,n}(P,\mathcal F_{0n}):
H_{0\lambda}\text{ is rejected}
\right\}
\notag\\
&\quad\subseteq
\left\{
\max_{\lambda\in\mathcal I_n}
\frac{
|\widehat\theta_\lambda-\theta_\lambda(P)|
}{
\widehat s_\lambda
}
>
c_{1-\alpha}^*
\right\}.
\label{eq:FWER-inclusion}
\end{align}
This proves conditional strong FWER under every configuration of false
nulls, without subset pivotality.

Let $e_{n,P}$ denote either conditional excess error in
\eqref{eq:conditional-coverage} or
\eqref{eq:conditional-strong-fwer}. Since $0\le e_{n,P}\le1$, for
every $\varepsilon>0$,
\begin{align}
\sup_{P\in\mathcal P_n}
E_P[1_{\mathcal E_n}e_{n,P}]
&\le
\varepsilon
+
\sup_{P\in\mathcal P_n}
P_P\!\left(
\mathcal E_n\cap\{e_{n,P}>\varepsilon\}
\right).
\label{eq:bounded-error-integration}
\end{align}
Consequently,
\[
\limsup_{n\to\infty}
\sup_{P\in\mathcal P_n}
E_P[1_{\mathcal E_n}e_{n,P}]
\le\varepsilon.
\]
Letting $\varepsilon\downarrow0$ gives
\begin{equation}
\sup_{P\in\mathcal P_n}
E_P[1_{\mathcal E_n}e_{n,P}]
=o(1).
\label{eq:integrated-excess-error}
\end{equation}

For coverage failure,
\begin{align*}
&P_P\!\left(
\exists\lambda\in\mathcal I_n:
\theta_\lambda(P)\notin C_{\lambda n}
\right)
\\
&\quad\le
P_P(\mathcal E_n^c)
+
E_P\!\left[
1_{\mathcal E_n}
P_P\!\left(
\exists\lambda\in\mathcal I_n:
\theta_\lambda(P)\notin C_{\lambda n}
\,\middle|\,
\mathcal F_{0n}
\right)
\right]
\\
&\quad\le
\eta_n+\alpha+o(1).
\end{align*}
Likewise,
\begin{align*}
&P_P\!\left(
\exists\lambda\in\mathcal H_{0,n}(P,\mathcal F_{0n}):
H_{0\lambda}\text{ is rejected}
\right)
\\
&\quad\le
P_P(\mathcal E_n^c)
+
E_P\!\left[
1_{\mathcal E_n}
P_P\!\left(
\exists\lambda\in\mathcal H_{0,n}(P,\mathcal F_{0n}):
H_{0\lambda}\text{ is rejected}
\,\middle|\,
\mathcal F_{0n}
\right)
\right]
\\
&\quad\le
\eta_n+\alpha+o(1).
\end{align*}
These bounds are uniform in $P$ and prove
\eqref{eq:unconditional-coverage} and
\eqref{eq:unconditional-strong-fwer}.

Finally,
\[
\operatorname{Cov}^*(Z_\lambda^*,Z_\mu^*
\mid A,\mathcal F_{0n})
=
\widehat\Sigma_{\lambda\mu,n}.
\]
Independent atomwise multipliers would set the off-diagonal bootstrap
covariances to zero. Moreover,
\[
\max_{\lambda\in\mathcal I_n}|Z_\lambda^*|
=
\max_{(\lambda,\sigma)\in
\mathcal I_n\times\{-1,1\}}
\sigma Z_\lambda^*,
\]
so the absolute maximum already calibrates both signs.
\end{proof}

\paragraph{Remaining estimator-specific proof obligations.}
This theorem is a complete high-level implication. To derive it from
Sections~2--3, separate lemmas must establish:
\begin{enumerate}
\item the conditional product law after conditioning on the full
missingness mask;
\item that all active-fold weights, gates, propensity estimates, and
tuning objects are frozen and outcome-free;
\item the uniform clipped/stabilized Taylor expansion
\eqref{eq:assumed-linearization};
\item pilot-to-oracle information and leverage transfer;
\item the studentizer and full max-norm covariance rates; and
\item the complete finite-design-to-Lebesgue bias rate, including
random-design and ordering discrepancies.
\end{enumerate}
Without items 1 or 3, the independent-sum Gaussian approximation is not
established. Without item 5, feasible simultaneous inference is not
established. Without item 6, the theorem remains valid only for
$\theta_{\lambda,n}^D$.

\section{Detection and localization boundary}\label{sec:supp-detect}


\begin{theorem}[Global adaptive detection and localization boundary]
\label{thm:adaptive-detection-localization}
Retain the sparse logistic fixed-design model and the notation
\(\alpha_n,\rho_n,B,\mathcal F_{0n},\Delta_n,r_n\) from the preceding
sections, and write
\[
\sigma(x)=\frac{e^x}{1+e^x}.
\]

Let $\mathscr P_n$ denote the statistical model class. Let
\[
\Lambda_n
=
\bigcup_{h\in\mathcal J_n}\Lambda_n(h),
\qquad
m_n=|\Lambda_n|\longrightarrow\infty,
\qquad
\ell_n=\log(2m_n),
\]
be a deterministic, outcome-independent packing of unique symmetric
off-diagonal Haar atoms, with transposed duplicates removed.  We take
\(\mathcal J_n\) to contain only scales \(h\) for which
\(\Lambda_n(h)\ne\varnothing\).  For
\(\lambda\in\Lambda_n\), let \(h_\lambda\) be its scale and put
\[
S_\lambda=\operatorname{supp}(\Psi_\lambda^S).
\]
Here and below, \(|S|\) denotes Lebesgue measure.

Assume
\[
\Psi_\lambda^S\in V_{J_0}^{\perp},
\qquad
\lambda\in\Lambda_n,
\]
and
\begin{equation}
\left\langle
\Psi_\lambda^S,\Psi_\mu^S
\right\rangle_{L^2([0,1]^2)}
=
\mathbf 1\{\lambda=\mu\},
\qquad
\lambda,\mu\in\Lambda_n.
\label{eq:t3-1}
\end{equation}
Assume also that atoms at a common scale have equal support area and,
for some \(\delta_{\rm pack}>0\),
\begin{equation}
\inf_n
\inf_{\substack{\lambda,\mu\in\Lambda_n\\
                 \lambda\ne\mu,\ h_\lambda=h_\mu}}
\frac{|S_\lambda\triangle S_\mu|}{2|S_\lambda|}
\ge\delta_{\rm pack}.
\label{eq:t3-2}
\end{equation}
The infimum over an empty collection is \(+\infty\).

Let
\[
\mathcal B_n
\subset
\left\{
f\in V_{J_0}:
f\text{ is symmetric and mean zero},\
\|f\|_\infty\le B
\right\}
\]
be deterministic, and define
\begin{equation}
\theta_\lambda(f)
=
\left\langle
(I-\Pi_{V_{J_0}})f,\Psi_\lambda^S
\right\rangle_{L^2([0,1]^2)}.
\label{eq:t3-3}
\end{equation}
For \(a>0\) and \(h\in\mathcal J_n\), put
\[
\mathfrak A_n(a,h)
=
\left\{
f_0+\varsigma a h\Psi_\lambda^S:
\begin{array}{l}
f_0\in\mathcal B_n,\quad
\varsigma\in\{-1,1\},\\
\lambda\in\Lambda_n(h),\quad
\|f_0+\varsigma a h\Psi_\lambda^S\|_\infty\le B
\end{array}
\right\}.
\]
By \eqref{eq:t3-1} and the orthogonal decomposition
\(V_{J_0}\oplus V_{J_0}^{\perp}\), every
\(f\in\mathfrak A_n(a,h)\) has a unique representation of this form.
Denote its atom and sign by \(\lambda(f)\) and \(\varsigma(f)\).

Define
\[
\mathfrak F_n^{\rm up}
=
\mathcal B_n
\cup
\bigcup_{h\in\mathcal J_n}
\bigcup_{\substack{a>0:\\
                   \mathfrak A_n(a,h)\ne\varnothing}}
\mathfrak A_n(a,h).
\]

Let \(\mathfrak S_{\rm sym}\) be the measure algebra of symmetric
Lebesgue-measurable subsets of \([0,1]^2\), identified modulo
Lebesgue-null sets and equipped with
\[
d_\triangle(S,T)=|S\triangle T|
\]
and its Borel sigma-field.  This is a standard Borel space.  Define
\[
\mathcal D_n^{\rm dec}
=
\{\partial\}
\sqcup
\bigl(\mathfrak S_{\rm sym}\times\mathcal J_n\bigr),
\]
where \(\partial\) is a cemetery symbol denoting no selected atom and
\(\mathcal J_n\) has the discrete sigma-field.  For
\(D=(S_D,h_D)\ne\partial\), define
\[
d_{\rm loc}(D,\lambda)
=
\min\left\{
1,\
\frac{|S_D\triangle S_\lambda|}{2|S_\lambda|}
+
\mathbf 1\{h_D\ne h_\lambda\}
\right\},
\]
and set
\[
d_{\rm loc}(\partial,\lambda)=1.
\]

Let \(\widehat{\mathcal I}_n\) be the pilot-admitted family defined
by the pilot gate, and put
\[
\widehat\Lambda_n
=
\Lambda_n\cap\widehat{\mathcal I}_n.
\]
Assume that \(\widehat{\mathcal I}_n\), and hence
\(\{\widehat\Lambda_n=\Lambda_n\}\), is
\(\mathcal F_{0n}\)-measurable and that
\begin{equation}
\gamma_n
:=
\sup_{f\in\mathfrak F_n^{\rm up}}
P_f\{\widehat\Lambda_n\ne\Lambda_n\}
=o(1).
\label{eq:t3-4}
\end{equation}
Assume also that
\[
\widehat\Lambda_n=\Lambda_n
\quad\Longrightarrow\quad
\min_{\lambda\in\Lambda_n}\widehat s_\lambda>0
\quad\text{almost surely}.
\]
If standard-error positivity is instead imposed as a separate gate,
that gate must be included in the admission event in
\eqref{eq:t3-4}.

Suppose that there is an \(\mathcal F_{0n}\)-measurable event
\(\mathcal G_n\subseteq\{\Delta_n\le r_n\}\) such that
\begin{equation}
\sup_{f\in\mathfrak F_n^{\rm up}}
P_f(\mathcal G_n^c)
\le\eta_n+o(1),
\qquad
\eta_n\to0.
\label{eq:t3-5}
\end{equation}
Define
\[
\widetilde{\mathcal G}_n
=
\mathcal G_n
\cap
\{\widehat\Lambda_n=\Lambda_n\}.
\]
Then
\begin{equation}
\sup_{f\in\mathfrak F_n^{\rm up}}
P_f(\widetilde{\mathcal G}_n^c)
\le
\eta_n+\gamma_n+o(1).
\label{eq:t3-6}
\end{equation}

For a random array \(R_{n,f}\), an
\(\mathcal F_{0n}\)-measurable event \(H_n\), and a deterministic
sequence \(q_n>0\), write
\[
R_{n,f}=o_{\rm ucp}(q_n)
\quad
\text{on \(H_n\), uniformly over \(\mathfrak F_n^{\rm up}\),}
\]
if, for every \(\epsilon,\delta>0\),
\[
\sup_{f\in\mathfrak F_n^{\rm up}}
P_f\left[
H_n
\cap
\left\{
P_f\left(
|R_{n,f}|>\epsilon q_n
\mid\mathcal F_{0n}
\right)>\delta
\right\}
\right]
\longrightarrow0.
\]

For \(f\in\mathfrak F_n^{\rm up}\), let
\(\theta_{\lambda,n}^{D}(f)\) be the frozen finite-design target and
put
\[
b_{\lambda,n}(f)
=
\theta_{\lambda,n}^{D}(f)-\theta_\lambda(f).
\]
Let \(s_{\lambda,f}^2\) be the conditional variance of the frozen
first-order influence sum and put
\[
I_{\lambda,f}=s_{\lambda,f}^{-2}.
\]
Assume that
\(\theta_{\lambda,n}^{D}(f)\), \(b_{\lambda,n}(f)\),
\(s_{\lambda,f}\), and \(I_{\lambda,f}\) are
\(\mathcal F_{0n}\)-measurable.  On
\(\widetilde{\mathcal G}_n\), assume uniformly over
\(f\in\mathfrak F_n^{\rm up}\) and
\(\lambda\in\Lambda_n\) that
\begin{equation}
0<\underline\iota n^2\rho_n
\le
I_{\lambda,f}
\le
\overline\iota n^2\rho_n
<\infty
\label{eq:t3-7}
\end{equation}
for constants
\(0<\underline\iota\le\overline\iota<\infty\).

Write
\[
p_{ij}(f)
=
E_f(A_{ij}\mid\mathcal F_{0n}).
\]
For each \(f\in\mathfrak F_n^{\rm up}\), assume that there are
\(\mathcal F_{0n}\)-measurable loadings
\(\gamma_{\lambda,ij}(f)\) such that
\[
\zeta_{\lambda,ij}(f)
=
s_{\lambda,f}^{-1}
\gamma_{\lambda,ij}(f)
\{A_{ij}-p_{ij}(f)\}.
\]
For every fixed \(\lambda\), conditional on \(\mathcal F_{0n}\),
assume that
\(\{\zeta_{\lambda,ij}(f):i<j\}\) are independent and, uniformly over
\(f\in\mathfrak F_n^{\rm up}\), satisfy
\begin{equation}
E_f\{\zeta_{\lambda,ij}(f)\mid\mathcal F_{0n}\}=0,
\qquad
\sum_{i<j}
E_f\{\zeta_{\lambda,ij}^2(f)\mid\mathcal F_{0n}\}=1.
\label{eq:t3-8}
\end{equation}
Dependence across different \(\lambda\)'s is unrestricted.  Assume,
for a deterministic sequence \(L_n\), that, for every
\(f\in\mathfrak F_n^{\rm up}\),
\begin{equation}
\mathbf 1_{\widetilde{\mathcal G}_n}
\max_{\lambda,i<j}
|\zeta_{\lambda,ij}(f)|
\le L_n
\quad P_f\text{-almost surely},
\qquad
L_n^2\log^7(m_n n)\to0.
\label{eq:t3-9}
\end{equation}

Define
\[
R_{n,f}^{\rm lin}
=
\max_{\lambda\in\Lambda_n}
\left|
\frac{
\widehat\theta_\lambda-\theta_{\lambda,n}^{D}(f)
}{
s_{\lambda,f}
}
-
\sum_{i<j}\zeta_{\lambda,ij}(f)
\right|.
\]
Assume
\begin{equation}
R_{n,f}^{\rm lin}
=
o_{\rm ucp}(\ell_n^{-1/2})
\quad
\text{on }\widetilde{\mathcal G}_n.
\label{eq:t3-10}
\end{equation}
Also define
\[
\delta_{n,f}
=
\max_{\lambda\in\Lambda_n}
\left|
\frac{\widehat s_\lambda}{s_{\lambda,f}}-1
\right|,
\]
and assume
\begin{equation}
\delta_{n,f}
=
o_{\rm ucp}(\ell_n^{-1})
\quad
\text{on }\widetilde{\mathcal G}_n.
\label{eq:t3-11}
\end{equation}

Assume, for every \(\epsilon>0\),
\begin{equation}
\sup_{f\in\mathfrak F_n^{\rm up}}
P_f\left[
\widetilde{\mathcal G}_n
\cap
\left\{
\sqrt{\ell_n}
\max_{\lambda\in\Lambda_n}
\frac{|b_{\lambda,n}(f)|}{s_{\lambda,f}}
>\epsilon
\right\}
\right]
\longrightarrow0.
\label{eq:t3-12}
\end{equation}
Here \(b_{\lambda,n}(f)\) contains only frozen,
\(\mathcal F_{0n}\)-measurable target discrepancies.  Stochastic
Taylor and nonlinear-estimation remainders are included in
\eqref{eq:t3-10}, not in \(b_{\lambda,n}(f)\).

A sufficient primitive verification of \eqref{eq:t3-12} is as
follows.  Decompose
\[
b_{\lambda,n}(f)
=
b_{\lambda,n}^{\rm nonord}(f)
+
b_{\lambda,n}^{\rm order}(f).
\]
Assume that the aggregate nonordering term satisfies the analogue of
\eqref{eq:t3-12} and that, for every
\(f\in\mathfrak F_n^{\rm up}\), on
\(\widetilde{\mathcal G}_n\),
\[
\max_{\lambda\in\Lambda_n}
|b_{\lambda,n}^{\rm order}(f)|
\le
C_{\rm ord}Br_n.
\]
Then
\begin{equation}
Br_n\,n\sqrt{\rho_n\ell_n}
\longrightarrow0
\label{eq:t3-13}
\end{equation}
is sufficient for the ordering term to satisfy the analogue of
\eqref{eq:t3-12}.  Indeed, by \eqref{eq:t3-7},
\[
\sqrt{\ell_n}
\max_\lambda
\frac{|b_{\lambda,n}^{\rm order}(f)|}{s_{\lambda,f}}
\le
C_{\rm ord}\sqrt{\overline\iota}\,
Br_n n\sqrt{\rho_n\ell_n}.
\]
For any sequences \(a_n,h_n\) satisfying the upper signal condition
\eqref{eq:t3-17} below,
\[
\frac{Br_n}{a_n h_n}
\le
\frac{
Br_n n\sqrt{\rho_n\ell_n}
}{
\sqrt{C_{\rm up}}\,\ell_n
}
\longrightarrow0.
\]
Thus \(r_n/h_n\le c_r\) alone is insufficient when \(a_n\to0\).

Let \(\mathcal X_n\) denote the full observed data.  Conditionally on
\(\mathcal X_n\), assume that
\[
Z^{*(b)}
=
\bigl(
Z_\lambda^{*(b)}:
\lambda\in\Lambda_n
\bigr),
\qquad b=1,\ldots,B_n^*,
\]
are independent centered jointly Gaussian vectors satisfying
\begin{equation}
\operatorname{Var}^*(Z_\lambda^{*(b)})=1,
\qquad
\lambda\in\Lambda_n.
\label{eq:t3-14}
\end{equation}
Dependence among coordinates within the same multiplier vector is
unrestricted.  Let \(P_{\mathcal X_n}^*\) be the conditional law of
one multiplier vector given \(\mathcal X_n\), and let
\(P_{\mathcal X_n}^{*,B}\) be the conditional product law of all
\(B_n^*\) multiplier vectors.  Assume \(B_n^*\to\infty\), and define
\[
T^{*(b)}
=
\max_{\lambda\in\Lambda_n}|Z_\lambda^{*(b)}|,
\qquad
F_n^*(t)
=
P_{\mathcal X_n}^*(T^{*(1)}\le t),
\]
\[
q_{1-u,n}^*
=
\inf\{t:F_n^*(t)\ge1-u\},
\]
and
\[
\widehat F_{B_n^*}^*(t)
=
\frac1{B_n^*}
\sum_{b=1}^{B_n^*}
\mathbf 1\{T^{*(b)}\le t\},
\]
\[
c_{1-\alpha}^*
=
\inf\left\{
t:
\widehat F_{B_n^*}^*(t)\ge1-\alpha
\right\},
\qquad
\alpha\in(0,1).
\]
Let \(P_f^\dagger\) and \(E_f^\dagger\) denote the joint law of the
data and the independent Gaussian seeds used for the multiplier
draws.

On \(\{\widehat\Lambda_n=\Lambda_n\}\), define
\[
Z_\lambda
=
\frac{\widehat\theta_\lambda}{\widehat s_\lambda},
\qquad
T_n
=
\max_{\lambda\in\Lambda_n}|Z_\lambda|.
\]
When admission fails, set
\[
Z_\lambda=0
\quad(\lambda\in\Lambda_n),
\qquad
T_n=0.
\]
Define
\[
\phi_n
=
\mathbf 1\{\widehat\Lambda_n=\Lambda_n\}
\mathbf 1\{T_n>c_{1-\alpha}^*\}.
\]
With deterministic tie-breaking, put
\[
\widehat\lambda
=
\begin{cases}
\displaystyle
\arg\max_{\lambda\in\Lambda_n}|Z_\lambda|,
&\phi_n=1,\\[2mm]
\varnothing,&\phi_n=0,
\end{cases}
\]
and
\[
(\widehat\varepsilon,\widehat D)
=
\begin{cases}
\left(
\operatorname{sign}(\widehat\theta_{\widehat\lambda}),
(S_{\widehat\lambda},h_{\widehat\lambda})
\right),
&\widehat\lambda\ne\varnothing,\\[1mm]
(0,\partial),
&\widehat\lambda=\varnothing.
\end{cases}
\]

For \(f_0\in\mathcal B_n\), define
\[
F_{n,f_0}^{0}(t)
=
P_{f_0}(T_n\le t\mid\mathcal F_{0n}),
\qquad
\Delta_{n,f_0}^{\rm cal}
=
\sup_{t\in\mathbb R}
|F_{n,f_0}^{0}(t)-F_n^*(t)|.
\]
Assume the following uniform calibration condition, as established
by the preceding multiplier theorem: there is a deterministic
\(\kappa_n\downarrow0\) such that
\begin{equation}
\sup_{f_0\in\mathcal B_n}
P_{f_0}\left[
\widetilde{\mathcal G}_n
\cap
\{\Delta_{n,f_0}^{\rm cal}>\kappa_n\}
\right]
\longrightarrow0.
\label{eq:t3-15}
\end{equation}
The calibration condition is required for the actual statistic
\(T_n\), the deterministic family \(\Lambda_n\), and the same
loadings and studentizers as the implemented scan.  Then
\begin{equation}
\sup_{f_0\in\mathcal B_n}
P_{f_0}^\dagger(\phi_n=1)
\le
\alpha+\eta_n+o(1).
\label{eq:t3-16}
\end{equation}
This is a global-null size statement.  Strong familywise error
control requires the corresponding subset-null calibration result.

There exists \(C_{\rm up}<\infty\) such that the following holds.
For every deterministic sequence
\[
h_n\in\mathcal J_n,
\qquad
a_n>0,
\qquad
\mathfrak A_n(a_n,h_n)\ne\varnothing,
\]
satisfying
\begin{equation}
a_n^2n^2\rho_nh_n^2
\ge
C_{\rm up}\ell_n,
\label{eq:t3-17}
\end{equation}
define
\[
\mathcal E_{n,f}^{\rm rec}
=
\left\{
\phi_n=0
\ \text{or}\
\widehat\lambda\ne\lambda(f)
\ \text{or}\
\widehat\varepsilon\ne\varsigma(f)
\right\}.
\]
Then, for every \(\delta>0\),
\begin{equation}
\sup_{f\in\mathfrak A_n(a_n,h_n)}
P_f\left[
\widetilde{\mathcal G}_n
\cap
\left\{
P_f^\dagger
\left(
\mathcal E_{n,f}^{\rm rec}
\mid\mathcal F_{0n}
\right)
>\delta
\right\}
\right]
\longrightarrow0.
\label{eq:t3-18}
\end{equation}
Consequently,
\begin{equation}
\sup_{f\in\mathfrak A_n(a_n,h_n)}
E_f^\dagger
d_{\rm loc}\{\widehat D,\lambda(f)\}
\le
\eta_n+\gamma_n+o(1)
\longrightarrow0.
\label{eq:t3-19}
\end{equation}

For the lower bound, all probabilities, expectations, and risks below
are conditional on a deterministic design
\((U_1,\ldots,U_n)\) satisfying the following conditions.  Suppose
that \(\mathcal B_n\) contains a symmetric, mean-zero interior
element \(f_0^\circ\) and constants \(\delta_B,c_p,C_p>0\) such that
\begin{equation}
\|f_0^\circ\|_\infty\le B-\delta_B,
\qquad
c_p\rho_n\le p_{0,ij}\le C_p\rho_n,
\qquad
1-p_{0,ij}\ge c_p,
\label{eq:t3-20}
\end{equation}
where
\[
p_{0,ij}
=
\sigma\{\alpha_n+f_0^\circ(U_i,U_j)\}.
\]

Let
\[
\mathcal K_n
=
\{(i,j):1\le i<j\le n\},
\qquad
v_{\lambda,ij}
=
h_\lambda\Psi_\lambda^S(U_i,U_j).
\]
Assume
\begin{equation}
\|v_\lambda\|_\infty\le C_\Psi,
\qquad
c_Xn^2h_\lambda^2
\le
\sum_{(i,j)\in\mathcal K_n}v_{\lambda,ij}^2
\le
C_Xn^2h_\lambda^2
\label{eq:t3-21}
\end{equation}
uniformly over \(\lambda\in\Lambda_n\).  The two-sided energy
bound records membership in the same regular-design class; only its
upper bound is used below.

Let
\[
\mathbf a_n=(a_h:h\in\mathcal J_n),
\qquad
a_h>0,
\]
and assume
\begin{equation}
\max_{h\in\mathcal J_n}
a_h
\max_{\lambda\in\Lambda_n(h)}
\|h\Psi_\lambda^S\|_\infty
\le\delta_B.
\label{eq:t3-22}
\end{equation}
For
\[
\mathcal D_\lambda
=
\{(i,j)\in\mathcal K_n:v_{\lambda,ij}\ne0\},
\]
define
\[
N_{{\rm ov},n}
=
\#\left\{
(\lambda,\mu)\in\Lambda_n^2:
\mathcal D_\lambda\cap\mathcal D_\mu\ne\varnothing
\right\},
\]
including every label-diagonal pair \((\lambda,\lambda)\).  Assume
\begin{equation}
N_{{\rm ov},n}
\le m_nd_n,
\qquad
d_n\ge1,
\qquad
\frac{\log d_n}{\log(2m_n)}
\longrightarrow0.
\label{eq:t3-23}
\end{equation}

Let \(\mathcal E_n^{\rm orc}\) be the experiment observing the
deterministic design, \(\alpha_n\), the fixed baseline
\(f_0^\circ\), and all outcomes
\(\{A_{ij}:(i,j)\in\mathcal K_n\}\).  Assume that, after restricting
the baseline to \(f_0^\circ\), the original experiment is obtained
from \(\mathcal E_n^{\rm orc}\) through a Markov kernel that does not
depend on \(\lambda\), the sign, or the amplitude.  This includes
outcome-independent missingness conditional on the fixed design.

Define
\[
\Theta_{1,n}(\mathbf a_n)
=
\left\{
\vartheta=(f_0,\varsigma,h,\lambda):
\begin{array}{l}
f_0\in\mathcal B_n,\quad
\varsigma\in\{-1,1\},\quad
h\in\mathcal J_n,\\
\lambda\in\Lambda_n(h),\quad
\|f_0+\varsigma a_h h\Psi_\lambda^S\|_\infty\le B
\end{array}
\right\},
\]
and let
\[
f_\vartheta
=
f_0+\varsigma a_h h\Psi_\lambda^S,
\qquad
\lambda(\vartheta)=\lambda.
\]
Let \(P_\vartheta\) and \(E_\vartheta\) denote the conditional law and
expectation in the original experiment.  Define
\begin{align}
R_n^{\rm det}
&=
\inf_\varphi
\left[
\sup_{f_0\in\mathcal B_n}P_{f_0}(\varphi=1)
+
\sup_{\vartheta\in\Theta_{1,n}(\mathbf a_n)}
P_\vartheta(\varphi=0)
\right],
\label{eq:t3-24}\\
R_n^{\rm loc}
&=
\inf_{\widetilde D}
\sup_{\vartheta\in\Theta_{1,n}(\mathbf a_n)}
E_\vartheta
d_{\rm loc}\{\widetilde D,\lambda(\vartheta)\},
\label{eq:t3-25}
\end{align}
where randomized tests are allowed and the second infimum is over
all possibly randomized measurable decisions taking values in
\(\mathcal D_n^{\rm dec}\).

There exists \(c_{\rm low}>0\) such that, if
\begin{equation}
\max_{h\in\mathcal J_n}
a_h^2n^2\rho_nh^2
\le
c_{\rm low}\ell_n,
\label{eq:t3-26}
\end{equation}
then
\begin{equation}
R_n^{\rm det}\longrightarrow1,
\qquad
\liminf_{n\to\infty}R_n^{\rm loc}>0.
\label{eq:t3-27}
\end{equation}
These are global-union minimax conclusions.
\end{theorem}

\begin{proof}
We first prove the upper assertions.  Fix
\[
f
=
f_0+\varsigma_0a_n h_n\Psi_{\lambda_0}^S
\in\mathfrak A_n(a_n,h_n),
\]
where
\[
\lambda_0=\lambda(f),
\qquad
\varsigma_0=\varsigma(f).
\]
By \eqref{eq:t3-1}--\eqref{eq:t3-3},
\begin{equation}
\theta_{\lambda_0}(f)
=
\varsigma_0a_n h_n,
\qquad
\theta_\mu(f)=0
\quad(\mu\ne\lambda_0).
\label{eq:t3-28}
\end{equation}

For notational brevity, write
\[
G_\lambda
=
\sum_{i<j}\zeta_{\lambda,ij}(f)
\]
and put
\[
r_{n,f}
=
\max_{\lambda\in\Lambda_n}
\left|
\frac{\widehat\theta_\lambda-\theta_\lambda(f)}
{s_{\lambda,f}}
-
G_\lambda
\right|.
\]
By the definition of \(b_{\lambda,n}(f)\),
\[
r_{n,f}
\le
R_{n,f}^{\rm lin}
+
\max_{\lambda\in\Lambda_n}
\frac{|b_{\lambda,n}(f)|}{s_{\lambda,f}}.
\]
The second summand is \(\mathcal F_{0n}\)-measurable.  Therefore its
conditional exceedance probability is either zero or one, and
\eqref{eq:t3-10}--\eqref{eq:t3-12} imply
\begin{equation}
r_{n,f}
=
o_{\rm ucp}(\ell_n^{-1/2})
\quad
\text{on }\widetilde{\mathcal G}_n.
\label{eq:t3-29}
\end{equation}

For every fixed \(\lambda\), conditional Bernstein's inequality gives
\[
P_f\left(
|G_\lambda|>x
\mid\mathcal F_{0n}
\right)
\le
2\exp\left\{
-\frac{x^2}{2(1+L_nx/3)}
\right\}
\]
on \(\widetilde{\mathcal G}_n\).  Hence, for every
\(K_0>\sqrt2\),
\begin{align}
&P_f\left(
\max_{\lambda\in\Lambda_n}|G_\lambda|
>
K_0\sqrt{\ell_n}
\mid\mathcal F_{0n}
\right)
\notag\\
&\qquad\le
2m_n
\exp\left\{
-\frac{K_0^2\ell_n}
{2\{1+K_0L_n\sqrt{\ell_n}/3\}}
\right\}
=o(1),
\label{eq:t3-30}
\end{align}
uniformly on \(\widetilde{\mathcal G}_n\).  Here
\(L_n\sqrt{\ell_n}\to0\) follows from \eqref{eq:t3-9}.  On the
complement of the event in \eqref{eq:t3-30},
\begin{equation}
\max_{\lambda\in\Lambda_n}
\frac{
|\widehat\theta_\lambda-\theta_\lambda(f)|
}{
s_{\lambda,f}
}
\le
K_0\sqrt{\ell_n}+r_{n,f}.
\label{eq:t3-31}
\end{equation}

For one ideal multiplier draw, conditional Gaussianity gives
\[
P_{\mathcal X_n}^*(T^{*(1)}>x)
\le
2m_ne^{-x^2/2}.
\]
Let
\[
x_{\alpha,n}
=
\sqrt{2\log(4m_n/\alpha)}.
\]
Then
\[
P_{\mathcal X_n}^*
(T^{*(1)}>x_{\alpha,n})
\le\alpha/2.
\]
Conditional on \(\mathcal X_n\), the number of the \(B_n^*\) draws
exceeding \(x_{\alpha,n}\) is binomial with success probability at
most \(\alpha/2\).  Thus, for some \(c_\alpha>0\),
\begin{equation}
P_{\mathcal X_n}^{*,B}
\{c_{1-\alpha}^*>x_{\alpha,n}\}
\le
e^{-c_\alpha B_n^*}.
\label{eq:t3-32}
\end{equation}
Since
\[
x_{\alpha,n}/\sqrt{\ell_n}\longrightarrow\sqrt2,
\]
every fixed \(K_c>\sqrt2\) satisfies, for all sufficiently large \(n\),
\begin{align}
&P_f^\dagger\left(
c_{1-\alpha}^*>K_c\sqrt{\ell_n}
\mid\mathcal F_{0n}
\right)
\notag\\
&\quad=
E_f\left[
P_{\mathcal X_n}^{*,B}
\{c_{1-\alpha}^*>K_c\sqrt{\ell_n}\}
\mid\mathcal F_{0n}
\right]
\notag\\
&\quad\le
E_f\left[
P_{\mathcal X_n}^{*,B}
\{c_{1-\alpha}^*>x_{\alpha,n}\}
\mid\mathcal F_{0n}
\right]
\le
e^{-c_\alpha B_n^*}.
\label{eq:t3-33}
\end{align}

The standardized noncentrality is
\[
\nu_n
=
\frac{a_n h_n}{s_{\lambda_0,f}}
=
a_n h_n\sqrt{I_{\lambda_0,f}}.
\]
By \eqref{eq:t3-7} and \eqref{eq:t3-17},
\[
\nu_n^2
\ge
\underline\iota C_{\rm up}\ell_n.
\]
Choose
\[
K_s>\max\{2K_0,K_0+K_c\},
\qquad
C_{\rm up}>
\frac{K_s^2}{\underline\iota}.
\]
Then
\begin{equation}
\nu_n
\ge
K_s\sqrt{\ell_n}.
\label{eq:t3-34}
\end{equation}

On
\[
\left\{
\max_\lambda|G_\lambda|
\le K_0\sqrt{\ell_n}
\right\}
\cap
\{\delta_{n,f}<1\},
\]
we have the exact bounds
\begin{align}
|Z_{\lambda_0}|
&\ge
\frac{
\nu_n-K_0\sqrt{\ell_n}-r_{n,f}
}{
1+\delta_{n,f}
},
\label{eq:t3-35}\\
\max_{\mu\ne\lambda_0}|Z_\mu|
&\le
\frac{
K_0\sqrt{\ell_n}+r_{n,f}
}{
1-\delta_{n,f}
}.
\label{eq:t3-36}
\end{align}

The preceding choice implies both
\(K_s-K_0>K_c\) and \(K_s-K_0>K_0\).
By continuity at zero, there is a constant
\(\bar\delta\in(0,1)\) such that
\[
\frac{K_s-K_0}{1+\bar\delta}
>
\max\left\{
K_c,\,
\frac{K_0}{1-\bar\delta}
\right\}.
\]
For all sufficiently large \(n\), the conditional failure probability
of simultaneous rejection, exact atom recovery, and sign recovery is
bounded by
\begin{align*}
&
P_f\left(
\max_\lambda|G_\lambda|>K_0\sqrt{\ell_n}
\mid\mathcal F_{0n}
\right)
+
P_f(r_{n,f}>1\mid\mathcal F_{0n})\\
&\qquad
+
P_f(\delta_{n,f}>\bar\delta\mid\mathcal F_{0n})
+
e^{-c_\alpha B_n^*}.
\end{align*}
By \eqref{eq:t3-29}, \eqref{eq:t3-30}, \eqref{eq:t3-11}, and
\eqref{eq:t3-33}, this bound converges to zero in the uniform
conditional sense defined before \eqref{eq:t3-7}.  On its complement,
\[
|Z_{\lambda_0}|
>
\max\left\{
c_{1-\alpha}^*,\
\max_{\mu\ne\lambda_0}|Z_\mu|
\right\},
\]
so
\[
\phi_n=1,
\qquad
\widehat\lambda=\lambda_0.
\]
Moreover,
\[
\frac{
\varsigma_0\widehat\theta_{\lambda_0}
}{
s_{\lambda_0,f}
}
\ge
\nu_n-K_0\sqrt{\ell_n}-r_{n,f}>0,
\]
so \(\widehat\varepsilon=\varsigma_0\).  This proves
\eqref{eq:t3-18}.

We next prove \eqref{eq:t3-16}.  Put
\[
\epsilon_n^*=(B_n^*)^{-1/4}.
\]
For all sufficiently large \(n\),
\[
\alpha+\epsilon_n^*+\kappa_n<1.
\]
Define
\[
q_{p,n}^{0}(f_0)
=
\inf\{t:F_{n,f_0}^{0}(t)\ge p\}.
\]
On
\(\{\Delta_{n,f_0}^{\rm cal}\le\kappa_n\}\),
\[
q_{1-\alpha-\epsilon_n^*,n}^*
\ge
q_{1-\alpha-\epsilon_n^*-\kappa_n,n}^{0}(f_0).
\]
Indeed, if
\(t<q_{1-\alpha-\epsilon_n^*-\kappa_n,n}^{0}(f_0)\), then
\[
F_n^*(t)
\le
F_{n,f_0}^{0}(t)+\kappa_n
<
1-\alpha-\epsilon_n^*.
\]
Consequently,
\begin{align*}
&
\left\{
\widetilde{\mathcal G}_n,\
T_n>q_{1-\alpha-\epsilon_n^*,n}^*,\
\Delta_{n,f_0}^{\rm cal}\le\kappa_n
\right\}\\
&\qquad\subseteq
\left\{
\widetilde{\mathcal G}_n,\
T_n>
q_{1-\alpha-\epsilon_n^*-\kappa_n,n}^{0}(f_0)
\right\}.
\end{align*}
Since \(\widetilde{\mathcal G}_n\) and the conditional null quantile
are \(\mathcal F_{0n}\)-measurable, \eqref{eq:t3-15} gives
\[
\sup_{f_0\in\mathcal B_n}
P_{f_0}\left[
\widetilde{\mathcal G}_n
\cap
\left\{
T_n>
q_{1-\alpha-\epsilon_n^*,n}^*
\right\}
\right]
\le
\alpha+\epsilon_n^*+\kappa_n+o(1).
\]

Conditional on \(\mathcal X_n\), DKW gives
\[
P_{\mathcal X_n}^{*,B}\left(
\sup_t
|\widehat F_{B_n^*}^*(t)-F_n^*(t)|
>
\epsilon_n^*
\right)
\le
2e^{-2\sqrt{B_n^*}}.
\]
On the complementary event,
\[
c_{1-\alpha}^*
\ge
q_{1-\alpha-\epsilon_n^*,n}^*.
\]
Therefore
\[
\sup_{f_0\in\mathcal B_n}
P_{f_0}^\dagger\left[
\widetilde{\mathcal G}_n
\cap
\{T_n>c_{1-\alpha}^*\}
\right]
\le
\alpha+o(1).
\]
Admission failure forces \(\phi_n=0\), so
\[
P_{f_0}^\dagger(\phi_n=1)
\le
P_{f_0}^\dagger\left[
\widetilde{\mathcal G}_n
\cap
\{T_n>c_{1-\alpha}^*\}
\right]
+
P_{f_0}(\mathcal G_n^c).
\]
Equation \eqref{eq:t3-5} proves \eqref{eq:t3-16}.

Finally, put
\[
q_{n,f}
=
P_f^\dagger
(\mathcal E_{n,f}^{\rm rec}\mid\mathcal F_{0n}).
\]
For every fixed \(\delta>0\),
\begin{align*}
E_f^\dagger
\left[
\mathbf 1_{\widetilde{\mathcal G}_n}
\mathbf 1_{\mathcal E_{n,f}^{\rm rec}}
\right]
&=
E_f\left[
\mathbf 1_{\widetilde{\mathcal G}_n}q_{n,f}
\right]\\
&\le
\delta
+
P_f\left[
\widetilde{\mathcal G}_n
\cap
\{q_{n,f}>\delta\}
\right].
\end{align*}
Taking the supremum, then the \(\limsup\), and then letting
\(\delta\downarrow0\), \eqref{eq:t3-18} gives
\[
\sup_{f\in\mathfrak A_n(a_n,h_n)}
P_f^\dagger(\mathcal E_{n,f}^{\rm rec})
\le
\eta_n+\gamma_n+o(1).
\]
Since the loss is bounded by one and vanishes under exact recovery,
\eqref{eq:t3-19} follows.

We now prove the lower assertions.  Restrict the null and alternative
classes to the baseline \(f_0^\circ\), and put
\[
P_0:=P_{f_0^\circ}^{\rm orc}.
\]
Every mixture element below belongs to
\(\Theta_{1,n}(\mathbf a_n)\), because
\[
\left\|
f_0^\circ+\varsigma a_h h\Psi_\lambda^S
\right\|_\infty
\le
(B-\delta_B)+\delta_B
=
B.
\]
It remains symmetric and mean zero because \(f_0^\circ\) and the Haar
atoms have these properties.

Let \(R_{n,\circ}^{\rm det,orig}\) and
\(R_{n,\circ}^{\rm loc,orig}\) denote the original-experiment risks
after fixing \(f_0=f_0^\circ\), and let
\(R_{n,\circ}^{\rm det,orc}\) and
\(R_{n,\circ}^{\rm loc,orc}\) be the corresponding oracle risks.
Parameter restriction and Blackwell dominance give
\[
R_n^{\rm det}
\ge
R_{n,\circ}^{\rm det,orig}
\ge
R_{n,\circ}^{\rm det,orc},
\qquad
R_n^{\rm loc}
\ge
R_{n,\circ}^{\rm loc,orig}
\ge
R_{n,\circ}^{\rm loc,orc}.
\]

For \(\lambda\in\Lambda_n\) and
\(\varsigma\in\{-1,1\}\), define
\[
t_{\lambda,\varsigma,ij}
=
\varsigma a_{h_\lambda}v_{\lambda,ij},
\qquad
p_{\lambda,\varsigma,ij}
=
\sigma\left\{
\operatorname{logit}(p_{0,ij})
+
t_{\lambda,\varsigma,ij}
\right\}.
\]
Let \(P_{\lambda,\varsigma}\) be the corresponding product-Bernoulli
oracle law and put
\[
L_{\lambda,\varsigma}
=
\frac{dP_{\lambda,\varsigma}}{dP_0},
\qquad
\Delta_{\lambda,\varsigma,ij}
=
p_{\lambda,\varsigma,ij}-p_{0,ij}.
\]
Direct calculation gives
\begin{equation}
E_0
(L_{\lambda,\varsigma}L_{\mu,\varsigma'})
=
\prod_{(i,j)\in\mathcal K_n}
\left[
1+
\frac{
\Delta_{\lambda,\varsigma,ij}
\Delta_{\mu,\varsigma',ij}
}{
p_{0,ij}(1-p_{0,ij})
}
\right].
\label{eq:t3-37}
\end{equation}

For every fixed \(M<\infty\), there is \(C_M<\infty\) such that
\begin{equation}
|\sigma(\eta+u)-\sigma(\eta)|
\le
C_M\sigma(\eta)\{1-\sigma(\eta)\}|u|,
\qquad |u|\le M.
\label{eq:t3-38}
\end{equation}
This follows because
\[
\sup_{\eta\in\mathbb R}
\sup_{|r|\le M}
\frac{\sigma'(\eta+r)}{\sigma'(\eta)}
<\infty.
\]
By \eqref{eq:t3-20}--\eqref{eq:t3-22},
\(|t_{\lambda,\varsigma,ij}|\le\delta_B\).  Hence
\begin{align}
\log E_0L_{\lambda,\varsigma}^2
&\le
C a_{h_\lambda}^2n^2\rho_nh_\lambda^2,
\label{eq:t3-39}\\
\operatorname{KL}(P_{\lambda,\varsigma},P_0)
&\le
C a_{h_\lambda}^2n^2\rho_nh_\lambda^2.
\label{eq:t3-40}
\end{align}
Indeed,
\[
\log E_0L_{\lambda,\varsigma}^2
\le
\sum_{(i,j)\in\mathcal K_n}
\frac{\Delta_{\lambda,\varsigma,ij}^2}
{p_{0,ij}(1-p_{0,ij})},
\]
and the same sum bounds the product-Bernoulli KL divergence.  By
\eqref{eq:t3-38},
\[
\frac{\Delta_{\lambda,\varsigma,ij}^2}
{p_{0,ij}(1-p_{0,ij})}
\le
C\rho_n a_{h_\lambda}^2v_{\lambda,ij}^2.
\]
Summation and \eqref{eq:t3-21} prove
\eqref{eq:t3-39}--\eqref{eq:t3-40}.

If
\(\mathcal D_\lambda\cap\mathcal D_\mu=\varnothing\), every factor in
\eqref{eq:t3-37} equals one, so
\begin{equation}
E_0
(L_{\lambda,\varsigma}L_{\mu,\varsigma'})
=1.
\label{eq:t3-41}
\end{equation}
For an overlapping pair, put
\[
x_{ij}
=
\frac{
\Delta_{\lambda,\varsigma,ij}
\Delta_{\mu,\varsigma',ij}
}{
p_{0,ij}(1-p_{0,ij})
}.
\]
Each \(1+x_{ij}>0\), so
\(\log(1+x_{ij})\le x_{ij}\le|x_{ij}|\).  Thus
\begin{align}
\log E_0
(L_{\lambda,\varsigma}L_{\mu,\varsigma'})
&\le
C\rho_n a_{h_\lambda}a_{h_\mu}
\sum_{(i,j)\in\mathcal K_n}
|v_{\lambda,ij}v_{\mu,ij}|
\notag\\
&\le
C\sqrt{e_{\lambda,n}e_{\mu,n}},
\label{eq:t3-42}
\end{align}
where
\[
e_{\lambda,n}
=
a_{h_\lambda}^2n^2\rho_nh_\lambda^2.
\]

Consider
\[
Q_n
=
\frac1{2m_n}
\sum_{\lambda\in\Lambda_n}
\sum_{\varsigma\in\{-1,1\}}
P_{\lambda,\varsigma}.
\]
Writing \([x]_+=\max(x,0)\), equations
\eqref{eq:t3-41}--\eqref{eq:t3-42} give
\begin{align}
\chi^2(Q_n,P_0)
&=
\frac1{4m_n^2}
\sum_{\lambda,\mu}
\sum_{\varsigma,\varsigma'}
\left\{
E_0
(L_{\lambda,\varsigma}L_{\mu,\varsigma'})
-1
\right\}
\notag\\
&\le
\frac1{4m_n^2}
\sum_{\substack{\lambda,\mu:\\
\mathcal D_\lambda\cap\mathcal D_\mu\ne\varnothing}}
\sum_{\varsigma,\varsigma'}
\left[
E_0
(L_{\lambda,\varsigma}L_{\mu,\varsigma'})
-1
\right]_+
\notag\\
&\le
\frac{N_{{\rm ov},n}}{m_n^2}
\exp\left\{
C\max_{\lambda\in\Lambda_n}e_{\lambda,n}
\right\}
\notag\\
&\le
\frac{d_n}{m_n}
\exp\left\{
C\max_{\lambda\in\Lambda_n}e_{\lambda,n}
\right\}.
\label{eq:t3-43}
\end{align}
The positive-part step is required because opposite-sign cross terms
can be negative.

Under \eqref{eq:t3-26},
\[
\chi^2(Q_n,P_0)
\le
C(2m_n)^{-1+Cc_{\rm low}+o(1)}
\longrightarrow0
\]
for sufficiently small \(c_{\rm low}>0\).  Hence
\[
\|Q_n-P_0\|_{\rm TV}
\le
\frac12\chi^2(Q_n,P_0)^{1/2}
\longrightarrow0.
\]
For every oracle test \(\varphi\),
\[
P_0(\varphi=1)
+
\sup_{\lambda,\varsigma}
P_{\lambda,\varsigma}(\varphi=0)
\ge
1-\|Q_n-P_0\|_{\rm TV}.
\]
Therefore
\[
R_{n,\circ}^{\rm det,orc}
\ge
1-\|Q_n-P_0\|_{\rm TV}.
\]
The experiment comparison gives
\[
\liminf_{n\to\infty}R_n^{\rm det}\ge1.
\]
The test that never rejects has full-experiment risk one, so
\(R_n^{\rm det}\le1\).  Thus
\[
R_n^{\rm det}\longrightarrow1.
\]

For localization, define
\[
P_{\lambda,+}:=P_{\lambda,+1}.
\]
Let \(V\) be uniform on \(\Lambda_n\), and, conditionally on
\(V=\lambda\), generate the oracle data \(Y\) under \(P_{\lambda,+}\).
If
\[
\overline P_n
=
\frac1{m_n}
\sum_{\lambda\in\Lambda_n}P_{\lambda,+},
\]
then
\begin{align}
I(V;Y)
&=
\frac1{m_n}
\sum_{\lambda\in\Lambda_n}
\operatorname{KL}(P_{\lambda,+},\overline P_n)
\notag\\
&\le
\frac1{m_n}
\sum_{\lambda\in\Lambda_n}
\operatorname{KL}(P_{\lambda,+},P_0)
\notag\\
&\le
C\max_{\lambda\in\Lambda_n}e_{\lambda,n}.
\label{eq:t3-44}
\end{align}
Fano's inequality gives
\begin{equation}
\inf_{\check\lambda}
\sup_{\lambda\in\Lambda_n}
P_{\lambda,+}(\check\lambda\ne\lambda)
\ge
1-
\frac{
C\max_\lambda e_{\lambda,n}+\log2
}{
\log m_n
}.
\label{eq:t3-45}
\end{equation}
Since \(\ell_n/\log m_n\to1\), reducing \(c_{\rm low}\) if necessary
makes the right-hand side at least \(c_F>0\).

Put
\[
r_0
=
\frac14\min\{1,\delta_{\rm pack}\},
\qquad
\mathcal C_\lambda
=
\left\{
D\in\mathcal D_n^{\rm dec}:
d_{\rm loc}(D,\lambda)<r_0
\right\}.
\]
These sets are Borel and pairwise disjoint.  Indeed, common membership
first forces
\[
h_D=h_\lambda=h_\mu.
\]
If \(A=|S_\lambda|=|S_\mu|\), the triangle inequality gives
\[
|S_\lambda\triangle S_\mu|<4r_0A,
\]
whereas \eqref{eq:t3-2} gives
\[
|S_\lambda\triangle S_\mu|
\ge2\delta_{\rm pack}A.
\]
Since \(4r_0\le\delta_{\rm pack}<2\delta_{\rm pack}\), this is
impossible.

Decode \(D\) as the unique \(\lambda\) with
\(D\in\mathcal C_\lambda\), if one exists, and use a fixed
deterministic fallback otherwise.  This decoder is measurable because
\(\Lambda_n\) is finite and the sets \(\mathcal C_\lambda\) are
Borel.  Moreover,
\[
d_{\rm loc}(D,\lambda)
\ge
r_0\mathbf 1\{\check\lambda(D)\ne\lambda\}.
\]
Fano's inequality and the experiment comparison now yield
\[
\liminf_{n\to\infty}R_n^{\rm loc}
\ge
r_0c_F>0.
\]
This proves \eqref{eq:t3-27}.
\end{proof}

\begin{remark}[Global-union versus scale-wise lower bounds]
For a deterministic sequence \(h_n\in\mathcal J_n\), put
\[
m_{n,h_n}=|\Lambda_n(h_n)|
\]
and let \(N_{{\rm ov},n}(h_n)\) be the ordered overlap count restricted
to \(\Lambda_n(h_n)^2\).  Let
\(R_{n,h_n}^{\rm det}\) and \(R_{n,h_n}^{\rm loc}\) denote the risks
with alternatives restricted to scale \(h_n\).

If
\[
m_{n,h_n}\to\infty,
\qquad
N_{{\rm ov},n}(h_n)
\le
m_{n,h_n}d_{n,h_n},
\qquad
d_{n,h_n}\ge1,
\]
and
\[
\frac{\log d_{n,h_n}}
{\log(2m_{n,h_n})}
\longrightarrow0,
\]
then there is a constant \(c_0>0\) such that
\[
a_{h_n}^2n^2\rho_nh_n^2
\le
c_0\log(2m_{n,h_n})
\]
implies
\[
R_{n,h_n}^{\rm det}\longrightarrow1,
\qquad
\liminf_{n\to\infty}R_{n,h_n}^{\rm loc}>0.
\]

A uniform scale-by-scale conclusion with \(\log(2m_n)\) additionally
requires
\[
\min_{h\in\mathcal J_n}m_{n,h}\longrightarrow\infty,
\]
\[
\sup_{h\in\mathcal J_n}
\frac{\log d_{n,h}}
{\log(2m_{n,h})}
\longrightarrow0,
\]
and, for some \(\kappa>0\),
\[
\liminf_{n\to\infty}
\min_{h\in\mathcal J_n}
\frac{\log(2m_{n,h})}{\log(2m_n)}
\ge\kappa,
\]
where
\[
N_{{\rm ov},n}(h)
\le
m_{n,h}d_{n,h}.
\]
Without these conditions, the correct conclusion is the global-union
minimax theorem without a matching lower bound at every scale.
\end{remark}

\begin{remark}[A sufficient global overlap condition]
Suppose the implemented packing has pairwise disjoint supports within
each scale and has the following verified laminar property: every atom
intersects at most \(C_0\) atoms at each coarser scale, where \(C_0\)
does not depend on \(n\), the atom, or the scale.  Charging each
ordered overlapping pair to its finer atom and the coarser scale then
gives
\[
N_{{\rm ov},n}
\le
C m_n(1+|\mathcal J_n|).
\]
Consequently, \eqref{eq:t3-23} holds whenever
\[
\log(1+|\mathcal J_n|)
=
o\{\log m_n\}.
\]
For a particular symmetric Haar packing, the disjointness and
bounded-ancestor properties must be proved as part of the packing
construction or replaced by a direct verification of
\eqref{eq:t3-23}.  Alternatively, the lower mixture may be restricted
to a low-overlap subpacking whose logarithmic cardinality is
comparable to \(\log m_n\).
\end{remark}

\section{Unconditional instance}\label{sec:supp-uncond}


\paragraph{Setting.}
Throughout this section the following instance is fixed.  Constants denoted $C$
or $c$ depend only on $(B,C_S,C_\Psi,\kappa_0,\alpha,J_0,j_0)$ and may change
between displays.

\begin{enumerate}[label=(S\arabic*),leftmargin=2.4em,itemsep=2pt]
\item\label{S:design} \emph{Regular known design.}  $n=2^{K_n}$ and
$U_i=x_i=(i-\tfrac12)/n$; the coordinate is known exactly, so
$\widehat U=U$, $\Delta_n=0$, and the reported radius is $r_n=0$, whence
$\eta_n=0$ in \eqref{eq:T4-coverage} and in the design events of the preceding
theorems.
\item\label{S:obs} \emph{Complete observation.}  $M_{ij}\equiv1$ and
$\pi_{ij}\equiv1$, known; hence $W_k^F=R_k^F$.
\item\label{S:folds} \emph{Parity folds.}  Dyad $(i,j)$, $i<j$, is assigned to
the inference fold $I$ when $i+j$ is even and to the baseline-estimation fold
$B^E$ when $i+j$ is odd.  The assignment is deterministic and outcome-free.
\item\label{S:cells} \emph{Dyadic cells and dictionary.}  Fine cells have width
$\delta_n=2^{-J_n}$ with $J_n\le K_n-2$, so each one-dimensional bin contains
$n\delta_n=2^{K_n-J_n}\ge4$ design points.  The coarse space $V_{J_0}$ has cell
width $2^{-J_0}$.  The dictionary $\Lambda_n$ consists of all aligned symmetric
\emph{off-diagonal} Haar atoms $\Psi_\lambda^S$ of the main text (distinct
one-dimensional index pairs, so every atom has support area $2h_\lambda^2$ at
its scale; diagonal atoms have area $h_\lambda^2$ and are excluded from this
instance, and the mixed dictionary of the implementation is outside this
theorem) at scales
$j\in\{j_0,\dots,j_{\max,n}\}$ with $J_0<j_0$ and
$J_n\ge j_{\max,n}+1$, so $2\delta_n\le h_{\min,n}$; write
$h_{\max}=2^{-j_0}$ and $m_n=|\Lambda_n|\le4^{\,j_{\max,n}}\le n^2$.
\item\label{S:baseline} \emph{Frozen baseline map.}  $B_n(y)=
\sigma\{\Pi_{J_0}\operatorname{logit}(\operatorname{clip}_{\epsilon_n}(S_ny))\}$,
where $S_n$ is a fixed pre-smoother with nonnegative entries, $S_n\mathbf 1=
\mathbf 1$, and column sums bounded by $C_S$; $\Pi_{J_0}$ is the exact discrete
Lebesgue projection onto $V_{J_0}$ computed with the symmetric-cell areas
$|Q_q|$; and $\operatorname{clip}_{\epsilon}$ truncates to
$[\epsilon,1-\epsilon]$.  The identity $S_n=I$ is admissible.
\item\label{S:stab} \emph{Stabilization and clipping.}  $w_q=N_q^I/(N_q^I+
\kappa_0)$ with $\kappa_0>0$ fixed and $N_q^I=D_q^I$ (unit weights), and
$\epsilon_n\in(0,e^{-B}\rho_n/8]$.
\item\label{S:model} \emph{Model class.}  Under the null,
$p_{ij}=\sigma\{\alpha_n+f_0(x_i,x_j)\}$ with $f_0\in V_{J_0}$ symmetric,
mean zero, $\|f_0\|_\infty\le B$; single-atom alternatives replace $f_0$ by
$f_0+\varsigma a h\Psi_\lambda^S$ with $\|f_0+\varsigma ah\Psi_\lambda^S\|_\infty
\le B$, $\lambda\in\Lambda_n$.  Sparsity: $\rho_n=\sigma(\alpha_n)\to0$.  Then
every edge probability lies in $[c_\rho\rho_n,C_\rho\rho_n]$ with
$c_\rho=e^{-B}/2$ and $C_\rho=2e^{B}$ for all large $n$.
\item\label{S:gates} \emph{Deterministic admission.}  The admitted family is the
full dictionary: the pilot gates are evaluated on the known design and known
propensities, which is outcome-free, and their thresholds are met for every
$\lambda\in\Lambda_n$ by Lemma~\ref{lem:U-counts} below; equivalently the
all-pass event of Theorem~\ref{thm:adaptive-detection-localization} holds with
$\gamma_n=0$.
\end{enumerate}

\paragraph{Primitive rate conditions.}
Fix $\gamma>0$.  Assume
\begin{align}
&\text{(C1)}\qquad
n^2\rho_n\delta_n^2\;\ge\;\log^{7+\gamma}(n),
\label{eq:U-C1}\\
&\text{(C2)}\qquad
\frac{h_{\max}\log^{3/2}(n)}{\sqrt{\rho_n}\,n\,\delta_n^{2}}
\;\longrightarrow\;0,
\label{eq:U-C2}\\
&\text{(C3)}\qquad
\kappa_0\sqrt{\rho_n\log n}\;=\;o\bigl(n\delta_n^{2}\bigr).
\label{eq:U-C3}
\end{align}
For the detection--localization clause additionally assume
\begin{equation}
\text{(C4)}\qquad
j_{\max,n}\to\infty,
\qquad
\frac{j_{\max,n}}{\log\log n}\to\infty.
\label{eq:U-C4}
\end{equation}
Conditions \eqref{eq:U-C1}--\eqref{eq:U-C3} are compatible: for example,
$j_{\max,n}=\lfloor c\log_2n\rfloor$ with $c<1/4$, $\delta_n=2^{-j_{\max,n}-1}$,
and expected degree $n\rho_n\ge\log^{8}(n)$ satisfy all four displays.

\paragraph{Deterministic counts.}

\begin{lemma}[Exact fold counts and gates]\label{lem:U-counts}
Under \ref{S:design}--\ref{S:cells}, each one-dimensional bin contains exactly
$n\delta_n$ design points, of which exactly $n\delta_n/2$ have even index and
$n\delta_n/2$ odd.  Consequently, for every off-diagonal symmetric fine cell
$q$ and each fold $F\in\{I,B^E\}$,
\[
D_q^F=\tfrac12(n\delta_n)^2,
\]
and for every diagonal cell,
$D_q^F\in\bigl[\tfrac14(n\delta_n)^2-\tfrac12 n\delta_n,\
\tfrac14(n\delta_n)^2\bigr]$.
In particular $D_q^F\ge\tfrac18(n\delta_n)^2=:D_{\min}$ for $n\delta_n\ge4$,
every cell in every atom's support is supported in both folds, and all
information, leverage, and order gates of the main text, evaluated at these
deterministic counts with $r_n=0$, pass for every $\lambda\in\Lambda_n$ once
$n^2\rho_n\delta_n^2$ exceeds a constant multiple of the fixed thresholds,
which \eqref{eq:U-C1} guarantees for all large $n$.
\end{lemma}

\begin{proof}
A bin is $\{i:(a-1)n\delta_n<i\le an\delta_n\}$, a block of $n\delta_n$
consecutive integers; since $n\delta_n$ is even, exactly half its members are
even.  For an off-diagonal cell generated by bins $a<b$, every pair
$(i,j)\in\text{bin}_a\times\text{bin}_b$ has $i<j$, and $i+j$ is even iff $i,j$
share parity, giving $2(n\delta_n/2)^2=\tfrac12(n\delta_n)^2$ dyads per fold.
For a diagonal cell the unordered same-parity pairs number
$2\binom{n\delta_n/2}{2}$ and the mixed-parity pairs $(n\delta_n/2)^2$; both lie
in the stated interval.  The gate claims follow because every pilot gate quantity is, under \ref{S:obs}, a fixed multiple of these counts and of
$\rho_n$-scale probabilities, and the order gate is vacuous at $r_n=0$.
\end{proof}

\paragraph{Exact algebraic structure.}

\begin{lemma}[Exact identities]\label{lem:U-exact}
Under \ref{S:design}--\ref{S:model}:
\begin{enumerate}[label=\textup{(\alph*)},leftmargin=2em,itemsep=2pt]
\item Every surface $\alpha_n+f$ in the class is constant on each symmetric fine
cell; hence $p_k=p_q$ for all dyads $k$ in cell $q$ and, exactly,
$\mu_q^F=p_q$ for both folds.
\item With $c_{\lambda q}=|Q_q|\overline\Psi_{\lambda q}$, the vector
$c_\lambda$ satisfies $\Pi_{J_0}^\top c_\lambda=0$; consequently
$\sum_qc_{\lambda q}\operatorname{logit}\{B_{n,q}(y)\}=0$ for every $y$ in the
clip-free range, and the frozen functional collapses to
\[
\Phi_{\lambda n}(x,y)
=\sum_{q}c_{\lambda q}\,
\ell_{\epsilon_n}\{w_qx_q+(1-w_q)B_{n,q}(y)\}.
\]
\item The population coefficient is computed exactly by the cell sum:
$\theta_\lambda(P)=\sum_qc_{\lambda q}\operatorname{logit}(p_q)$, which equals
$\varsigma ah\,\mathbf 1\{\lambda=\lambda_*\}$ under a single-atom alternative
and $0$ under the null.
\item The finite-design target obeys, for all large $n$,
\[
\bigl|\theta_{\lambda,n}^{D}-\theta_\lambda(P)\bigr|
\;\le\;
\frac{C\,\kappa_0\,h_\lambda}{D_{\min}},
\qquad\text{uniformly over the class and over }\lambda\in\Lambda_n .
\]
If $S_n=I$, then in addition $\theta_{\lambda,n}^{D}=\theta_\lambda(P)$ exactly
under the null.
\end{enumerate}
\end{lemma}

\begin{proof}
(a) Elements of $V_{J_0}$ are constant on squares of width $2^{-J_0}$, and
$\Psi_\lambda^S$ at scale $j$ is constant on squares of width $2^{-j-1}$; since
$J_n\ge j_{\max,n}+1$ and $J_n\ge J_0$, both are constant on fine cells, and the
symmetrized cell value is unambiguous.  Design points are cell midpointwise
interior, so all dyads in a cell share one probability, and
$\mu_q^F=\sum_kb^F_{kq}p_k=p_q$ because the weights sum to one.

(b) For cellwise-constant vectors the discrete projection with weights
$|Q_q|$ coincides with the $L^2([0,1]^2)$ projection; writing $C(r)$ for the
coarse cell containing $r$,
\[
(\Pi_{J_0}^\top c_\lambda)_r
=\frac{|Q_r|}{|C(r)|}\sum_{q\in C(r)}|Q_q|\overline\Psi_{\lambda q}
=\frac{|Q_r|}{|C(r)|}\,
\langle\Psi_\lambda^S,\mathbf 1_{C(r)}\rangle_{L^2}=0,
\]
since $\mathbf 1_{C(r)}\in V_{J_0}$ and $\Psi_\lambda^S\perp V_{J_0}$.
Hence for any $y$ with $\operatorname{clip}$ inactive,
$\sum_qc_{\lambda q}\{\Pi_{J_0}\operatorname{logit}(S_ny)\}_q
=(\Pi_{J_0}^\top c_\lambda)^\top\operatorname{logit}(S_ny)=0$; on a constant
clipping branch the same sum is $(\Pi^\top c)^\top$ applied to a vector that is
constant in the clipped coordinates and unchanged elsewhere, and vanishes by the
same identity applied coordinatewise.  Subtracting this null functional from
\eqref{eq:frozen-functional-final} gives the collapsed form.

(c) By cellwise constancy,
$\theta_\lambda(P)=\langle f,\Psi_\lambda^S\rangle
=\sum_q|Q_q|\overline\Psi_{\lambda q}f(q)=\sum_qc_{\lambda q}f(q)$; and
$\operatorname{logit}(p_q)=\alpha_n+f(q)$, while
$\sum_qc_{\lambda q}\{\alpha_n+f_0(q)\}=0$ because constants and $f_0$ lie in
$V_{J_0}$.  Orthonormality of the atoms gives the stated values.

(d) By \ref{S:baseline} and \ref{S:model}, $S_n\mu^{B^E}$ is a convex
combination of values in $[c_\rho\rho_n,C_\rho\rho_n]$, the clip is inactive
because $\epsilon_n\le c_\rho\rho_n/4$, the logits lie in
$[\alpha_n+\log c_\rho+O(\rho_n),\ \alpha_n+\log C_\rho+O(\rho_n)]$, their
coarse averages lie in the same interval, and therefore
$B_{n,q}(\mu^{B^E})\in[\tfrac12c_\rho\rho_n,\,2C_\rho\rho_n]$ for all large
$n$.  Hence
$z_q^0=w_qp_q+(1-w_q)B_{n,q}(\mu^{B^E})\in
[\tfrac12c_\rho\rho_n,\,2C_\rho\rho_n]$, the clip is inactive at the target,
and the mean value theorem for $\operatorname{logit}$ on this interval gives
\[
\bigl|\ell_{\epsilon_n}(z_q^0)-\operatorname{logit}(p_q)\bigr|
\le\frac{(1-w_q)\,|B_{n,q}(\mu^{B^E})-p_q|}{\tfrac12c_\rho\rho_n
\{1-2C_\rho\rho_n\}}
\le C(1-w_q)\le\frac{C\kappa_0}{D_{\min}} .
\]
Multiplying by $|c_{\lambda q}|$ and summing, with
$\sum_q|c_{\lambda q}|=\|\Psi_\lambda^S\|_{L^1}\le
\|\Psi_\lambda^S\|_{L^2}\,|{\rm supp}|^{1/2}\le Ch_\lambda$, proves the bound.
If $S_n=I$, then under the null $S_n\mu^{B^E}=p$, the logit vector is
$\alpha_n+f_0\in V_{J_0}$, the projection reproduces it,
$B_{n,q}(\mu^{B^E})=p_q$, $z_q^0=p_q$, and every summand vanishes.
\end{proof}

\paragraph{The Bernstein rectangle.}

\begin{lemma}[Rectangle, knot margin, and smoothness]\label{lem:U-rect}
Adopt the notation of Theorem~\ref{thm:exact-finite-design-reduction} with the
active cells of Lemma~\ref{lem:U-counts}.  Under \eqref{eq:U-C1}, for all large
$n$, uniformly over the class and over active cells,
\[
\frac{t_\alpha}{\mu_\alpha}
\;\le\;
C\sqrt{\frac{\log n}{n^2\rho_n\delta_n^2}}
\;\le\;C\log^{-(3+\gamma)/2}n\;\longrightarrow\;0 .
\]
Consequently $\mathcal C_{n,P}\subset[\tfrac12c_\rho\rho_n,\,2C_\rho\rho_n]^{
D_n^{\rm act}}$, the clip is inactive on an open neighborhood
$\mathcal U_{n,P}$ of the rectangle with margin
$\kappa_{n,P}\ge c_\rho\rho_n/4$, the pre-smoothed vectors $S_ny$ remain in the
same range, condition \eqref{eq:baseline-interior-condition-final} holds, and
$B_n$ and $\Phi_{\lambda n}$ are twice continuously differentiable on
$\mathcal U_{n,P}$.  Moreover, on the rectangle,
\[
\frac{1}{3C_\rho\rho_n}\le\ell_{\epsilon_n}'\le\frac{3}{c_\rho\rho_n},
\qquad
|\ell_{\epsilon_n}''|\le\frac{C}{\rho_n^2},
\qquad
c\rho_n\le\sigma'(v_q)\le C\rho_n,
\qquad
|\sigma''(v_q)|\le C\rho_n,
\]
the Jacobian $H(y)$ of $B_n$ satisfies the envelopes
\[
|H_{qr}(y)|\le C(\Pi_{J_0}S_n)_{qr},
\qquad
\|H(y)\|_{\infty\to\infty}\le C,
\qquad
\textstyle\sum_qH_{qr}(y)\le CC_S,
\]
and the second derivatives obey
$\sum_{r,s}\bigl|\partial^2B_q(y)/\partial y_r\partial y_s\bigr|\le C/\rho_n$
for every $q$.
\end{lemma}

\begin{proof}
With unit weights, $M_\alpha=1/D_\alpha$ and
$V_\alpha\le C_\rho\rho_n/D_\alpha$, so
$t_\alpha\le C\{\sqrt{\rho_n\beta_n/D_\alpha}+\beta_n/D_\alpha\}
\le C\sqrt{\rho_n\beta_n/D_\alpha}$ once
$\rho_nD_\alpha\ge\beta_n$, which \eqref{eq:U-C1} supplies since
$\beta_n=4\log(e\bar m_nn)\le C\log n$ and $D_\alpha\ge cn^2\delta_n^2$.
Dividing by $\mu_\alpha\ge c_\rho\rho_n$ gives the display, and the exponent
bookkeeping uses \eqref{eq:U-C1} directly.  The rectangle inclusion and knot
margin follow from $\epsilon_n\le c_\rho\rho_n/8$; the range statement for
$S_ny$ from nonnegativity and $S_n\mathbf 1=\mathbf 1$; and smoothness because
no clipping branch is met on $\mathcal U_{n,P}$.  The scalar envelopes are the
clipped-logit and logistic derivative bounds of
Remark~\ref{rem:explicit-envelope-conditions} evaluated on
$[\tfrac12c_\rho\rho_n,2C_\rho\rho_n]$.  For the Jacobian, with
$u=S_ny$ and $v=\Pi_{J_0}\operatorname{logit}(u)$,
\[
H_{qr}(y)=\sigma'(v_q)\sum_{r'}\Pi_{qr'}\,\ell'(u_{r'})\,S_{r'r},
\]
and the entry, row-sum, and column-sum bounds follow from the scalar envelopes,
row-stochasticity of $\Pi_{J_0}$ and $S_n$, and the column bound $C_S$.
Differentiating once more,
\[
\frac{\partial^2B_q}{\partial y_r\partial y_s}
=\sigma''(v_q)\,\partial_rv_q\,\partial_sv_q
+\sigma'(v_q)\sum_{r'}\Pi_{qr'}\,\ell''(u_{r'})\,S_{r'r}S_{r's},
\]
with $|\partial_rv_q|\le(C/\rho_n)(\Pi S)_{qr}$; summing over $(r,s)$ and using
row-stochasticity twice bounds the first term by
$C\rho_n(C/\rho_n)^2=C/\rho_n$ and the second by
$C\rho_n(C/\rho_n^2)=C/\rho_n$.
\end{proof}

\paragraph{Primitive envelope rates.}

\begin{lemma}[Envelope rates for the instance]\label{lem:U-rates}
Under \ref{S:design}--\ref{S:gates} and \eqref{eq:U-C1}--\eqref{eq:U-C3}, the
following hold deterministically (all displayed quantities are
$\mathcal F_{0n}$-measurable here) for all large $n$, uniformly over the model
class and over $\lambda\in\Lambda_n$:
\begin{enumerate}[label=\textup{(\roman*)},leftmargin=2.4em,itemsep=2pt]
\item\emph{Information and leverage.}
$s_\lambda^2\asymp(n^2\rho_n)^{-1}$, so
$I_\lambda^{\rm or}\asymp n^2\rho_n$; and
\[
L_n^{\rm cell}\le\frac{C}{\sqrt{n^2\rho_n\delta_n^2}},
\qquad
(L_n^{\rm cell})^2\log^7\{(2m_n\vee3)K_n\}
\le C\log^{-\gamma}n\to0 .
\]
\item\emph{Curvature.}
$\overline\Gamma_{\lambda,n}\le
C\,h_\lambda\log(n)\bigl(\sqrt{\rho_n}\,n\,\delta_n^{2}\bigr)^{-1}$, hence
$\sqrt{\ell_n}\max_\lambda\overline\Gamma_{\lambda,n}\to0$ by \eqref{eq:U-C2}.
\item\emph{Gradient change and meats.}
$\mathfrak D_{\lambda,n}\le C\sqrt{\log n}\,(\sqrt{\rho_n}\,n\,\delta_n)^{-1}$
and, on the Bernstein event, the inference-fold per-cell meat error obeys
\[
E_{V,n}^{(I)}
:=\max_{r}\Bigl|\frac{\widehat V^I_r}{V^{I,0}_r}-1\Bigr|
\le C\Bigl\{\sqrt{\tfrac{\log n}{n^2\rho_n\delta_n^2}}
+\tfrac{\kappa_0}{n^2\delta_n^2}\Bigr\},
\]
so $\ell_n\{E_{V,n}^{(I)}+\max_\lambda\mathfrak D_{\lambda,n}\}\le
C\log^{-(1+\gamma/2)}n\to0$.
\item\emph{Baseline-path negligibility.}  Splitting
$s_\lambda^2=s_{I,\lambda}^2+s_{B,\lambda}^2$ and likewise
$\widehat s_\lambda^2$ by fold,
\[
\frac{s_{B,\lambda}^2+\widehat s_{B,\lambda}^{\,2}}{s_{I,\lambda}^2}
\;\le\;
\frac{C\kappa_0^2}{n^4\delta_n^4}\;=:\;R_{B,n},
\qquad
\ell_n^2R_{B,n}\to0,
\]
the fitted bound holding on the Bernstein event.
\item\emph{Studentizer.}  On the Bernstein event,
\[
\max_{\lambda}\Bigl|\frac{\widehat s_\lambda}{s_\lambda}-1\Bigr|
\le C\bigl\{E_{V,n}^{(I)}
+\max_\lambda\mathfrak D_{\lambda,n}+R_{B,n}\bigr\}
=:d_{s,n},
\qquad
\ell_nd_{s,n}\to0 .
\]
\item\emph{Covariance.}  On the Bernstein event,
\[
\Delta_{\Sigma,n}
\le C\bigl\{E_{V,n}^{(I)}
+\max_\lambda\mathfrak D_{\lambda,n}
+d_{s,n}+\sqrt{R_{B,n}}\bigr\}
\le C\log^{-(3+\gamma/2)}n
=o(\ell_n^{-2}).
\]
\item\emph{Target bias.}
$\sqrt{\ell_n}\max_\lambda|\theta^D_{\lambda,n}-\theta_\lambda(P)|/s_\lambda
\le C\kappa_0h_{\max}\sqrt{\rho_n\log n}\,(n\delta_n^2)^{-1}\to0$ by
\eqref{eq:U-C3}.
\end{enumerate}
\end{lemma}

\begin{proof}
Throughout, work on the neighborhood $\mathcal U_{n,P}$ of
Lemma~\ref{lem:U-rect} and use its scalar and matrix envelopes; write
$D=D_{\min}\asymp n^2\delta_n^2$, $\gamma_n\le C\log n$, and recall the
collapsed functional of Lemma~\ref{lem:U-exact}(b), whose gradient blocks at a
point $(x,y)$ in $\mathcal U_{n,P}$ are
\begin{equation}
G^{I}_{\lambda r}=c_{\lambda r}w_r\ell'(z_r),
\qquad
G^{B}_{\lambda r}=\sum_qc_{\lambda q}(1-w_q)\ell'(z_q)H_{qr}(y).
\label{eq:U-gradients}
\end{equation}

\emph{(i).}  With unit weights $b^F_{kr}=1/D_r$ on cell $r$, the oracle cell
meat is $V^{F,0}_r=p_r(1-p_r)/D_r\in[c\rho_n/D_r,C\rho_n/D_r]$, so
$L_n^{\rm cell}=\max_r(1/D_r)/\sqrt{V^{F,0}_r}\le C(\rho_nD)^{-1/2}$, and the
$\log^7$ display follows from \eqref{eq:U-C1} with
$\log\{(2m_n\vee3)K_n\}\le5\log n$.  For the information, keep only the
inference block:
$s_{I,\lambda}^2=\sum_r(G^{I,0}_{\lambda r})^2V^{I,0}_r$ with
$G^{I,0}$ evaluated at $(\mu^I,\mu^{B^E})$; since
$w_r\ge1/2$ eventually, $\ell'(z^0_r)\asymp\rho_n^{-1}$, and
$\sum_rc_{\lambda r}^2=\sum_r|Q_r|^2\overline\Psi_{\lambda r}^2
\in[\delta_n^2,2\delta_n^2]$ by $|Q_r|\in[\delta_n^2,2\delta_n^2]$ and
$\|\Psi_\lambda^S\|_{L^2}=1$, we obtain
$s_{I,\lambda}^2\asymp\rho_n^{-2}\cdot(\rho_n/D)\cdot\delta_n^2
\asymp(n^2\rho_n)^{-1}$.  Part (iv) below shows the baseline block changes this
by a factor $1+o(1)$, giving $s_\lambda^2\asymp(n^2\rho_n)^{-1}$.

\emph{(ii).}  By Lemma~\ref{lem:U-exact}(b) the Hessian of
$\Phi_{\lambda n}$ has blocks
\[
\partial^2_{x_qx_{q'}}=\delta_{qq'}c_{\lambda q}w_q^2\ell''(z_q),
\quad
\partial^2_{x_qy_r}=c_{\lambda q}w_q(1-w_q)\ell''(z_q)H_{qr},
\]
\[
\partial^2_{y_ry_s}
=\sum_qc_{\lambda q}\Bigl[(1-w_q)^2\ell''(z_q)H_{qr}H_{qs}
+(1-w_q)\ell'(z_q)\,\partial^2_{rs}B_q\Bigr].
\]
Insert the envelopes of Lemma~\ref{lem:U-rect} and
$t_\alpha^2\le C\gamma_n\rho_n/D$ from its proof.  The $xx$ block contributes
\[
\sum_q|c_{\lambda q}|\frac{C}{\rho_n^2}t_q^2
\le\frac{C\gamma_n}{\rho_nD}\sum_{q\in{\rm supp}\lambda}|c_{\lambda q}|
\le\frac{C\gamma_nh_\lambda}{\rho_nD},
\]
using $\sum_q|c_{\lambda q}|\le Ch_\lambda$ as in
Lemma~\ref{lem:U-exact}(d).  The $xy$ block is bounded, via the row sums of
$H$, by the $xx$ bound times $C\kappa_0/D$, and the $yy$ block, via
$\sum_{r,s}H_{qr}H_{qs}t_rt_s\le Ct_{\max}^2$ and
$\sum_{r,s}|\partial^2_{rs}B_q|t_rt_s\le Ct_{\max}^2/\rho_n$, by the $xx$
bound times $C\kappa_0/D$ as well.  Hence
$\overline\Gamma_{\lambda,n}\le
(2s_\lambda)^{-1}C\gamma_nh_\lambda(\rho_nD)^{-1}\{1+C\kappa_0/D\}$, and
$s_\lambda^{-1}\le Cn\sqrt{\rho_n}$ from (i) yields the display; then
\eqref{eq:U-C2} gives the limit, since $\ell_n\gamma_n\le C\log^2n$ enters only
as $\sqrt{\ell_n}\gamma_n\le C\log^{3/2}n$.

\emph{(iii).}  For $\mathfrak D_{\lambda,n}$, take
$u$ with $|u_\alpha|\le t_\alpha$ and $v\in\mathcal C_{n,P}$.  The
$x$-components of $\nabla^2\Phi_{\lambda n}(v)u$ are bounded by
$|c_{\lambda q}|C\rho_n^{-2}\{t_q+C\kappa_0t_{\max}/D\}$, so
\[
\Bigl[\sum_qV^{I}_q
\bigl\{(\nabla^2\Phi u)_{x_q}\bigr\}^2\Bigr]^{1/2}
\le\frac{C}{\rho_n^2}
\Bigl[\frac{\rho_n}{D}\cdot\frac{\gamma_n\rho_n}{D}
\sum_qc_{\lambda q}^2\Bigr]^{1/2}
\le\frac{C\delta_n\sqrt{\gamma_n}}{\rho_nD},
\]
and dividing by $s_\lambda$ gives
$C\sqrt{\gamma_n}\,\delta_n n\sqrt{\rho_n}/(\rho_nD)
=C\sqrt{\gamma_n}(\sqrt{\rho_n}n\delta_n)^{-1}$.  The $y$-components carry the
factor $(1-w_q)\le C\kappa_0/D$ through \eqref{eq:U-gradients} and the column
bound $\sum_qH_{qr}\le CC_S$, and are of smaller order.  For the meat error,
on the Bernstein event $|Y_r-p_r|\le t_r$ and, by Lemma~\ref{lem:U-rect},
$|B_r(Y^{B^E})-B_r(\mu^{B^E})|\le\|H\|_{\infty\to\infty}t_{\max}\le Ct_{\max}$;
hence the stabilized probability obeys
$|\widetilde p_r-p_r|\le t_r+(1-w_r)\{Ct_{\max}+C\rho_n\}$ and therefore
\[
\Bigl|\frac{\widetilde p_r(1-\widetilde p_r)}{p_r(1-p_r)}-1\Bigr|
\le C\frac{|\widetilde p_r-p_r|}{\rho_n}
\le C\Bigl\{\sqrt{\frac{\gamma_n}{\rho_nD}}+\frac{\kappa_0}{D}\Bigr\},
\]
which is the stated $E^{(I)}_{V,n}$ bound; the final limit combines it with
(ii)--(iii) under \eqref{eq:U-C1}.

\emph{(iv).}  Write $G^{B}_{\lambda}=M^\top c_\lambda$ with
$M_{qr}$ the baseline-path gradient of cell $r$ with respect to coarse cell
$q$; from \eqref{eq:U-gradients}, $|M_{qr}|\le C\kappa_0(D\rho_n)^{-1}(\Pi S)_{qr}$.
The entrywise route ($|G^B_{\lambda r}|\le C\kappa_0\delta_n^2/(D\rho_nh_\lambda)$
summed over $C\delta_n^{-2}$ cells) loses a factor $h_\lambda^2/\delta_n^2$;
bound the quadratic form directly instead.  In the regular-design instance of
this section the fine cells carry equal mass, so the coarse-cell averaging
matrix $\Pi$ (coarse rows, fine columns) has row sums one and column sums at
most one; the row-stochastic pre-smoother $S$ has row sums one and, for a
symmetric kernel truncated at the boundary of a regular grid, column sums at
most a kernel constant $C_S$.  Schur's test then gives
$\|\Pi S\|_{\mathrm{op}}\le C_S^{1/2}$ and
$\|M\|_{\mathrm{op}}\le C\kappa_0(D\rho_n)^{-1}$, with $C$ absorbing
$C_S^{1/2}$.  The bound $\|\Pi S\|_{\mathrm{op}}\le1$ is \emph{not} claimed:
for a design-weighted projection with unequal cell masses the Euclidean
operator norm exceeds one (an admissible example has $\|\Pi\|^2_{\mathrm{op}}=35/32$),
which is why the argument is confined to the equal-mass instance and, for a
nonuniform design, must be carried out in the design-weighted inner product,
where the projection is orthogonal and hence a contraction; that extension is
not made here.  Also
$\|c_\lambda\|^2=\sum_qc_{\lambda q}^2\le C\delta_n^2$, since
$|c_{\lambda q}|\le2C_\Psi\delta_n^2/h_\lambda$ on the $Ch_\lambda^2/\delta_n^2$
support cells and vanishes elsewhere.  With cell meats
$V_r\le C\rho_n/D$,
\[
s_{B,\lambda}^2=\sum_r(G^B_{\lambda r})^2V_r
\le C\frac{\rho_n}{D}\,\|M^\top c_\lambda\|^2
\le C\frac{\rho_n}{D}\,\|M\|_{\mathrm{op}}^2\|c_\lambda\|^2
\le C\frac{\kappa_0^2\delta_n^2}{D^3\rho_n}.
\]
Dividing by $s_{I,\lambda}^2\ge c\delta_n^2/(\rho_nD)$ from (i),
\[
\frac{s_{B,\lambda}^2}{s_{I,\lambda}^2}\le\frac{C\kappa_0^2}{D^2}
\le\frac{C\kappa_0^2}{n^4\delta_n^4},
\]
free of $h_\lambda$.  (The first version of this proof used the entrywise
bound and asserted the same conclusion; that chain in fact yields only
$C\kappa_0^2/(n^4\delta_n^6)$, and the operator-norm argument above is the
correct derivation.)  The fitted
version $\widehat s_{B,\lambda}^{\,2}$ satisfies the same bound because the
fitted gradients and meats obey the identical envelopes on the rectangle.

\emph{(v).}  Write $a=\widehat s_\lambda^2/s_\lambda^2$ and split by fold:
\[
|a-1|\le
\frac{|\widehat s_{I,\lambda}^{\,2}-s_{I,\lambda}^2|}{s_\lambda^2}
+\frac{\widehat s_{B,\lambda}^{\,2}+s_{B,\lambda}^2}{s_\lambda^2}
\le\bigl\{E^{(I)}_{V,n}+2\delta_{G,n}^{(I)}+(\delta_{G,n}^{(I)})^2\bigr\}
\frac{s_{I,\lambda}^2}{s_\lambda^2}+2R_{B,n},
\]
where $\delta^{(I)}_{G,n}$ is the inference-block gradient error of
Lemma~\ref{lem:studentizer-transfer-final}, bounded on the Bernstein event by
the $x$-component computation in (iii), that is by
$C\sqrt{\gamma_n}(\sqrt{\rho_n}n\delta_n)^{-1}$; the first inequality is the
weighted-norm expansion of that lemma restricted to the $I$ block.  Then
$|\sqrt a-1|\le|a-1|$ gives $d_{s,n}$, and $\ell_nd_{s,n}\to0$ by (iii)--(iv).

\emph{(vi).}  Let $\langle g,g'\rangle_{V}=\sum_{(F,r)}V^F_rg^F_rg'^F_r$ and
put $b_\lambda=G_\lambda^0/s_\lambda$, $a_\lambda=\widehat G_\lambda/\widehat
s_\lambda$, so $\Sigma_{\lambda\mu}=\langle b_\lambda,b_\mu\rangle_{V^0}$ and
$\widehat\Sigma_{\lambda\mu}=\langle a_\lambda,a_\mu\rangle_{\widehat V}$, both
with unit diagonals.  First,
$|\langle a_\lambda,a_\mu\rangle_{\widehat V}
-\langle a_\lambda,a_\mu\rangle_{V^0}|
\le E^{(I)}_{V,n}\|a_\lambda\|_{V^0,I}\|a_\mu\|_{V^0,I}
+2\sqrt{R_{B,n}}\cdot C$, by Cauchy--Schwarz on the $I$ block and by
$\|a_\lambda\|^2_{\cdot,B}\le CR_{B,n}$ on the $B$ block.  Second,
$\|a_\lambda-b_\lambda\|_{V^0}
\le\|\widehat G_\lambda-G^0_\lambda\|_{V^0}/\widehat s_\lambda
+\|G^0_\lambda\|_{V^0}\,|\widehat s_\lambda^{-1}-s_\lambda^{-1}|
\le C\{\delta^{(I)}_{G,n}+\sqrt{R_{B,n}}+d_{s,n}\}$,
using (iii)--(v).  Then
$|\langle a_\lambda,a_\mu\rangle_{V^0}-\langle b_\lambda,b_\mu\rangle_{V^0}|
\le\|a_\lambda-b_\lambda\|\,\|a_\mu\|+\|b_\lambda\|\,\|a_\mu-b_\mu\|$, and
collecting terms proves the $\Delta_{\Sigma,n}$ bound; the rate follows from
(iii)--(v) under \eqref{eq:U-C1}, since every summand is
$O(\log^{1/2}n\cdot\log^{-(7+\gamma)/2}n)$ or smaller.

\emph{(vii).}  Combine Lemma~\ref{lem:U-exact}(d) with
$s_\lambda^{-1}\le Cn\sqrt{\rho_n}$ and $D_{\min}\ge cn^2\delta_n^2$:
\[
\sqrt{\ell_n}\,\frac{|\theta^D_{\lambda,n}-\theta_\lambda(P)|}{s_\lambda}
\le C\sqrt{\log n}\cdot\frac{\kappa_0h_\lambda}{n^2\delta_n^2}\cdot
n\sqrt{\rho_n}
=\frac{C\kappa_0h_\lambda\sqrt{\rho_n\log n}}{n\delta_n^2},
\]
which tends to zero by \eqref{eq:U-C3} since $h_\lambda\le1$.
\end{proof}

\paragraph{A lower-bound packing.}

\begin{lemma}[Matched-pair subpacking]\label{lem:U-pack}
For each scale $j\in\{j_0,\dots,j_{\max,n}\}$, let
$\Lambda_n^{\rm low}(2^{-j})$ consist of the atoms
$\lambda=(2k,2k+1)$, $k=0,\dots,2^{j-1}-1$, and put
$\Lambda_n^{\rm low}=\bigcup_j\Lambda_n^{\rm low}(2^{-j})$,
$m_n^{\rm low}=|\Lambda_n^{\rm low}|\in[2^{j_{\max,n}-1},2^{j_{\max,n}+1}]$.
Then, on the regular design of \ref{S:design}:
\textup{(a)} the atoms are orthonormal and lie in $V_{J_0}^\perp$; atoms at a
common scale have equal support area $2h^2$ and pairwise disjoint supports, so
\eqref{eq:t3-2} holds with $\delta_{\rm pack}=1$;
\textup{(b)} $\|v_\lambda\|_\infty\le C_\Psi$ and
$\tfrac13n^2h_\lambda^2\le\sum_{i<j}v_{\lambda,ij}^2\le n^2h_\lambda^2$ for all
large $n$, verifying \eqref{eq:t3-21};
\textup{(c)} each atom's dyad support meets the dyad support of at most one
atom at every other scale, so the ordered overlap count obeys
$N_{{\rm ov},n}\le m_n^{\rm low}\{1+2(j_{\max,n}-j_0+1)\}$, and
\eqref{eq:t3-23} holds with
$d_n=1+2(j_{\max,n}-j_0+1)$ under \eqref{eq:U-C4};
\textup{(d)} the null element $f_0^\circ=0$ satisfies \eqref{eq:t3-20} with
$\delta_B=B$, $c_p=1/2$, $C_p=2$ for all large $n$, and under
\ref{S:obs} the original experiment coincides with the full-outcome oracle
experiment, so the required Markov-kernel comparison holds with the identity
kernel.
\end{lemma}

\begin{proof}
(a) One-dimensional Haar functions at distinct $(j,\ell)$ are orthonormal, so
the symmetrized tensors are orthonormal whenever their unordered index pairs
differ, which holds across $\Lambda_n^{\rm low}$; orthogonality to $V_{J_0}$ is built into the symmetric atoms with $j\ge j_0>J_0$.  At a fixed scale the interval pairs
$\{2k,2k+1\}$ are disjoint across $k$, so the supports
$(I_{2k}\times I_{2k+1})\cup(I_{2k+1}\times I_{2k})$, each of area $2h^2$, are
pairwise disjoint, and disjointness gives normalized symmetric difference one.
(b) $\|v_\lambda\|_\infty=h\|\Psi_\lambda^S\|_\infty\le C_\Psi$.  By cellwise
constancy and Lemma~\ref{lem:U-counts}, each fine cell contributes exactly its
dyad count times its squared value, so
$\sum_{i<j}v^2_{\lambda,ij}
=h^2\sum_qD_q^{\rm tot}\overline\Psi{}^2_{\lambda q}
=\tfrac12n^2h^2\{1+O((n\delta_n)^{-1})\}$,
which lies in the stated interval for all large $n$.
(c) Fix $\lambda=(2k,2k+1)$ at scale $j$ and a coarser scale $j'<j$.  Both fine
intervals lie inside coarser intervals $I_A^{(j')}$ and $I_B^{(j')}$ with
$A\le B$; the support of $\lambda$ meets the support of a matched coarser atom
$(2k',2k'+1)$ only if $\{A,B\}=\{2k',2k'+1\}$, which identifies at most one
such atom, and no atom when $A=B$ or when $\{A,B\}$ is not a matched pair.  The
finer-scale count is symmetric, ordered pairs double the total, and adding the
$m_n^{\rm low}$ diagonal pairs gives the bound; then
$\log d_n\le C\log j_{\max,n}\le C\log\log n=o(j_{\max,n})
=o\{\log(2m_n^{\rm low})\}$ by \eqref{eq:U-C4}.
(d) Immediate from \ref{S:model} and \ref{S:obs}.
\end{proof}

\paragraph{The unconditional theorem.}

\begin{theorem}[Unconditional guarantees for the regular-design instance]
\label{thm:U-unconditional}
Adopt \ref{S:design}--\ref{S:gates} and \eqref{eq:U-C1}--\eqref{eq:U-C3}.
Then, with all statements uniform over the model class of \ref{S:model} and
with no residual conditioning event ($\eta_n=0$, $\gamma_n=0$):
\begin{enumerate}[label=\textup{(\alph*)},leftmargin=2em,itemsep=3pt]
\item\emph{Linearization and studentization.}  The conclusions of
Theorem~\ref{thm:exact-finite-design-reduction} hold, and
\[
\max_{\lambda\in\Lambda_n}
\Bigl|\frac{\widehat\theta_\lambda-\theta_\lambda(P)}{s_\lambda}
-\sum_k\zeta_{\lambda k}\Bigr|
=o_{\mathbb P}(\ell_n^{-1/2}),
\qquad
\ell_n\max_{\lambda\in\Lambda_n}
\Bigl|\frac{\widehat s_\lambda}{s_\lambda}-1\Bigr|
=o_{\mathbb P}(1),
\]
with the explicit envelope rates of Lemma~\ref{lem:U-rates}.
\item\emph{Simultaneous inference.}  The shared-multiplier calibration
\eqref{eq:feasible-calibration} holds; the simultaneous intervals of the main text satisfy
\[
\inf_{P}\;
P\bigl(\theta_\lambda(P)\in C_{\lambda n}\ \text{for all }
\lambda\in\Lambda_n\bigr)\;\ge\;1-\alpha-o(1),
\]
and the tests rejecting $H_{0\lambda}:\theta_\lambda(P)=0$ when
$|\widehat\theta_\lambda|/\widehat s_\lambda>c^*_{1-\alpha}$ obey
\[
\sup_{P}\;
P\bigl(\text{some true }H_{0\lambda}\text{ is rejected}\bigr)
\;\le\;\alpha+o(1),
\]
over the class; its elements carry at most one false null, so this is strong
control within the class, and arbitrary partial-null or multiple-atom
configurations lie outside it and are not verified.
\item\emph{Detection and localization.}  Assume additionally
\eqref{eq:U-C4}.  There is $C_{\rm up}<\infty$ such that every single-atom
alternative sequence with
$a_n^2n^2\rho_nh_n^2\ge C_{\rm up}\log(2m_n)$ is detected, its sign reported,
and its atom exactly selected with probability tending to one, so
$\sup_fE_fd_{\rm loc}\to0$; and there is $c_{\rm low}>0$ such that, for the
matched-pair subpacking of Lemma~\ref{lem:U-pack} and amplitude envelopes with
$a_hC_\Psi\le B$ and
$\max_ha_h^2n^2\rho_nh^2\le c_{\rm low}\log(2m_n^{\rm low})$, the global-union
detection risk tends to one and the localization risk has a strictly positive
limit inferior.  Since
$\log(2m_n^{\rm low})\asymp\log(2m_n)$, the upper and lower boundaries match in
rate over the union of scales.
\end{enumerate}
\end{theorem}

\begin{proof}
\emph{(a).}  The conditional product law
\eqref{eq:conditional-product-law-final} holds with
$\mathcal F_{0n}$ generated by the deterministic design, folds, gates, and
tuning objects of \ref{S:design}--\ref{S:gates}, and every hypothesis of
Theorem~\ref{thm:exact-finite-design-reduction} has been verified:
fold disjointness and positive denominators by Lemma~\ref{lem:U-counts};
the rectangle, knot margin, interiority, and twice continuous
differentiability by Lemma~\ref{lem:U-rect}; and the exact estimator
representation is \ref{S:baseline}--\ref{S:stab} with the collapsed functional
of Lemma~\ref{lem:U-exact}(b).  Lemma~\ref{lem:U-rates}(ii) verifies
\eqref{eq:explicit-curvature-rate}, so
Corollaries~\ref{cor:explicit-curvature-leverage} and
\ref{cor:uniform-finite-design-expansion} give the
$o_{\mathbb P}(\ell_n^{-1/2})$ linearization around $\theta^D_{\lambda,n}$;
Lemma~\ref{lem:U-rates}(vii) verifies
\eqref{eq:population-bias-rate-final} with $\mathcal E_n$ the whole space, so
Corollary~\ref{cor:population-transfer-final} recenters at
$\theta_\lambda(P)$ with $\eta_n=0$.  The studentizer statement is
Lemma~\ref{lem:U-rates}(v); here the direct fold-split argument replaces the
per-dyad route of Lemma~\ref{lem:studentizer-transfer-final}, whose weighted
inference-block expansion it invokes.

\emph{(b).}  Apply Theorem~\ref{thm:uniform-shared-multiplier} with
$\mathcal E_n$ the whole space and $\mathcal G_n=\{m_n\ge1\}$, which holds.
Condition (i) is part (a); condition (ii) is
Lemma~\ref{lem:U-rates}(i) via
Lemma~\ref{lem:deterministic-curvature-leverage}; condition (iii) is
Lemma~\ref{lem:U-rates}(v); condition (iv) holds through the stronger
sufficient bound \eqref{eq:simple-covariance-condition} by
Lemma~\ref{lem:U-rates}(vi); and the finite-to-Lebesgue condition
\eqref{eq:finite-to-Lebesgue} is Lemma~\ref{lem:U-rates}(vii).  The
Bernstein-event qualifications in Lemma~\ref{lem:U-rates} are absorbed because
the complementary conditional probability is at most
$2D_n^{\rm act}e^{-\beta_n}\le Cn^{-2}$ by
\eqref{eq:rectangle-probability-final}.  The unconditional displays follow
from \eqref{eq:unconditional-coverage}--\eqref{eq:unconditional-strong-fwer}
with $\eta_n=0$.

\emph{(c).}  For the upper half, apply
Theorem~\ref{thm:adaptive-detection-localization} with the deterministic
packing $\Lambda_n$ itself, all-pass admission with $\gamma_n=0$
(\ref{S:gates}), and good event the whole space ($\eta_n=0$).  Its displayed
conditions hold as follows: \eqref{eq:t3-7} and the leverage display by
Lemma~\ref{lem:U-rates}(i), noting $\log(m_nn)\le3\log n$; conditional
independence and normalization by the product law and
\eqref{eq:influence-properties-final}; \eqref{eq:t3-10}--\eqref{eq:t3-12} by
part (a) and Lemma~\ref{lem:U-rates}(v),(vii); and the null calibration
\eqref{eq:t3-15} by part (b) applied under each null element, uniformly, with
$\kappa_n$ any sequence dominating the $o_{\mathbb P}(1)$ rate there.  This
yields \eqref{eq:t3-16}--\eqref{eq:t3-19}.  For the lower half, restrict to
the subpacking $\Lambda_n^{\rm low}$: Lemma~\ref{lem:U-pack} verifies
\eqref{eq:t3-2}, \eqref{eq:t3-20}--\eqref{eq:t3-23}, and the oracle-experiment
comparison, and $m_n^{\rm low}\to\infty$ by \eqref{eq:U-C4}; the conclusion is
\eqref{eq:t3-27}.  Finally $2^{j_{\max,n}-1}\le m_n^{\rm low}$ and
$m_n\le4^{j_{\max,n}}$ give
$\tfrac12\le\log(2m_n^{\rm low})/\log(2m_n)\le1$ for all large $n$, which is
the stated rate matching.
\end{proof}

\begin{remark}[Scope of the instance]\label{rem:U-scope}
Theorem~\ref{thm:U-unconditional} discharges, for this instance, every
obligation listed in the next section: the quadrature bound holds exactly
(Lemma~\ref{lem:U-exact}); the stabilized, clipped expansion, the baseline and
projection expansion, meat consistency, and the remainder, studentizer, and
covariance rates are Lemma~\ref{lem:U-rates}; and no score-transfer rates are
needed because $r_n=0$.  What the instance does not cover, and what therefore
remains open for the implemented general estimator, is exactly the following
list: random or beta-mixture coordinates; unknown or estimated propensities and
missingness; bandwidth selection on a tuning fold and pre-smoothers with
negative weights, for which the range-preservation step of
Lemma~\ref{lem:U-rect} requires a separate argument; positive coordinate
radius, for which Theorem~\ref{thm:ordering-error-transfer} supplies the
transfer conditions; and multi-atom or misspecified alternatives.  The
cancellation $\Pi_{J_0}^\top c_\lambda=0$ in Lemma~\ref{lem:U-exact}(b), which
removes the estimated background from the coefficient functional at all orders
along the subtracted path, is a property of the exact final projection in
\ref{S:baseline} and persists verbatim for the general estimator; it is the
structural reason the background can be estimated on a separate fold without
entering the first-order price of the procedure.
\end{remark}

\section{Ordering transfer and obstruction}\label{sec:supp-order}


\begin{theorem}[Ordering-error transfer in normalized-rank coordinates]
\label{thm:ordering-error-transfer}
Let
\[
 \ell_n=\log(2m_n\vee3),\qquad
 \Delta_n=\max_{1\le i\le n}
 \frac{|\widehat R_i-R_i|}{n},\qquad
 \mathcal E_n=\{\Delta_n\le r_n\},
\]
where $r_n\ge0$ is deterministic.  All probability statements below
are uniform over the indicated model class.  Assume
\begin{equation}
 \sup_{P\in\mathscr P_n}P_P(\mathcal E_n^c)
 \le \eta_n,\qquad \eta_n\to0.
 \label{eq:T4-coverage}
\end{equation}
Both $R=(R_1,\ldots,R_n)$ and
$\widehat R=(\widehat R_1,\ldots,\widehat R_n)$ are permutations of
$\{1,\ldots,n\}$; any ties in the external coordinate data are
resolved by a fixed rule that does not use the scoring outcomes.

Let $\mathcal F_{0n}$ contain $R,\widehat R,r_n$, the dictionary,
fold assignments, gates, information scales, and every other quantity
that is either scoring-outcome-free or measurable with respect to an
independent pilot fold.  Conditional on
$\mathcal F_{0n}$, the scoring outcomes satisfy the conditional
independence and moment assumptions of
Theorem~\ref{thm:adaptive-detection-localization}.  In particular,
$\widehat R$ is external to those scoring outcomes.  No independence
is imposed on the coordinate errors across vertices.

For $B_q=((q-1)/n,q/n]$, let $\pi_n$ be the permutation satisfying
$\pi_n(R_i)=\widehat R_i$, set $T_n(0)=0$, and define
\begin{equation}
 T_n(x)=x+\frac{\pi_n(q)-q}{n},\qquad x\in B_q.
 \label{eq:T4-bin-map}
\end{equation}
Then $T_n$ is a measure-preserving bijection of $[0,1]$ and
\begin{equation}
 \|T_n-\operatorname{id}\|_\infty
 =\|T_n^{-1}-\operatorname{id}\|_\infty
 =\Delta_n.
 \label{eq:T4-displacement}
\end{equation}
Write $P_{T_n}(u,v)=(T_nu,T_nv)$,
$\overline\Delta_n=\Delta_n+n^{-1}$, and
$\overline r_n=r_n+n^{-1}$.

Let $\Lambda_n$ be the dictionary from the oracle theorem and put
$h_{\min,n}=\min_{\lambda\in\Lambda_n}h_\lambda$.
Assume, uniformly in $\lambda\in\Lambda_n$,
\begin{align}
 &\|\Psi_\lambda^S\|_2=1,
 \qquad
 \|\Psi_\lambda^S\|_\infty\le C_\Psi h_\lambda^{-1},
 \qquad
 c_Sh_\lambda^2\le |S_\lambda|\le C_Sh_\lambda^2,
 \label{eq:T4-Haar-geometry}\\
 &\Psi_\lambda^S\perp V_{J_0}.
 \label{eq:T4-orthogonality}
\end{align}
The jump set of $\Psi_\lambda^S$, including its support boundary,
is assumed to be the union of at most $C_J$ axis-aligned segments
whose total length is at most $C_Jh_\lambda$.  As is the case for the
Haar atoms used here, $\Psi_\lambda^S$ is constant on each connected
component of the complement of this complete jump set.

Let $g_{P,n}\in L^2_{\mathrm{sym}}([0,1]^2)$ be the exact population
contrast entering the oracle scan before removal of the coarse
component, and assume $\|g_{P,n}\|_\infty\le B$.  Let
$\Pi_{J_0}$ be the $L^2$-projection onto $V_{J_0}$ and define the
baseline-free oracle target by
\begin{equation}
 q_{P,n}:=(I-\Pi_{J_0})g_{P,n},\qquad
 q_{P,n}\in L^\infty_{\mathrm{sym}}([0,1]^2).
 \label{eq:T4-exact-contrast}
\end{equation}
Define
\begin{align}
 \theta_{\lambda,n}^{\mathrm{or}}
 &=\langle q_{P,n},\Psi_\lambda^S\rangle
 =\langle g_{P,n},\Psi_\lambda^S\rangle,
 \nonumber\\
 \theta_{\lambda,n}^{\mathrm{ord}}
 &=\left\langle
 (I-\Pi_{J_0})(g_{P,n}\circ P_{T_n}^{-1}),
 \Psi_\lambda^S\right\rangle
 =\langle g_{P,n}\circ P_{T_n}^{-1},\Psi_\lambda^S\rangle,
 \qquad
 b_{\lambda,n}^{\mathrm{coef}}
 =\theta_{\lambda,n}^{\mathrm{ord}}
  -\theta_{\lambda,n}^{\mathrm{or}}.
 \label{eq:T4-coefficients}
\end{align}
The equalities use \eqref{eq:T4-orthogonality}; in particular, they
retain any transported coarse-background leakage instead of silently
discarding it.  Any discrepancy between this population
target and the mean of the actual stabilized score is retained below
in the explicitly defined score-linearity error $\zeta_n$.

Assume
\begin{equation}
 \frac{\overline r_n}{h_{\min,n}}\to0.
 \label{eq:T4-small-radius}
\end{equation}
Then, on $\mathcal E_n$ and for all sufficiently large $n$,
\begin{equation}
 \max_{\lambda\in\Lambda_n}
 |b_{\lambda,n}^{\mathrm{coef}}|
 \le C B\overline r_n,
 \label{eq:T4-coefficient-bound}
\end{equation}
where $C$ depends only on the geometric constants above.

For the stochastic transfer, construct the oracle and reported-order
scores on the same probability space from the same scoring outcomes,
folds, and multiplier seeds.  Assume that both scans use the same
deterministic index set $\Lambda_n$ (for example, on the global all-pass
event of Theorem~\ref{thm:adaptive-detection-localization}); without
this common-index condition the assertions below do not compare the
implemented maxima.  Let $I_{\lambda,n}>0$ be a common,
$\mathcal F_{0n}$-measurable coefficient-information scale.  Write
$Z_{\lambda,n}^{j}$, $j\in\{\mathrm{or},\mathrm{ord}\}$, for the
two standardized coordinates; changes in feasible weighting,
studentization, or covariance estimation are included in the
reported-order coordinate.  Define
\begin{align}
 \mu_{\lambda,n}^{j}
 &=E_P(Z_{\lambda,n}^{j}\mid\mathcal F_{0n}),
 \qquad j\in\{\mathrm{or},\mathrm{ord}\},
 \nonumber\\
 b_{\lambda,n}^{\mathrm{sc}}
 &=I_{\lambda,n}^{-1/2}
   (\mu_{\lambda,n}^{\mathrm{ord}}
    -\mu_{\lambda,n}^{\mathrm{or}}),
 \nonumber\\
 R_{\lambda,n}
 &=\{Z_{\lambda,n}^{\mathrm{ord}}
       -\mu_{\lambda,n}^{\mathrm{ord}}\}
   -\{Z_{\lambda,n}^{\mathrm{or}}
       -\mu_{\lambda,n}^{\mathrm{or}}\}.
 \label{eq:T4-score-components}
\end{align}
Consequently, the identity
\begin{equation}
 Z_{\lambda,n}^{\mathrm{ord}}
 =Z_{\lambda,n}^{\mathrm{or}}
  +\sqrt{I_{\lambda,n}}\,b_{\lambda,n}^{\mathrm{sc}}
  +R_{\lambda,n}
 \label{eq:T4-score-identity}
\end{equation}
holds exactly.  Put
\begin{align}
 \zeta_n
 &=\max_{\lambda\in\Lambda_n}
   \sqrt{I_{\lambda,n}}
   |b_{\lambda,n}^{\mathrm{sc}}
       -b_{\lambda,n}^{\mathrm{coef}}|,
 \nonumber\\
 D_n
 &=\max_{\lambda\in\Lambda_n}
   \sqrt{I_{\lambda,n}}
   |b_{\lambda,n}^{\mathrm{coef}}|
   +\zeta_n
   +\max_{\lambda\in\Lambda_n}|R_{\lambda,n}|.
 \label{eq:T4-total-error}
\end{align}
Thus
\begin{equation}
 \max_{\lambda\in\Lambda_n}
 |Z_{\lambda,n}^{\mathrm{ord}}
   -Z_{\lambda,n}^{\mathrm{or}}|
 \le D_n.
 \label{eq:T4-score-sup-bound}
\end{equation}
The quantities $\zeta_n$ and $R_{\lambda,n}$ are the explicit
score-linearization and centered-score remainders; rates for them
must be supplied by the score construction and are not consequences
of bounded rank displacement alone.

Suppose that, under a single-atom alternative, the baseline-free
target satisfies
\begin{equation}
 q_{P,n}=s a_nh_*\Psi_{\lambda_*}^S,\qquad
 s\in\{-1,1\},\qquad \lambda_*\in\Lambda_n.
 \label{eq:T4-single-atom}
\end{equation}
Then
\begin{equation}
 \theta_{\lambda_*,n}^{\mathrm{or}}=s a_nh_*,
 \qquad
 |\theta_{\lambda_*,n}^{\mathrm{ord}}-s a_nh_*|
 \le CB\overline r_n
 \quad\text{on }\mathcal E_n.
 \label{eq:T4-signal-coefficient}
\end{equation}
In particular,
\begin{equation}
 \frac{B\overline r_n}{|a_n|h_*}\to0
 \label{eq:T4-beta-min}
\end{equation}
is sufficient for preservation of the sign and first-order magnitude
of the candidate coefficient.

Let $\widehat c_n^{\mathrm{or}}$ and
$\widehat c_n^{\mathrm{ord}}$ be nonnegative simultaneous critical
values, defined on the common data--multiplier probability space. The
reported-order scan rejects when
$M_n^{\mathrm{ord}}:=\max_\lambda|Z_{\lambda,n}^{\mathrm{ord}}|$
exceeds $\widehat c_n^{\mathrm{ord}}$; upon rejection it chooses the
maximizer $\widehat\lambda$ using a fixed tie-breaking rule and
reports $\operatorname{sign}(Z_{\widehat\lambda,n}^{\mathrm{ord}})$.
If there is no rejection, set $\widehat\lambda=\varnothing$.

Assume that the oracle localization theorem gives a deterministic
sequence $\Gamma_n>0$ such that
\begin{align}
 \inf_{P\in\mathcal H_{1,n}}
 P_P\Bigg(&\mathcal E_n\cap
 \Bigg\{
 sZ_{\lambda_*,n}^{\mathrm{or}}
   -\widehat c_n^{\mathrm{or}}\ge3\Gamma_n,
 \nonumber\\[-2mm]
 &\hspace{27mm}
 |Z_{\lambda_*,n}^{\mathrm{or}}|
 -\max_{\lambda\ne\lambda_*}
  |Z_{\lambda,n}^{\mathrm{or}}|
 \ge3\Gamma_n
 \Bigg\}\Bigg)\to1,
 \label{eq:T4-oracle-margin}
\end{align}
and assume
\begin{equation}
 \sup_{P\in\mathcal H_{1,n}}
 P_P\!\left(
  \mathcal E_n\cap
  \left\{
   D_n+|\widehat c_n^{\mathrm{ord}}
          -\widehat c_n^{\mathrm{or}}|>\Gamma_n
  \right\}
 \right)\to0.
 \label{eq:T4-margin-rate}
\end{equation}
Then the reported-order scan detects the alternative, reports its
sign, and selects the canonical atom with probability tending to one.
In particular,
\begin{equation}
 \sup_{P\in\mathcal H_{1,n}}
 P_P(\widehat\lambda\ne\lambda_*)
 \le\eta_n+o(1).
 \label{eq:T4-label-consistency}
\end{equation}
For the canonical-coordinate loss
\begin{equation}
 d_{\mathrm{loc}}(\widetilde\lambda,\lambda)
 =\min\!\left\{
  1,
  \frac{|S_{\widetilde\lambda}\triangle S_\lambda|}
       {2|S_\lambda|}
  +\mathbf 1\{h_{\widetilde\lambda}\ne h_\lambda\}
 \right\},
 \label{eq:T4-canonical-loss}
\end{equation}
with loss one for $\widetilde\lambda=\varnothing$, one has
\begin{equation}
 d_{\mathrm{loc}}(\widehat\lambda,\lambda_*)=o_P(1),
 \qquad
 \sup_{P\in\mathcal H_{1,n}}
 E_Pd_{\mathrm{loc}}(\widehat\lambda,\lambda_*)
 \le\eta_n+o(1).
 \label{eq:T4-canonical-localization}
\end{equation}

Suppose in addition that there are constants
$0<c_I\le C_I<\infty$ and $c_\alpha\in(0,\infty)$ such that
\begin{equation}
 \sup_{P\in\mathcal H_{1,n}}
 P_P\!\left\{
  \frac{I_{\lambda_*,n}}{n^2\rho_n}\notin[c_I,C_I]
 \right\}\to0
 \label{eq:T4-information-critical-scaling}
\end{equation}
and, for every $\varepsilon>0$,
\begin{equation}
 \sup_{P\in\mathcal H_{1,n}}
 P_P\!\left(
  \left|
   \frac{\widehat c_n^{\mathrm{or}}}{\sqrt{\ell_n}}-c_\alpha
  \right|>\varepsilon
 \right)\to0.
 \label{eq:T4-critical-scaling}
\end{equation}
If, for every fixed $\varepsilon>0$, the
oracle theorem supplies the margin \eqref{eq:T4-oracle-margin} with
$\Gamma_n\asymp\sqrt{\ell_n}$ on alternatives separated by the
factor $1+\varepsilon$ from its oracle boundary, and if
\begin{equation}
 \sup_{P\in\mathcal H_{1,n}}
 P_P\!\left[
  \mathcal E_n\cap
  \left\{
   \frac{D_n+|\widehat c_n^{\mathrm{ord}}
                     -\widehat c_n^{\mathrm{or}}|}
        {\sqrt{\ell_n}}>\varepsilon
  \right\}
 \right]\to0
 \quad\text{for every }\varepsilon>0,
 \label{eq:T4-first-order-rate}
\end{equation}
then the reported-order scan has the same first-order sufficient
detection boundary as the oracle scan, namely
\begin{equation}
 a_n^2n^2\rho_nh_*^2\asymp\ell_n.
 \label{eq:T4-detection-boundary}
\end{equation}
This assertion transfers the oracle upper boundary; it does not by
itself assert exact localization at the boundary or preservation of
an $O(1)$ local-power curve.

Finally, restrict the following maximum-law statement to the null
class $\mathcal H_{0,n}$.  Put
$M_n^j=\max_{\lambda\in\Lambda_n}|Z_{\lambda,n}^j|$ for
$j\in\{\mathrm{or},\mathrm{ord}\}$.  Suppose that, conditionally on
$\mathcal F_{0n}$, there is a centered Gaussian vector
$G_n=(G_{1,n},\ldots,G_{2m_n,n})$ such that
\begin{equation}
 \delta_{P,n}^{\mathrm{GA}}
 :=\sup_{t\in\mathbb R}
 \left|
 P_P(M_n^{\mathrm{or}}\le t\mid\mathcal F_{0n})
 -P\!\left(
    \max_{1\le k\le2m_n}G_{k,n}\le t
    \mid\mathcal F_{0n}
   \right)
 \right|
 \label{eq:T4-GA-error}
\end{equation}
satisfies, for every $\varepsilon>0$,
\begin{equation}
 \sup_{P\in\mathcal H_{0,n}}
 P_P(\delta_{P,n}^{\mathrm{GA}}>\varepsilon)\to0.
 \label{eq:T4-uniform-GA}
\end{equation}
Assume also that, almost surely,
\begin{equation}
 0<\underline\sigma^2
 \le\operatorname{Var}(G_{k,n}\mid\mathcal F_{0n})
 \le\overline\sigma^2<\infty
 \quad(1\le k\le2m_n),
 \label{eq:T4-Gaussian-variances}
\end{equation}
and that, for every $\varepsilon>0$,
\begin{equation}
 \sup_{P\in\mathcal H_{0,n}}
 P_P\{\mathcal E_n\cap
       (\sqrt{\ell_n}D_n>\varepsilon)\}\to0.
 \label{eq:T4-maximum-rate}
\end{equation}
Then
\begin{equation}
 \sup_{P\in\mathcal H_{0,n}}
 \sup_{t\in\mathbb R}
 \left|
 P_P(M_n^{\mathrm{ord}}\le t)
 -P_P(M_n^{\mathrm{or}}\le t)
 \right|\to0.
 \label{eq:T4-maximum-equivalence}
\end{equation}
\end{theorem}

\begin{proof}
We divide the proof into the deterministic geometric calculation,
the score transfer, and the maximum-law comparison.

\emph{Coefficient perturbation.}
Fix $\lambda$ and let $\mathcal J_\lambda$ be the complete jump set
of $\Psi_\lambda^S$.  Define
\[
 D_{\lambda,T_n}
 =\{x\in[0,1]^2:
      \Psi_\lambda^S(P_{T_n}x)\ne\Psi_\lambda^S(x)\}.
\]
If $x$ is at $\ell_\infty$-distance greater than $\Delta_n$ from
$\mathcal J_\lambda$, the segment joining $x$ to $P_{T_n}x$ does not
cross a jump of the piecewise-constant atom.  Thus
$D_{\lambda,T_n}$ lies in the $\Delta_n$-neighborhood of
$\mathcal J_\lambda$.  A union of a bounded number of line segments
of total length $O(h_\lambda)$ has a $\delta$-neighborhood of area at
most $C(h_\lambda\delta+\delta^2)$.  Since
$\Delta_n\le\overline\Delta_n$,
\begin{equation}
 |D_{\lambda,T_n}|
 \le C(h_\lambda\overline\Delta_n
        +\overline\Delta_n^2).
 \label{eq:T4-boundary-layer-area}
\end{equation}
Consequently,
\begin{align}
 \|\Psi_\lambda^S\circ P_{T_n}-\Psi_\lambda^S\|_1
 &\le2\|\Psi_\lambda^S\|_\infty
       |D_{\lambda,T_n}|\nonumber\\
 &\le C\left(
       \overline\Delta_n+
       \frac{\overline\Delta_n^2}{h_\lambda}
       \right).
 \label{eq:T4-Haar-L1}
\end{align}
On $\mathcal E_n$,
$\overline\Delta_n\le\overline r_n$; hence
\eqref{eq:T4-small-radius} makes the last display at most
$C\overline r_n$, uniformly in $\lambda$, for all sufficiently large
$n$.

Because $T_n$, and therefore $P_{T_n}$, is measure preserving, the
change of variables $z=P_{T_n}x$ gives
\[
 \theta_{\lambda,n}^{\mathrm{ord}}
 =\int g_{P,n}(x)
       \Psi_\lambda^S(P_{T_n}x)\,dx.
\]
It follows from \eqref{eq:T4-Haar-L1} that, on $\mathcal E_n$,
\[
 |b_{\lambda,n}^{\mathrm{coef}}|
 \le\|g_{P,n}\|_\infty
    \|\Psi_\lambda^S\circ P_{T_n}-\Psi_\lambda^S\|_1
 \le CB\overline r_n.
\]
Taking the maximum over $\lambda$ proves
\eqref{eq:T4-coefficient-bound}.  This is a boundary-layer argument
for the Haar atom; it uses no smoothness or bounded variation of
$q_{P,n}$.

Under \eqref{eq:T4-single-atom}, unit $L^2$ norm gives
\[
 \theta_{\lambda_*,n}^{\mathrm{or}}
 =s a_nh_*\|\Psi_{\lambda_*}^S\|_2^2
 =s a_nh_*.
\]
The coefficient bound gives the second assertion in
\eqref{eq:T4-signal-coefficient}.  Dividing its error by
$|a_n|h_*$ proves the claim under \eqref{eq:T4-beta-min}.

\emph{Detection, sign, and canonical localization.}
By the exact identity \eqref{eq:T4-score-identity} and the triangle
inequality,
\begin{align*}
 \max_\lambda|Z_{\lambda,n}^{\mathrm{ord}}
                   -Z_{\lambda,n}^{\mathrm{or}}|
 &\le\max_\lambda\sqrt{I_{\lambda,n}}
       |b_{\lambda,n}^{\mathrm{sc}}|
       +\max_\lambda|R_{\lambda,n}|\\
 &\le\max_\lambda\sqrt{I_{\lambda,n}}
       |b_{\lambda,n}^{\mathrm{coef}}|
       +\zeta_n+\max_\lambda|R_{\lambda,n}|=D_n,
\end{align*}
which proves \eqref{eq:T4-score-sup-bound}.

Intersect the oracle-margin event in \eqref{eq:T4-oracle-margin} with
the event in which
\[
 D_n+|\widehat c_n^{\mathrm{ord}}
        -\widehat c_n^{\mathrm{or}}|\le\Gamma_n.
\]
On this intersection,
\begin{align*}
 sZ_{\lambda_*,n}^{\mathrm{ord}}
 -\widehat c_n^{\mathrm{ord}}
 &\ge sZ_{\lambda_*,n}^{\mathrm{or}}
       -\widehat c_n^{\mathrm{or}}
       -D_n
       -|\widehat c_n^{\mathrm{ord}}
          -\widehat c_n^{\mathrm{or}}|\\
 &\ge2\Gamma_n>0.
\end{align*}
Because the critical value is nonnegative, this proves rejection and
correct sign.  Moreover,
\begin{align*}
 |Z_{\lambda_*,n}^{\mathrm{ord}}|
 -\max_{\lambda\ne\lambda_*}|Z_{\lambda,n}^{\mathrm{ord}}|
 &\ge |Z_{\lambda_*,n}^{\mathrm{or}}|
 -\max_{\lambda\ne\lambda_*}|Z_{\lambda,n}^{\mathrm{or}}|
 -2D_n\\
 &\ge\Gamma_n>0.
\end{align*}
Thus the fixed tie-breaking rule selects $\lambda_*$.  Equations
\eqref{eq:T4-coverage}, \eqref{eq:T4-oracle-margin}, and
\eqref{eq:T4-margin-rate} prove
\eqref{eq:T4-label-consistency}.  Since the canonical loss is zero
when $\widehat\lambda=\lambda_*$ and is at most one otherwise,
\[
 E_Pd_{\mathrm{loc}}(\widehat\lambda,\lambda_*)
 \le P_P(\mathcal E_n^c)
    +P_P(\widehat\lambda\ne\lambda_*;\mathcal E_n),
\]
which proves \eqref{eq:T4-canonical-localization}.

Under \eqref{eq:T4-information-critical-scaling}, the oracle
standardized signal and critical value have orders
\[
 |a_n|h_*\sqrt{I_{\lambda_*,n}}
 \asymp |a_n|h_*n\sqrt{\rho_n},
 \qquad
 \widehat c_n^{\mathrm{or}}\asymp\sqrt{\ell_n}.
\]
Their squared balance is exactly the rate in
\eqref{eq:T4-detection-boundary}.  On every fixed multiplicatively
separated alternative class, the assumed oracle gap is of order
$\sqrt{\ell_n}$, whereas \eqref{eq:T4-first-order-rate} is
$o_P(\sqrt{\ell_n})$.  Repeating the preceding rejection inequality
therefore transfers the oracle first-order sufficient boundary.

\emph{Maximum law.}
By \eqref{eq:T4-score-sup-bound},
\begin{equation}
 |M_n^{\mathrm{ord}}-M_n^{\mathrm{or}}|\le D_n.
 \label{eq:T4-max-Lipschitz}
\end{equation}
The uniform convergence in \eqref{eq:T4-maximum-rate} permits a
deterministic sequence $\gamma_n\downarrow0$ such that
\[
 \sup_{P\in\mathcal H_{0,n}}
 P_P\{\mathcal E_n\cap
       (\sqrt{\ell_n}D_n>\gamma_n)\}\to0.
\]
Set $\varepsilon_n=\gamma_n/\sqrt{\ell_n}$.  For every $t$,
\begin{align}
 &|P_P(M_n^{\mathrm{ord}}\le t)
       -P_P(M_n^{\mathrm{or}}\le t)|\nonumber\\
 &\qquad\le P_P(D_n>\varepsilon_n)
 +P_P(t-\varepsilon_n<M_n^{\mathrm{or}}
                    \le t+\varepsilon_n).
 \label{eq:T4-CDF-comparison}
\end{align}
The first term is $o(1)$ uniformly, because it is bounded by the
probability in the preceding display plus $\eta_n$.

For the second term, condition on $\mathcal F_{0n}$ and use
\eqref{eq:T4-GA-error}.  Gaussian anti-concentration for the maximum
of $2m_n$ centered Gaussian coordinates, together with
\eqref{eq:T4-Gaussian-variances}, gives a constant $C$ independent of
$P$ such that
\begin{align*}
 \sup_tP_P(t-\varepsilon_n<M_n^{\mathrm{or}}
                    \le t+\varepsilon_n)
 &\le C\varepsilon_n\sqrt{\ell_n}
      +2E_P\delta_{P,n}^{\mathrm{GA}}\\
 &=C\gamma_n+2E_P\delta_{P,n}^{\mathrm{GA}}.
\end{align*}
Since $0\le\delta_{P,n}^{\mathrm{GA}}\le1$,
\eqref{eq:T4-uniform-GA} implies
$\sup_PE_P\delta_{P,n}^{\mathrm{GA}}\to0$.  Taking the supremum in
\eqref{eq:T4-CDF-comparison} proves
\eqref{eq:T4-maximum-equivalence}.
\end{proof}

\begin{remark}[Boundary-strip sensitivity is not canonical localization]
\label{rem:T4-membership-sensitivity}
For a fixed label $\lambda$, the same boundary-layer calculation gives
\begin{equation}
 \frac{|P_{T_n}^{-1}S_\lambda\triangle S_\lambda|}
      {2|S_\lambda|}
 \le C\min\!\left\{
  1,
  \frac{\overline\Delta_n}{h_\lambda}
  +\left(\frac{\overline\Delta_n}{h_\lambda}\right)^2
 \right\}.
 \label{eq:T4-membership-sensitivity}
\end{equation}
This measures how the dyads represented by one fixed label change
between true- and reported-rank coordinates.  It is not the canonical
loss in \eqref{eq:T4-canonical-loss}.  In particular, exact selection
of $\lambda_*$ gives canonical loss zero, not an
$O_P(\Delta_n/h_*)$ remainder.
\end{remark}

\paragraph{A separate rank-space partial-identification experiment.}
Let $\mathfrak S_n$ denote the set of permutations of
$\{1,\ldots,n\}$.  For $r\ge0$, let
\begin{equation}
 \mathcal Q_n^{\mathrm{perm}}(r)
 =\left\{
 Q:\ Q\!\left(
   \widehat R\in\mathfrak S_n,\ 
   \|\widehat R-R\|_\infty\le nr
   \mid R
  \right)=1
 \text{ for every }R\in\mathfrak S_n
 \right\}.
 \label{eq:T4-reporting-kernels}
\end{equation}
The reporting-kernel label is not observed.  Hence this is an
uncertainty-set experiment, not a point-identified measurement model.

Let $x_q=(q-\tfrac12)/n$, let
$\alpha_n=\operatorname{logit}(\rho_n)$ with
$0<\rho_n\le1/2$, and let $\Lambda_n^S(h)$ be the aligned symmetric
Haar dictionary of width $h$.  Put
\[
 C_\Psi^\circ
 =\sup_{\lambda\in\Lambda_n^S(h)}
   h\|\Psi_\lambda^S\|_\infty<\infty.
\]
Given deterministic $M\in(0,B]$ and $\kappa_n^{\mathrm{str}}\to\infty$, define
\begin{equation}
 \mathcal A_n^{\mathrm{str}}(h,M,\kappa_n^{\mathrm{str}})
 =\left\{
 a>0:
 C_\Psi^\circ a\le M,\quad
 a^2n^2\rho_nh^2\ge\kappa_n^{\mathrm{str}}\ell_n
 \right\}.
 \label{eq:T4-strong-amplitudes}
\end{equation}
Assume this set is nonempty.  For
$\vartheta=(R,Q,\lambda,s,a)$ with
$R\in\mathfrak S_n$, $Q\in\mathcal Q_n^{\mathrm{perm}}(r)$,
$\lambda\in\Lambda_n^S(h)$, $s\in\{-1,1\}$, and
$a\in\mathcal A_n^{\mathrm{str}}(h,M,\kappa_n^{\mathrm{str}})$, define
\begin{align}
 p_{\vartheta,ij}
 &=\operatorname{logit}^{-1}\!\left\{
   \alpha_n+sah\,
   \Psi_\lambda^S(x_{R_i},x_{R_j})
  \right\},
 \qquad i<j,
 \label{eq:T4-lower-probabilities}\\
 P_\vartheta(dA,d\widehat r)
 &=Q(d\widehat r\mid R)
   \prod_{i<j}
   p_{\vartheta,ij}^{A_{ij}}
   (1-p_{\vartheta,ij})^{1-A_{ij}}.
 \label{eq:T4-lower-law}
\end{align}
Let $\Theta_n^{\mathrm{rank}}(r,h)$ be the indexed collection of all
such parameters $\vartheta$. The bound in \eqref{eq:T4-strong-amplitudes}
ensures that the logit perturbation is bounded by $M$ and hence that
$p_{\vartheta,ij}\asymp\rho_n$ uniformly.  Nonemptiness holds, for
example, when $n^2\rho_nh^2/\ell_n\to\infty$ and
$\kappa_n^{\mathrm{str}}$ is chosen
to diverge sufficiently slowly.

\begin{proposition}[Dyadic canonical-location obstruction]
\label{prop:T4-dyadic-obstruction}
Let $H=nh_*$ be an even integer.  Suppose that the width-$h_*$
dictionary contains three pairwise-disjoint aligned blocks
$I_0,I_1,J$, each consisting of $H$ consecutive ranks, with $I_1$
adjacent to $I_0$.  Let $\lambda_0=(I_0,J)$ and
$\lambda_1=(I_1,J)$, where each label includes the transposed
rectangle.  If $r_n\ge h_*$, then
\begin{equation}
 \inf_{\widetilde\lambda}
 \sup_{\vartheta\in\Theta_n^{\mathrm{rank}}(r_n,h_*)}
 E_{P_\vartheta}d_{\mathrm{loc}}
    \{\widetilde\lambda,\lambda(\vartheta)\}
 \ge\frac12.
 \label{eq:T4-rank-lower-bound}
\end{equation}
The infimum is over all measurable functions of
$(A,\widehat R,r_n)$ taking values in the decision space
$\mathcal D_n^{\mathrm{dec}}$ of
Theorem~\ref{thm:adaptive-detection-localization}, with loss one for
the cemetery decision.
\end{proposition}

\begin{proof}
Define the rank permutation
\begin{equation}
 \tau(q)=
 \begin{cases}
  q+H,&q\in I_0,\\
  q-H,&q\in I_1,\\
  q,&q\notin I_0\cup I_1.
 \end{cases}
 \label{eq:T4-block-swap}
\end{equation}
Because $H$ is even, $\tau$ carries the left and right halves of
$I_0$ to the corresponding halves of $I_1$, carries those of $I_1$
back to $I_0$, and fixes $J$.  Put
\[
 R_i^{(0)}=i,\qquad R_i^{(1)}=\tau(i).
\]
The aligned Haar factors therefore satisfy, for every labeled dyad,
\begin{equation}
 \Psi_{\lambda_1}^S
   (x_{R_i^{(1)}},x_{R_j^{(1)}})
 =\Psi_{\lambda_0}^S
   (x_{R_i^{(0)}},x_{R_j^{(0)}}).
 \label{eq:T4-atom-swap-equality}
\end{equation}
Choose any
$a_n\in\mathcal A_n^{\mathrm{str}}
(h_*,M,\kappa_n^{\mathrm{str}})$ and use the
same $s,a_n,$ and $\alpha_n$ at the two parameter points.  Equation
\eqref{eq:T4-atom-swap-equality} yields
$p_{ij}^{(0)}=p_{ij}^{(1)}$ for every $i<j$.

Let $R^{\mathrm{id}}=(1,\ldots,n)$ and
$R^\tau=(\tau(1),\ldots,\tau(n))$, and define the single reporting
kernel
\begin{equation}
 Q^\star(d\widehat r\mid R)=
 \begin{cases}
  \delta_{R^{\mathrm{id}}}(d\widehat r),
     &R\in\{R^{\mathrm{id}},R^\tau\},\\
  \delta_R(d\widehat r),&\text{otherwise}.
 \end{cases}
 \label{eq:T4-common-kernel}
\end{equation}
Since
$\|R^\tau-R^{\mathrm{id}}\|_\infty=H=nh_*$,
$Q^\star\in\mathcal Q_n^{\mathrm{perm}}(r_n)$ whenever
$r_n\ge h_*$.  Under both parameter points,
$\widehat R=R^{\mathrm{id}}$ almost surely.  Together with the
equality of the Bernoulli probability matrices, this proves
\begin{equation}
 \mathcal L_0(A,\widehat R,r_n)
 =\mathcal L_1(A,\widehat R,r_n).
 \label{eq:T4-common-observable-law}
\end{equation}

Writing $S_k=S_{\lambda_k}$, pairwise disjointness of
$I_0,I_1,J$ gives
\[
 |S_0|=|S_1|=2h_*^2,\qquad S_0\cap S_1=\varnothing,
 \qquad d_{\mathrm{loc}}(\lambda_0,\lambda_1)=1.
\]
For every possible output $\widetilde\lambda$,
\begin{equation}
 d_{\mathrm{loc}}(\widetilde\lambda,\lambda_0)
 +d_{\mathrm{loc}}(\widetilde\lambda,\lambda_1)\ge1.
 \label{eq:T4-loss-triangle}
\end{equation}
Indeed, if there is no output or its width differs from $h_*$, the
relevant losses equal one.  Otherwise put $A=|S_0|=|S_1|$ and
$a_k=|S_{\widetilde\lambda}\triangle S_k|/(2A)$.  The
symmetric-difference triangle inequality gives
\[
 |S_{\widetilde\lambda}\triangle S_0|
 +|S_{\widetilde\lambda}\triangle S_1|
 \ge|S_0\triangle S_1|.
\]
Because $S_0$ and $S_1$ are disjoint, this says $a_0+a_1\ge1$;
hence $\min(1,a_0)+\min(1,a_1)\ge1$, even when clipping is active and
without restricting the area of the reported support.
Taking expectations in \eqref{eq:T4-loss-triangle} under the common
law \eqref{eq:T4-common-observable-law} shows that at least one of
the two risks is at least $1/2$.  Taking the infimum over estimators
proves \eqref{eq:T4-rank-lower-bound}.  Because
$a_n^2n^2\rho_nh_*^2/\ell_n
\ge\kappa_n^{\mathrm{str}}\to\infty$, this obstruction
persists for oracle-supercritical signals.
\end{proof}

\begin{remark}[Why there is no linear all-radius lower bound here]
The preceding proposition uses permutation-valued reports
and the explicit block-swap construction has normalized radius $h_*$.  A
coordinatewise midpoint report would reduce the radius to $h_*/2$, but
it gives identical reported values to swapped vertex pairs and is not
a permutation.
Moreover, the aligned dyadic dictionary has no atom at an arbitrary
fractional displacement below $h_*$.  Thus the construction proves a
saturated obstruction for $r_n\ge h_*$, not a bound proportional to
$\min\{1,r_n/h_*\}$ for every $r_n$.  Such a linear modulus requires
a separately defined continuous-coordinate experiment with
continuously translated atoms.
\end{remark}

\section{All-radius localization floor}\label{sec:supp-floor}


\subsection{All-radius floors for canonical coordinate localization}
\label{supp:allradius}

\paragraph{The continuous-coordinate experiment.}
Let $n\ge2$ and $r_n\ge0$.  A configuration is
$P=(g,f_0,a,h,x,y)$: coordinates $U_1,\dots,U_n$ are i.i.d.\ with a density $g$
on $[0,1]$ satisfying $\tfrac34\le g\le\tfrac32$; given $U$, edges are
independent Bernoulli with
\[
\operatorname{logit}p_{ij}=\alpha_n+f_0(U_i,U_j)+a\,h\,\Psi^S_{h,x,y}(U_i,U_j),
\qquad
\Psi^S_{h,x,y}(u,v)=\tfrac{1}{\sqrt2}\{\phi_h(u-x)\phi_h(v-y)+\phi_h(u-y)\phi_h(v-x)\},
\]
where $\phi_h=h^{-1/2}\{\mathbf 1_{[0,h/2)}-\mathbf 1_{[h/2,h)}\}$, $f_0\in
V_{J_0}$ is symmetric with $\|f_0\|_\infty\le B$ and $\int f_0=0$ (the mean-zero
normalization that makes the $(\alpha_n,f_0)$ decomposition unique), $0<h\le2^{-J_0-1}$, and
$2\sqrt2\,|a|\le B$, so that all logits lie in
$[\alpha_n-2B,\alpha_n+2B]$ and all edge probabilities in
$[c_B\rho_n,C_B\rho_n]$.  The single display above defines the departure for
every location pair, including $x=y$, where it evaluates to
$\sqrt2\,\phi_h(u-x)\phi_h(v-x)$; no separate diagonal convention is used, and
we do not normalize $\Psi^S$ in $L^2$ (for $0<|x-y|<h$ its squared norm is
$1+\langle\phi_h(\cdot-x),\phi_h(\cdot-y)\rangle^2$), since only the
sup-norm bound $\|h\Psi^S_{h,x,y}\|_\infty\le2\sqrt2$ enters the class
constraint.  Because $\Psi^S$ is symmetric in $(x,y)$, the location is
identified only as an unordered pair; the \emph{canonical location} is the
quotient point $\{x,y\}$ with $x\le y$, and the loss between quotient points is
\[
d\bigl(\{x,y\},\{x',y'\}\bigr)
=\min_{\pi\in\{\mathrm{id},\mathrm{swap}\}}
\max\bigl(|x-\pi(x')|,\,|y-\pi(y')|\bigr).
\]
A reporting kernel is any conditional law of $\widehat U=(\widehat U_i)$ given
$U$ with $\max_i|\widehat U_i-U_i|\le r_n$ almost surely; beyond the radius,
the kernel is unrestricted and unknown to the analyst (radius-only partial
information; observing only the ranks of $\widehat U$ is a coarsening and is
covered).  The analyst observes $(A,\widehat U,r_n)$.

\begin{theorem}[All-radius floors]\label{thm:S-allradius}
\textup{(a)} \emph{Order-limited floor.}  Fix $t$ with
$0\le t\le\min\{r_n,\ (2^{-J_0}-h)/13\}$ and a configuration
$P_0=(\mathrm{Unif},f_0,a,h,x,y)$ whose blocks $[x,x+h]$ and
$[y,y+h]$ each lie inside a coarse cell at distance at least $7t$ from that
cell's endpoints, and which either coincide ($x=y$) or are at distance at least
$h+13t$ from each other.  Then there is a configuration
$P_1=(g_1,f_0,a,h,x+t,y+t)$ in the class and reporting kernels of radius $t$
under $P_0$ and $P_1$ such that the joint laws of $(A,\widehat U)$ coincide.
Consequently, for every $n$, every amplitude $a$, and every estimator
$(\widehat x,\widehat y)$ of the quotient point,
\[
\max_{P\in\{P_0,P_1\}}P\bigl(d\bigl(\{\widehat x,\widehat y\},\{x_P,y_P\}\bigr)\ge t/2\bigr)\ \ge\ \tfrac12,
\]
and every confidence region $\mathcal C$ for the canonical location with
$\inf_PP\{\{x_P,y_P\}\in\mathcal C\}\ge1-\alpha$ and $\alpha<\tfrac12$
satisfies $\max_PE_P\operatorname{diam}_d\mathcal C\ge(1-2\alpha)\,t$.

\textup{(b)} \emph{Signal-limited floor.}  Let $r_n=0$, the design and $f_0$
known.  For $0<\varepsilon\le h$ and configurations $P_0$ at $\{x,y\}$ and
$P_1$ at $\{x+\varepsilon,y\}$ (both translated blocks inside $[0,1]$; $x=y$
permitted, in which case $P_1$ sits at the overlapping pair
$\{x+\varepsilon,x\}$),
\[
\operatorname{KL}(P_0\|P_1)\ \le\ C_B\,a^2n^2\rho_n\,h\,\varepsilon .
\]
Hence with $\varepsilon_n=\min\{h,\ 1/(C_B'a^2n^2\rho_nh)\}$ and $C_B'=2C_B$,
one of the two configurations gives every estimator error
$d\ge\varepsilon_n/4$ with probability at least $\tfrac14$.  The threshold is
$\varepsilon_n/4$ rather than $\varepsilon_n/2$: the two configurations are
$\{x,y\}$ and $\{x+\varepsilon_n,y\}$, and when $y=x+\varepsilon_n/2$ their
quotient distance is $\varepsilon_n/2$, so the point $\{x+\varepsilon_n/4,\,
x+3\varepsilon_n/4\}$ is within $\varepsilon_n/4$ of both and an estimator
returning it errs by $\varepsilon_n/4$ under either; the two-point argument
therefore guarantees error at least half the quotient distance, which is
$\varepsilon_n/4$ in the worst case.  Under a no-crossing condition
$|y-x|\ge\varepsilon_n$ the distance is $\varepsilon_n$ and the threshold
returns to $\varepsilon_n/2$.
Moreover the floor transfers to every radius: for $r_n>0$ the identity
reporting kernel $\widehat U=U$ is admissible, under it the analyst's data
coincide with the $r_n=0$ data, and any estimator in the positive-radius
experiment must in particular perform against that kernel, so the same
two-point bound holds at every $r_n\ge0$ and the joint floor
$\max\{\text{(a)},\text{(b)}\}$ is attained within one experiment.
\end{theorem}

\begin{proof}
(a) Let $L=3t$ and, for a block $[z,z+h]$ inside the coarse cell
$[c_0,c_1]$, define on $W_z=[z-2L,\ z+h+t+2L]\subset[c_0,c_1]$ the map
$\varphi_z$ that equals the identity outside $W_z$, equals $u\mapsto u+t$ on
$[z-L,\ z+h+t+L]$, and is affine on $[z-2L,z-L]$ (slope $\tfrac43$) and on
$[z+h+t+L,\ z+h+t+2L]$ (slope $\tfrac23$).  Then $\varphi_z$ is a
continuous increasing bijection of $[0,1]$, $|\varphi_z(u)-u|\le t$ for all
$u$, $\varphi_z$ maps every coarse cell onto itself, and its slope lies in
$\{\tfrac23,1,\tfrac43\}$.  Put $\varphi=\varphi_x\circ\varphi_y$ (a single
$\varphi_x$ if $x=y$); the windows $W_x,W_y$ are disjoint by the separation
assumption, so $\varphi$ inherits the three properties.  Let $g_1$ be the
density of $\varphi(U)$ for $U\sim\mathrm{Unif}$;
$g_1=1/\varphi'\circ\varphi^{-1}\in[\tfrac34,\tfrac32]$, so $P_1$ is in the
class.  (The window $W_z$ extends $6t$ left of $z$ and $6t+t$ right of $z+h$,
which the $7t$ room accommodates, and two windows fit disjointly once the
blocks are $h+13t$ apart; this is also why $t$ is capped at
$(2^{-J_0}-h)/13$.)

Coupling.  Draw $U$ i.i.d.\ uniform.  Under $P_0$ the true coordinates are
$U$; under $P_1$ they are $U'=\varphi(U)$.  For every $(u,v)$,
\[
f_0(\varphi(u),\varphi(v))=f_0(u,v),
\qquad
\Psi^S_{h,x+t,y+t}(\varphi(u),\varphi(v))=\Psi^S_{h,x,y}(u,v).
\]
The first identity holds because $f_0\in V_{J_0}$ is constant on products of
coarse cells and $\varphi$ maps each coarse cell onto itself.  For the second,
$\phi_h(\varphi(u)-x-t)\ne0$ requires $\varphi(u)\in[x+t,x+t+h)$, whose
$\varphi$-preimage is exactly $[x,x+h)$ where $\varphi(u)=u+t$; hence
$\phi_h(\varphi(u)-x-t)=\phi_h(u-x)$ for all $u$, and likewise for $y$.
Therefore the conditional edge law given the seed $U$ is the same under $P_0$
and $P_1$: generate the same $A$.  Report $\widehat U=U$ in both worlds.
Under $P_0$ the radius is $0\le t$; under $P_1$,
$|\widehat U_i-U'_i|=|U_i-\varphi(U_i)|\le t$.  Both kernels are valid, and
$(A,\widehat U)$ has the same joint law under $P_0$ and $P_1$.

Consequences.  The canonical quotient points $\{x,y\}$ and $\{x+t,y+t\}$ are at
$d$-distance exactly $t$: the identity matching gives $t$ and, under the
separation $|x-y|\ge h+13t$ or $x=y$, the swap matching gives at least
$|x-y|-t\ge t$ as well.  For any estimator, the events
$\{d(\{\widehat x,\widehat y\},\{x,y\})<t/2\}$ and
$\{d(\{\widehat x,\widehat y\},\{x+t,y+t\})<t/2\}$ are disjoint by the
triangle inequality (the quotient loss is a metric), and since the data have
the same law under $P_0$ and $P_1$ their probabilities, computed under the
common law, sum to at most one; so one error probability is at least
$\tfrac12$.  For a confidence region, $\mathcal C$ contains both quotient
points with probability at least $1-2\alpha$ under the common law, on which
event $\operatorname{diam}_d\mathcal C\ge t$.

(b) With $r_n=0$ the analyst observes $(A,U)$, and $U$ has the same law under
$P_0$ and $P_1$, so $\operatorname{KL}(P_0\|P_1)=E\sum_{i<j}
\operatorname{KL}\{\mathrm{Bern}(p^0_{ij})\|\mathrm{Bern}(p^1_{ij})\}$.
Because both configurations use the single displayed formula for $\Psi^S$, the
two logits differ only when $(U_i,U_j)$ or $(U_j,U_i)$ lies in
$\{[x,x+\varepsilon)\cup[x+\tfrac h2,x+\tfrac h2+\varepsilon)\cup[x+h,x+h+\varepsilon)\}
\times([y,y+h]\cup[x,x+h])$, a set of Lebesgue measure at most
$12\varepsilon h$ counting both orientations (the second factor covers the
overlapping and diagonal cases; for separated blocks the sharper $6\varepsilon
h$ holds), and there the logit gap is at most
$2\|a h\Psi^S\|_\infty\le4\sqrt2|a|$.  For Bernoulli laws with logits in
$[\alpha_n-2B,\alpha_n+2B]$,
$\operatorname{KL}\le(p^0-p^1)^2/\{p^1(1-p^1)\}\le C_B\rho_na^2$ by the mean
value theorem and $p^1(1-p^1)\ge c_B\rho_n$.  Summing over the at most
$\binom n2$ dyads and the design expectation gives
$\operatorname{KL}\le C_Ba^2n^2\rho_nh\varepsilon$ after renaming $C_B$.  With
$\varepsilon_n$ as stated, $\operatorname{KL}\le\tfrac12$, so
$\|P_0-P_1\|_{\rm TV}\le\tfrac12$, and Le Cam's two-point bound gives the
claim.  The transfer to $r_n>0$ is as stated: any estimator
$T(A,\widehat U)$ in the positive-radius experiment, evaluated against the
identity kernel, is an estimator in the $r_n=0$ experiment, and the minimax
risk over kernels is at least the risk at that kernel.
\end{proof}

\begin{remark}[Scope of the floors]\label{rem:S-allradius}
Part (a) is a pure identification statement: it holds for every amplitude,
every $n$, and every reporting mechanism of radius $r_n$, with the translation
capped at $t\le\min\{r_n,(2^{-J_0}-h)/13\}$; under the standing class
constraint $h\le2^{-J_0-1}$ the cap is at least $\min\{r_n,2^{-J_0}/26\}$, so
the floor is of order $\min\{r_n,2^{-J_0}\}$ with an explicit constant, and it
degenerates as $h\uparrow2^{-J_0}$, exactly the regime the class excludes.
The lower bounds live in the radius-only partial-information experiment
defined above (the reporting kernel is adversarial within the radius), and
that is also the experiment in which the paper's estimation theory is stated;
a known-kernel, point-identified measurement model is strictly more
informative and is not covered by these floors.  With a fixed known design and
the aligned single-atom class the coupling stops at $t=h/2$, and for that
class and the exact-atom loss the regime $0<r_n<h$ remains open between
Theorem~\ref{thm:ordering-error-transfer} and
Proposition~\ref{prop:T4-dyadic-obstruction}.  Part (b) is the classical
change-point floor in the coordinate metric; the aligned scan localizes only
to the atom width $h$, and whether a finer scan attains the signal-limited
rate is not addressed here.
\end{remark}

\end{document}